\documentclass[twoside,11pt]{article}
\usepackage{jmlr2e}

\usepackage{amsmath, amsfonts, amssymb, amsbsy, bigstrut, graphicx, enumerate,  upref, longtable, comment, comment, csquotes, booktabs, array, bbm}
\usepackage{mathtools}
\usepackage{booktabs}
\usepackage{tabularx}
\usepackage{tikz}
\usepackage{float}
\usepackage{subcaption}
\usepackage{pstricks}

\newcommand{\eqd}{\stackrel{d}{=}}
\newcommand{\RR}{\mathbb{R}}

\newcommand{\PP}{{\mathbb{P}}}
\newcommand{\DD}{{\mathbb{D}}}
\newcommand{\EE}{{\mathbb{E}}}
\newcommand{\one}{\mathbbm{1}}
\newcommand{\var}{\mathrm{Var}}

\DeclareMathOperator*{\argmin}{argmin} 

\numberwithin{equation}{section}

\usepackage{lastpage}
\jmlrheading{}{}{1-\pageref{LastPage}}{09/23}{}{}{Mona Azadkia, and Fadoua Balabdaoui}

\firstpageno{1}

\begin{document}
\title{Linear Regression With Unmatched Data: A Deconvolution Perspective}
\author{\name Mona Azadkia \email m.azadkia@lse.ac.uk \\
       \addr Department of Statistics\\
       London School of Economics\\
       London, United Kingdom
       \AND
       \name Fadoua Balabdaoui \email fadoua.balabdaoui@stat.math.ethz.ch \\
       \addr Department of Mathematics\\
       ETH Z\"{u}rich \\
       Z\"{u}rich, Switzerland}


\maketitle

\begin{abstract}%
Consider the regression problem where the response $Y\in\RR$ and the covariate $X\in\RR^d$ for $d\geq 1$ are \textit{unmatched}.  Under this scenario, we do not have access to pairs of observations from the distribution of $(X, Y)$, but instead, we have separate data sets $\{Y_i\}_{i=1}^{n_Y}$ and $\{X_j\}_{j=1}^{n_X}$, possibly collected from different sources. We study this problem assuming that the regression function is linear and the noise distribution is known, an assumption that we relax in the applications we consider. We introduce an estimator of the regression vector based on deconvolution and demonstrate its consistency and asymptotic normality under identifiability. Even when identifiability does not hold, we show in some cases that our estimator, the DLSE (Deconvolution Least Squared Estimator), is consistent in terms of an extended $\ell_2$ norm. Using this observation, we devise a method for semi-supervised learning, i.e., when we have access to a small sample of matched pairs $\{(X_k, Y_k)\}_{k=1}^m$. Several applications with synthetic and real data sets are considered to illustrate the theory.
\end{abstract}


\section{Introduction}
Consider the standard regression setting
\begin{eqnarray}\label{matchedData}
Y = m_0(X) + \epsilon
\end{eqnarray}
where $\epsilon\in\RR$ is the noise variable, $X\in\RR^d$ is the vector of features and $m_0:\RR^d\rightarrow\RR$ is a measurable function. Given an i.i.d. sample of $(Y_i, X_i)$, the problem of estimating $m_0$ has been vastly studied in the literature. The estimation problem is much harder when we do not have access to \textit{matched}/\textit{linked} data, i.e. the pairs $(Y_i, X_i)$ from~\eqref{matchedData}, but instead we have separate samples $\{Y_i\}_{i=1}^{n_Y}$ and $\{X_j\}_{j=1}^{n_X}$ such that $Y_i$ and $X_j$'s have the same distribution as $Y$ and $X$. Under this setting, the exact equality in (\ref{matchedData}) is replaced by the following equality in distribution
\begin{eqnarray}\label{unlinkedData}
Y \stackrel{d}{=} m_0(X) + \epsilon.
\end{eqnarray}
This type of data commonly arises in applications when the data has been collected through different sources. It can also occur when the link between the covariates and the response has been deleted because of privacy concerns. This problem is known as \textit{unmatched} or \textit{unlinked} regression in the literature and has been studied in different variations. 

Unmatched regression can be seen as a generalization of a problem known as \textit{permuted} regression or \textit{unlabeled sensing}. In permuted regression, one again observes only the variables $\{Y_i\}_{i=1}^n$ and $\{X_i\}_{i=1}^n$ separately. The difference with unmatched regression is that in permuted regression, it is assumed that the link is lost because of the existence of some permutation $\pi^*$ such that 
\begin{eqnarray}\label{permutedRegression}
Y_i = m_0(X_{\pi^*(i)}) + \epsilon_i.
\end{eqnarray}
Clearly, permuted regression is a special case of unmatched regression. However, note that the equality \eqref{permutedRegression} is stronger than the equality in the distribution in \eqref{unlinkedData} for unmatched regression since, in unmatched regression, such a permutation may not exist. Permuted regression has received considerable attention, especially when $m_0$ is a linear transformation, e.g., \cite{Unnikrishnan:2018gp, Pananjady:2018hd, abid2017linear, hsu2017linear, Slawski19, Tsakiris20, slawski2020two, Tsakiris19, zhang2021benefits, slawski2021pseudo, Rigollet2018}. In permuted regression, finding the corresponding permutation has also been the focus of many research works. Note that this is a much harder problem than just estimating $m_0$. If one estimates $\pi^*$, by some estimator $\hat{\pi}$ say, then estimating $m_0$ boils down to a simple estimation with matched data $(X_{\hat{\pi}(i)}, Y_i)$. For a detailed summary of the related work on permuted regression, we refer the reader to \cite{slawski2022permuted}.

The problem of unmatched regression is of special interest in areas such as microeconomics, where the variable of interest has not been observed jointly with some of the covariates. This problem is also known as data fusion, where multiple observational and experimental data sets exist and the links are unavailable. Popular approaches in this situation include methods based on \textit{matching}, see e.g.~\cite{cohen2002learning, monge1996field, walter1999matching}. Matching-based methods rely on the access to an extra \textit{contexual} variable such that they can use this variable to pair the response variable and the covariates. 

In the absence of such contextual variables, there is little hope for pairing the variables. In a recent work~\cite{carp16}, estimating $m_0$ in unmatched regression was studied for $d=1$, assuming that the distribution of $\epsilon$ is known and $m_0$ is monotone without assuming any contextual variable. Note that when $m_0$ is monotone, this problem is also known as unmatched \textit{isotonic} regression. In~\cite{carp16}, authors have shown the close connection between unmatched isotonic regression and deconvolution. In~\cite{carp16} for estimating $m_0$, authors have resorted to kernel deconvolution as the main estimation technique. Under some smoothness assumptions, the authors have provided the rate of convergence for the proposed kernel estimator obtained with available deconvolution methods in the literature. In~\cite{Balabdaoui21}, the authors made use of the tight relationship between the unmatched isotonic regression and deconvolution and provided an estimator for $m_0$ under the assumption that the noise distribution is known. Their method follows the idea of estimation of the mixing distribution for a normal mean, as done by~\cite{edelman}. The authors provide a rate of convergence of a weighted $\ell_1$-distance of their estimator of $m_0$ under several smoothness regimes of the distribution of the noise as well as some regularity conditions satisfied by $m_0$.  For example, in the case where the noise distribution is ordinary smooth with a smoothness parameter equal to $\beta$, it follows from their Theorem 1 that the convergence rate of their monotone estimator to the truth when both are restricted on certain compacts cannot be slower than $n^{-1/(2(2\beta+1))}$.  In~\cite{Rigollet2018}, a minimum Wasserstein deconvolution estimator was suggested that achieves the rate of $O(\log(\log n) / \log n)$ and it was shown
that for normally distributed errors, this rate is optimal. In \cite{MM} the minimum Wasserstein deconvolution estimator achieves much better rate of $O(n^{-1/(2p)})$ for $\ell_p$ risks for discrete noise which is optimal. Later in~\cite{slawski2022permuted}, the problem was considered for $d\geq 1$ where the authors have shown that a generalized notion of monotonicity, called cyclical monotonicity, of the regression function, is sufficient for estimation of the regression function $m_0$. Their method leverages ideas from the theory of optimal transport, specifically the Kiefer-Wolfowitz nonparametric maximum likelihood estimator. \cite{slawski2022permuted} provide the rate of convergence of their estimator in terms of $\ell_2$-distance under some smoothness assumption on $m_0$ and assuming that the noise is Gaussian.

\textbf{Our Contribution:} In this work, we consider the problem of estimating the regression vector in unmatched linear regression under the assumption that the noise distribution is known. Our work conceptually follows~\cite{Balabdaoui21}, although there is a considerable difference between the linear model and the monotone one. Compared to~\cite{Balabdaoui21, carp16}, the main contribution of our work is that our theory is valid for any given dimension $d\geq 1$ (though not depending on the sample size of the observations). This framework, to the best of our knowledge, has been only considered in~\cite{slawski2022permuted}. Our proposed deconvolution least squares estimator (DLSE) is not necessarily consistent; therefore, we cannot directly compare our results to those obtained in \cite{slawski2022permuted}. On the other hand, when we are in a setting in which our estimator is consistent, we do have the faster rate of convergence $O(n^{-1/2})$ compared to the rate of convergence $O(\log(\log n)/\log n)$ in ~\cite{slawski2022permuted}, without having to assume that noise is Gaussian. Such a fast rate may come as a big surprise given the non-standard situation of lacking any knowledge of the existing link, even partial, between the responses and covariates. The explanation is that, under the identifiability of the model, the estimation problem can be cast in the usual scope of the theory of $M$-estimators in parametric models. In the settings where the DLSE is not consistent, we show that a generalized notion of the estimator's norm is consistent, which can be used in semi-supervised scenarios. 

\textbf{Outline of the paper:} In Section~\ref{Methodology} we introduce our methodology for studying the problem. In Section~\ref{MainResult}, we provide the main results of the paper by first introducing the estimator and its properties and then applying it in a semi-supervised setting. Section~\ref{AppSynRealData} provides simulations on synthetic and real data. We provide the proofs of all lemmas and theorems in the supplementary material. 

\section{Methodology}\label{Methodology}
Let $X$ be a random vector in $\RR^d$ for $d\geq 1$. Consider $Y$ to be a random variable such that $Y\eqd\beta^T_0X+\epsilon$ where $\epsilon$ is the noise random variable, and $\beta_0\in\RR^d$ is a deterministic vector of coefficients. In this work we assume that we have access to two independent data sets $\{Y_i\}_{i=1}^{n_Y}$ and $\{X_i\}_{i=1}^{n_X}$ such that $Y_i, i=1, \ldots, n_Y$ and $X_i, i=1, \dots, n_X$ are i.i.d random variables which are distributed as $Y$ and $X$ respectively. 

For a random variable or vector $Z$, we denote by $F^Z$ its cumulative distribution function, that is $F^Z(z) = \PP(Z \le z), \ z \in \mathbb R$. In the case of a vector, the inequality should be considered component-wise. The equality in distribution $Y\eqd\beta^T_0 X  + \epsilon$ yields
\begin{eqnarray}\label{TheModel}
F^Y = F^{\beta^T_0 X} \star F^\epsilon
\end{eqnarray}
where $\star$ denotes the convolution operator. Note that $\beta_0$ minimizes the function
\begin{eqnarray}\label{D_population}
\beta\mapsto\mathcal{D}_{p, \mathbb{Q}}(\beta):= \int\Big\vert F^Y(y) - \big[F^{\beta^T_0 X} \star F^\epsilon\big](y)\Big\vert^p d\mathbb{Q}(y)
\end{eqnarray}
where $p \ge 1$ is an integer and $\mathbb{Q}$ is a positive measure on $\RR$. In this work, we focus on the case where $X$ has a continuous distribution. Also, we assume that $F^\epsilon$ is known. However, we do not assume that $\epsilon$ is necessarily Gaussian. Note that the assumption that the noise distribution is known was made in all the prior works on unmatched regression. In the applications below, we relax this assumption by letting the scale parameter of the distribution of the noise unspecified. Deriving the theory for such a relaxation is still needed as the estimation procedure needs to be extended to allow for the additional estimator of the scale parameter.

A natural choice of $\mathbb{Q}$ is $F^Y$, the probability measure induced by the distribution of $Y$ on $\RR$. Given a sample $\{Y_1, \cdots, Y_{n_Y}\}$ we can use the empirical estimate of $F^Y$
\begin{eqnarray*}
F^Y_{n_Y}(y) := {n_Y}^{-1} \sum_{i=1}^{n_Y} \one_{(-\infty, y]}(Y_i).
\end{eqnarray*}
By using $F_{n_Y}^Y$ as a surrogate of $\mathbb{Q}$, we get the empirical version of \eqref{D_population}
\begin{eqnarray}\label{Dnp}
\DD_{n_Y, n_X,p}(\beta) := \int\Big\vert F^Y_{n_Y}(y)-{n_X}^{-1}\sum_{i=1}^{n_X} F^\epsilon(y - \beta^T X_i)\Big\vert^p d F^Y_{n_Y}(y).
\end{eqnarray}
Without loss of generality and for the sake of a simpler notation, we assume that $n_Y = n_X = n$. Also, we let $\DD_{n,p}(\beta) := \DD_{n_Y, n_X,p}(\beta)$. The main idea of this paper is to take the solution of the optimization problem $\min_{\beta \in \RR^d}\DD_{n,p}(\beta)$ as an estimate of $\beta_0$. In Section \ref{MainResult}, we study this optimization problem in detail.

\section{Main Result}\label{MainResult}
\noindent \textbf{Notation and Definition.}
Let $C^k$ be the class of all functions on $\RR$ such that their first $k$ derivatives exist and are continuous. Let $\mathcal{B}_0$ be the set of all $\beta\in\RR^d$ such that the unlinked linear regression model in \eqref{TheModel} holds, i.e. $\mathcal{B}_0 =\{\beta^T X \stackrel{d}{=} \beta_0^T X\}$. Also, let us write $\mathcal{B}_{n, p} := \{\beta\in\RR^d\text{ s.t. }\DD_{n, p}(\beta) = \min\limits_{\beta}\DD_{n, p}(\beta)\}$. For $\beta\in\RR^d$, the Euclidean norm of $\beta$ is defined as $\Vert \beta \Vert_2 = \sqrt{\sum_{i=1}^d \beta^2_i}$. Let $S^{d-1} = \{x\in\RR^d:\|x\|_2 = 1\}$ be the $d$-sphere. For $\Sigma$ a positive-definite matrix, $\Vert x\Vert_{2, \Sigma}:=\sqrt{x^T\Sigma x}$ for $x\in\RR^d$ is a norm. 

\subsection{The Deconvolution Least Squares Estimator}
Given samples $\{Y_i\}_{i=1}^n$ and $\{X_i\}_{i=1}^n$ let $\hat{\beta}_n := \argmin_{\beta\in\RR^d} \DD_{n,p}(\beta)$. In the following proposition, we show that when $\epsilon$ admits a continuous distribution, there exists at least a $\hat{\beta}_n\in\RR^d$ such that it minimizes $\DD_{n,p}(\beta)$.

\begin{proposition}\label{Existence}
Assume that $F^\epsilon$ is continuous on $\mathbb R$. For a given integer $p \ge 1$, and for $n \ge 8d$ with probability $1$ we have $|\mathcal{B}_{n, p}|\geq 1$.
\end{proposition}

Note that $\mathcal{B}_0$ may have more than one element. For example for $X = [X^1, \cdots, X^d]^T\in\RR^d$ where $X^i$ are exchangeable random variables it is easy to see that we have $\beta_0^\pi:= [\beta_0^{\pi(1)}, \cdots, \beta_0^{\pi(d)}]^T \in \mathcal{B}_0$ where $\pi$ is any permutation on $[d]$. The other interesting case is when $X = [X^1, \cdots, X^d]^T\in\RR^d$ is distributed as $\mathcal{N}(0, \Sigma)$. In Lemma~\ref{Charbeta} we show that in this case $\mathcal{B}_0$ can be characterized in terms of $\Vert\cdot\Vert_{2, \Sigma}$.

\begin{lemma}\label{Charbeta}
Suppose the unlinked linear regression model in \eqref{TheModel} holds with $X\sim\mathcal{N}(0, \Sigma)$ where $\Sigma$ is non-singular. Then, there exists a constant $c > 0$ such that
\begin{eqnarray*}
\mathcal{B}_0 = \left\{\beta \in \RR^d: \Vert\beta\Vert_{2, \Sigma} = c \right\}.
\end{eqnarray*}
\end{lemma}
Since $\mathcal{B}_0$ is not necessarily a singleton, it is impossible to have a consistency result in the classical sense. On the other hand, with the characterization of $\mathcal{B}_0$ via the norm constraint, Theorem \ref{convNorm} shows that $\hat{\beta}_n$ is consistent in terms of such a norm constraint. From now on, we will take $p=2$. This means that our deconvolution least squares estimator is any vector $\hat{\beta}_n$ which minimizes the criterion \eqref{Dnp} for $p = 2$.
\begin{theorem}\label{convNorm}
Suppose that $F^\epsilon$ is continuous on $\mathbb R$. If there exists a non-singular positive-definite $d \times d$ matrix $\Gamma$ such that $\Vert\beta\Vert_{2, \Gamma} = c$ for all $\beta \in \mathcal{B}_0$ and some $c\in [0, \infty)$, then if $\Vert \hat{\beta}_n\Vert_2 = O_{\PP}(1)$ we have
\begin{eqnarray*}
\Vert \hat{\beta}_n \Vert_{2, \Gamma} \stackrel{\PP}{\rightarrow} c.
\end{eqnarray*}
\end{theorem}
Note that the result in Theorem~\ref{convNorm} requires that $\Vert\hat{\beta}_n\Vert_2 = O_{\PP}(1)$. The following proposition shows that when $X \sim \mathcal{N}(0, \Sigma)$ this requirement is met, and therefore we have the consistency of the DLSE $\hat{\beta}_n$ in terms of $\Vert\cdot\Vert_{2, \Sigma}$.
\begin{proposition}\label{Gaussian}
Suppose that $X \sim \mathcal{N}(0, \Sigma)$ for some non-singular $\Sigma$. Then, $\Vert \hat{\beta}_n \Vert_2 = O_{\PP}(1)$ and 
\begin{eqnarray*}
\Vert\hat{\beta}_n\Vert_{2, \Sigma}\stackrel{\PP}{\rightarrow} c, \    \  \textrm{and}  \  \   \Vert\hat{\beta}_n\Vert_{2, \widehat{\Sigma}_n}\stackrel{\PP}{\rightarrow} c,
\end{eqnarray*} 
where $c = \sqrt{\beta^T \Sigma \beta}$ for any $\beta \in\mathcal{B}_0$ for which the unlinked linear regression model in (\ref{TheModel}) holds, and $\widehat{\Sigma}_n$ is any consistent estimator of $\Sigma$.
\end{proposition}

In the case where $|\mathcal{B}_0| = 1$, meaning that there exists a unique $\beta_0 \in \mathbb{R}^d $ such that $Y \stackrel{d}{=} \beta^T_0  X + \epsilon$ we show that $\hat{\beta}_n$ is a consistent estimator in the classical sense. In Theorem \ref{conv2}, we show that under this uniqueness assumption $\hat{\beta}_n\rightarrow\beta_0$ in probability. 
\begin{theorem}\label{conv2}
Assume that $F^\epsilon$ admits a density $f^\epsilon$ and $|\mathcal{B}_0| = 1$. In addition assume that $\PP(\alpha^T X = 0) = 0$ for any $\alpha\in \RR^d$ such that $\alpha\neq 0$. Then for all $\hat{\beta}_n\in\mathcal{B}_{n, 2}$ we have 
\begin{eqnarray*}
\hat{\beta}_n  \stackrel{\PP}{\rightarrow} \beta_0.
\end{eqnarray*}
\end{theorem}

When $|\mathcal{B}_0| > 1$, then depending on the distribution of $X$, we might be able to narrow down $\mathcal{B}_{n, 2}$ by more than just specifying the norm of the members. In Theorem~\ref{pole} below, we show that for a certain family of probability distributions, $\mathcal{B}_0$ consists of all vectors that result from permuting the components of a member $\beta_0$.
\begin{theorem}\label{pole}
Suppose that the components of the covariate $X$, $X^1, \cdots, X^d$ are i.i.d. with moment generating function that takes the form
\begin{eqnarray*}
G_{X^1}(t) := E[e^{t X^1}] =  \frac{1}{(\lambda - t)^{\alpha}} H(t), \  \ \text{for $t < \lambda$},
\end{eqnarray*}
where $0 < \lambda$, $\alpha \in (0, \infty)$, and $H$ is a continuous function such that
\begin{eqnarray*}
\lim_{t \nearrow \lambda}  H(t)  \in (0, \infty).
\end{eqnarray*}  
Then the set of all $\beta$ such that $Y \stackrel{d}{=} \beta_0^T X + \epsilon$ is 
\begin{eqnarray*}
\mathcal B_0 = \left \{\beta \in \RR^d:  \exists \ \text{permutation} \ \pi \ \text{such that}  \ \beta_{0, i} = \beta_{\pi(i)}, \  i=1, \ldots, d  \right \}.
\end{eqnarray*}
\end{theorem}
Note that this family includes exponential and Gamma distributions as well as any convolution thereof.   

Finally, in Theorem~\ref{limitDist}, where there exists a unique regression vector for which the unmatched regression model holds we are able to show that under some smoothness assumptions the DLSE $\hat{\beta}_n$ is asymptotically normal with mean $\beta_0$ and covariance matrix that depends on the distribution of $X$, that of $\epsilon$ and the true regression vector $\beta_0$.

\begin{theorem}\label{limitDist}
Suppose $F^\epsilon\in C^2$ and there exists $M > 0$ such that $f^\epsilon \le M$ and $\vert (f^\epsilon)^{(1)} \vert \le M$. Also assume that there exists an integer $m \ge 1$ and $a_1 < \cdots < a_m$ such that $f^\epsilon$ is monotone on $(-\infty, a_1], (a_1, a_2), \cdots, [a_m, \infty)$. Furthermore, we assume that $\EE\Vert X \Vert^2_2< \infty$, and the matrix
\begin{eqnarray*}
U = \int \left(\int x f^\epsilon(y - \beta^T_0 x)  dF^X(x)\right) \left(\int x^T f^\epsilon(y - \beta^T_0 x)  dF^X(x)\right)  dF^Y(y)
\end{eqnarray*}
is positive-definite. Then when $|\mathcal{B}_0| = 1$ we have 
\begin{eqnarray*}
\sqrt n (\hat{\beta}_n - \beta_0) \stackrel{d}{\rightarrow} U^{-1} \mathbb{V} 
\end{eqnarray*}
where
\begin{eqnarray*}
\MoveEqLeft \mathbb{V} = - \int_{\RR}\left(\mathbb B_1 \circ F^Y(y) - \int\mathbb  B_2 \circ F^Z(z) f^\epsilon(y - z) dz \right) \\
& \left(\int x f^\epsilon(y - \beta^T_0 x)  d F^X(x)\right) d F^Y(y), \ 
\end{eqnarray*}
with $Z = \beta^T_0 X$, and $\mathbb{B}_1$ and $\mathbb{B}_2$ are two independent standard Brownian Bridges from $(0,0)$ to $(1, 0)$. The random vector $U^{-1}\mathbb{V}$ is distributed as $\mathcal{N}(\mathbf{0}, U^{-1} \EE[\mathbb{V} \mathbb{V} ^T]U^{-1})$.
\end{theorem}
The conditions required for $f^\epsilon$ are satisfied by most of the well-known densities, including normal, Laplace, symmetric Gamma densities and any finite convolution thereof. The asymptotic covariance matrix has a complicated dependence on the parameters of the model. However, the exact knowledge of this asymptotic covariance matrix is not at all necessary to make useful inferences. In fact, the asymptotic result of Theorem~\ref{limitDist} allows us to use re-sampling techniques to infer the true regression vector. Therefore, one can resort to bootstrap in order to find an approximation of asymptotic confidence bounds for $\beta_0$. 

Theorem~\ref{limitDist} provides the asymptotic normality of the DLSE only under the identifiability of the unmatched linear regression model, and hence it cannot be used beyond the scope where $|\mathcal{B}_0| = 1$.

When we are in the situation of Theorem~\ref{pole}, and if the components of $\beta_0 \in \mathcal{B}_0$ are all distinct, then the model becomes identifiable if the components of $\beta_0$ are ordered from smallest to largest. In fact, the arguments in the proof of Theorem~\ref{orderedBeta} can be used again to show the following result. 
\begin{theorem}\label{orderedBeta}
Suppose that $F^\epsilon$ satisfies the same conditions as in Theorem~\ref{pole}. Suppose also that the assumption of Theorem~\ref{pole} holds and that the components of the vector $\beta_0 \in \mathcal{B}_0$ are all distinct. Denote $b_0 = [\beta^{(1)}_0, \cdots, \beta^{(d)}_0]^T$ the vector of ordered components of $\beta_0$ from smallest to largest, and $\hat{b}_n = [\hat{\beta}^{(1)}_n, \cdots, \hat{\beta}^{(d)}_n]^T$ the vector of ordered components of the DLSE, $\hat{\beta}_n$. Then,
\begin{eqnarray*}
\sqrt{n}(\hat{b}_n - b_0)\stackrel{d}{=}U^{-1}\mathbb{V},
\end{eqnarray*}
where $U$ and $\mathbb{V}$ are defined in Theorem~\ref{pole}.
\end{theorem}

To close this section, we would like to stress the fact that the fast rate of convergence, $1/\sqrt n$, is obtained under the very important condition of identifiability of $\beta_0$, or in other words its uniqueness. When this condition is satisfied,  arguments from the theory of M-estimators can be used. In this case, the $1/\sqrt n$-rate can be obtained since the estimation problem is fully parametric. In case identifiability is not satisfied, Theorem \ref{pole} provides a sufficient condition on the distribution of the covariates which guarantees identifiability of the ordered values of the regression vector. Although the condition given in that theorem encompasses many distributions, a unified result with  more general identifiability conditions is yet to be established.

\subsection{Application to Semi-supervised Learning}
As established above, the estimator $\hat{\beta}_n$ is not consistent in case of non-identifiability. Therefore, one cannot possibly think of $\hat{\beta}_n$ as something that will be close to some unique truth since the latter does not even exist.  However, there are situations where some feature of the model is unique, and in this case, one would expect that the DLSE succeeds in remaining faithful to such features.  

Consider the case where $X\sim\mathcal{N}(0, \Sigma)$ with a non-singular covariance matrix $\Sigma$.  Lemma~\ref{Charbeta} describes the set $\mathcal{B}_0$ in terms of a norm with respect to $\Sigma$: There exists a constant $c > 0$ such that for any $\beta\in\mathcal{B}_0$ we have $\Vert\beta\Vert_{2, \Sigma} = c$. In addition, Proposition~\ref{Gaussian} guarantees that although we do not have consistency of $\hat{\beta}_n$, it holds that $\Vert\hat{\beta}_n\Vert_{2, \widehat{\Sigma}_n}\stackrel{\PP}{\rightarrow}c$. In this case if we have access to a third sample of matched data, we can benefit from the result of Proposition~\ref{Gaussian} together with an estimate of $\beta_0$ from the matched sample. 

Let $\{(\tilde{Y}_i, \tilde{X}_i)\}_{i=1}^m$ be a set of matched data such that $Y_i = \beta_0^T X_i + \epsilon_i$. Take $\tilde{\beta}_m$ to be an estimator of $\beta_0$ based on only $\{(\tilde{Y}_i, \tilde{X}_i)\}_{i=1}^m$. For example, one can take $\tilde{\beta}_m$ to be the ordinary least squared (OLS) estimator. Let $\widehat{\Sigma}_n$ be a consistent estimator of $\Sigma$, the covariance matrix of $X$ using the sample $\{X_i\}_{i=1}^n$ (here, an estimate based on this sample is more accurate since $n$ is much bigger than $m$). Now consider the following modified estimator
\begin{eqnarray}\label{betaDagger}
\beta_{n, m}^\dagger = \tilde{\beta}_m\frac{\Vert\hat{\beta}_n\Vert_{2, \widehat{\Sigma}_n}}{\Vert\tilde{\beta}_m\Vert_{2, \widehat{\Sigma}_n}}.
\end{eqnarray}
When $n \gg m$, we expect that $\Vert\hat{\beta}_n\Vert_{2, \widehat{\Sigma}_n}$ to be a better estimate of $\Vert\beta_0\Vert_{2, \Sigma}$ compared to $\Vert\tilde{\beta}_m\Vert_{2, \widehat{\Sigma}_n}$. Therefore, while the estimate based on the unlinked data fails to provide a meaningful estimate of the direction of $\beta_0$, we can still use $\Vert\hat{\beta}_n\Vert_{2, \widehat{\Sigma}_n}$ to modify the norm of $\tilde{\beta}_m$ and hence improve the performance of the latter. 

Unfortunately, we do not yet have a rate of convergence for $\Vert\hat{\beta}_n\Vert_{2, \widehat{\Sigma}_n}$. In Section~\ref{AppSynRealData} we provide experimental results suggesting that for large enough $n$ one can use $\Vert\hat{\beta}_n\Vert_{2, \widehat{\Sigma}_n}$ as an estimate of $\Vert\beta_0\Vert_{2, \Sigma}$. As we do not have the right arguments which show theoretically that the modified estimator $\beta_{n, m}^\dagger$ in \eqref{betaDagger} improves the performance $\tilde{\beta}_m$, we resort to another estimator of the norm $\Vert\beta_0\Vert_{2, \Sigma}$ and whose convergence rate can be established. The estimator is based on the following simple observation:
\begin{eqnarray*}
\var(Y) = \beta_0^T \Sigma \beta_0 + \var(\epsilon) = \Vert\beta_0\Vert_{2, \Sigma}^2 + \var(\epsilon).
\end{eqnarray*}
Therefore, under the assumption of known distribution of $\epsilon$, it is possible to estimate $\Vert\beta_0\Vert_{2, \Sigma}$ by estimating $\var(Y)$. The following remark shows that when $m = o(n)$, using the information from unmatched data improves the performance of the OLS estimator.

\begin{remark}\label{improved}
Let $\tilde{\beta}_m$ be the ordinary least square estimate of $\beta_0$ using the matched data $\{(\tilde{Y}_i, \tilde{X}_i)\}_{i=1}^m$. Let $r_n^2 = \frac{1}{n-1}\sum_{i=1}^n(Y_i - \bar{Y}_n)^2 - \var(\epsilon)$ where $\bar{Y}_n$ is the sample mean of $\{Y_i\}_{i=1}^n$. Define
\[
\tilde{\beta}_{n, m} = \tilde{\beta}_m \frac{r_n}{\Vert\tilde{\beta}_m\Vert_{2, \widehat{\Sigma}_n}},
\]
where $\widehat{\Sigma}_n$ is the empirical covariance matrix of $\{X_i\}_{i=1}^n$. When $m = o(n)$
\[
\EE[\Vert\tilde{\beta}_{n, m} - \beta_0\Vert^2] < \EE[\Vert\tilde{\beta}_m - \beta_0\Vert^2],
\]
where expectations on the right and left sides are taken with respect to the distribution of the matched data and the product distribution of the matched and unmatched, respectively. 
\end{remark}

\section{Applications to Synthetic and Real Data}\label{AppSynRealData}
In this section, we present the result of several simulations with the goal of illustrating the theoretical results derived above. 
In addition, we shall consider real applications using two different real data sets: The (unmatched) inter-generational mobility data set already analyzed by several authors; e.g. \cite{olivetti2015name} and the (matched) Power Plant data set to which we apply ideas from semi-supervised learning. Computation of the DLSE, especially for large data sets (synthetic or not) and for several runs impose numerical challenges. One major issue is that the minimization problem is generally not convex and can be multi-modal especially in the case of non-identifiability.  We use the default setting of the function \lq\lq optim\rq\rq from the package \lq\lq stats\rq\rq  \ of the open software R. The default method of optimization is the method introduced by \cite{neldermead}. 

\subsection{Synthetic Data}
\begin{example}\label{normalExp}
For this example, we generated 1000 independent samples of $Y$ and $X$ where $Y\stackrel{d}{=} \beta_0^T X + \epsilon$ of size $n = 4000$ such that
\begin{eqnarray*}
&X = [X^1, X^2]^T\in\RR^2, \qquad X^1\sim\mathcal{N}(1, 1), \qquad  X^2\sim\text{Exp}(1), \\ &\epsilon\sim\mathcal{N}(0, 1), \qquad\beta_0^T = [1, 2],
\end{eqnarray*}
where $X^1, X^2$ and $\epsilon $ are independent. In this case, and as shown in the supplementary material, $\mathcal{B}_0$ contains only $\beta_0$. Therefore, we expect that $\hat{\beta}_n$ is a consistent estimator of $\beta_0$. Figure~\ref{fig:boxplot_NormalExp_mu1_n4000_it100} shows the boxplot of $\hat{\beta}_n^1$ and $\hat{\beta}_n^2$ based on the 100 replications.
\begin{figure}
    \centering
    \begin{subfigure}[b]{0.4\textwidth}
        \includegraphics[width=\textwidth]{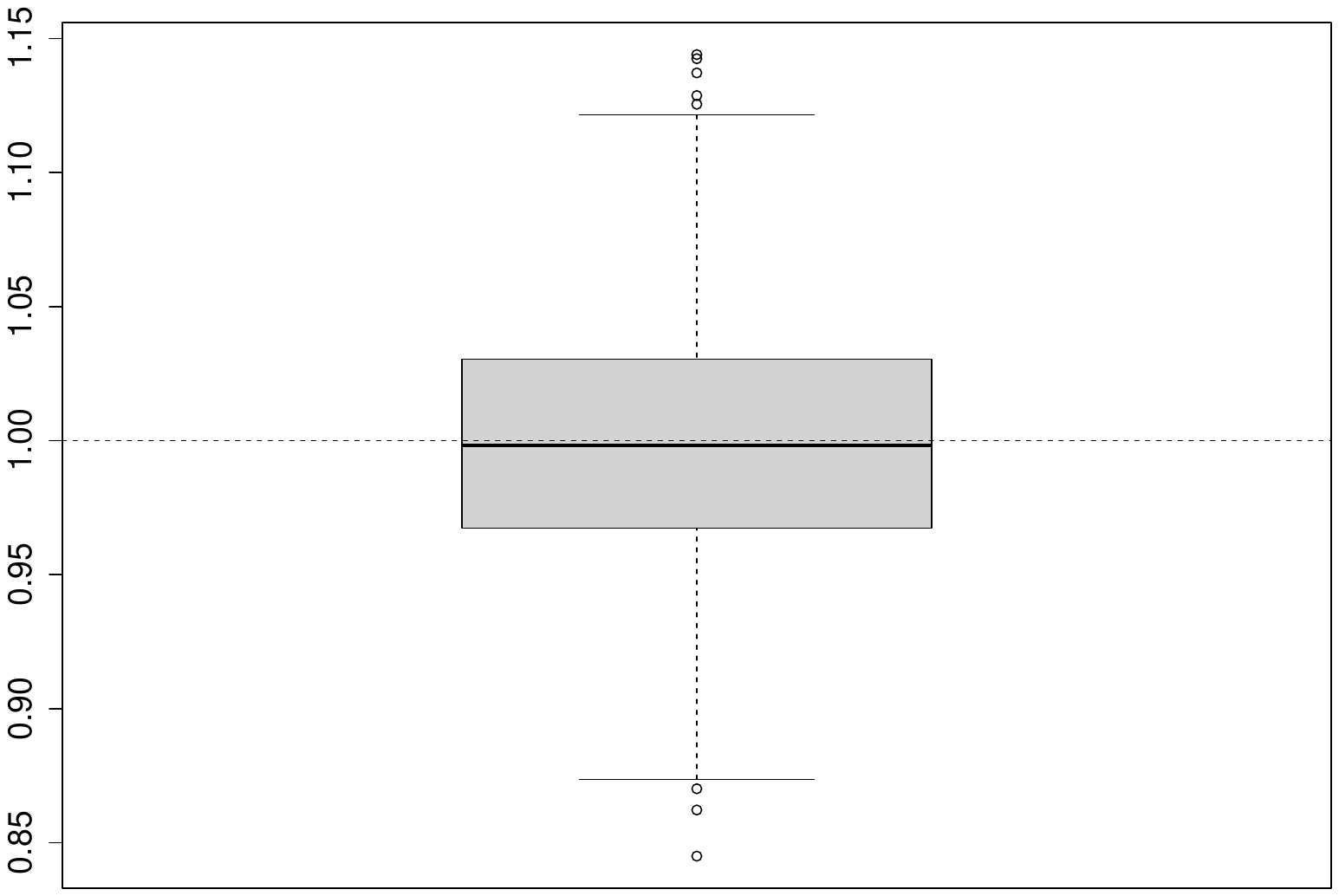}
        \caption{$\hat{\beta}_n^1$}
    \end{subfigure}
    ~ 
    \begin{subfigure}[b]{0.4\textwidth}
        \includegraphics[width=\textwidth]{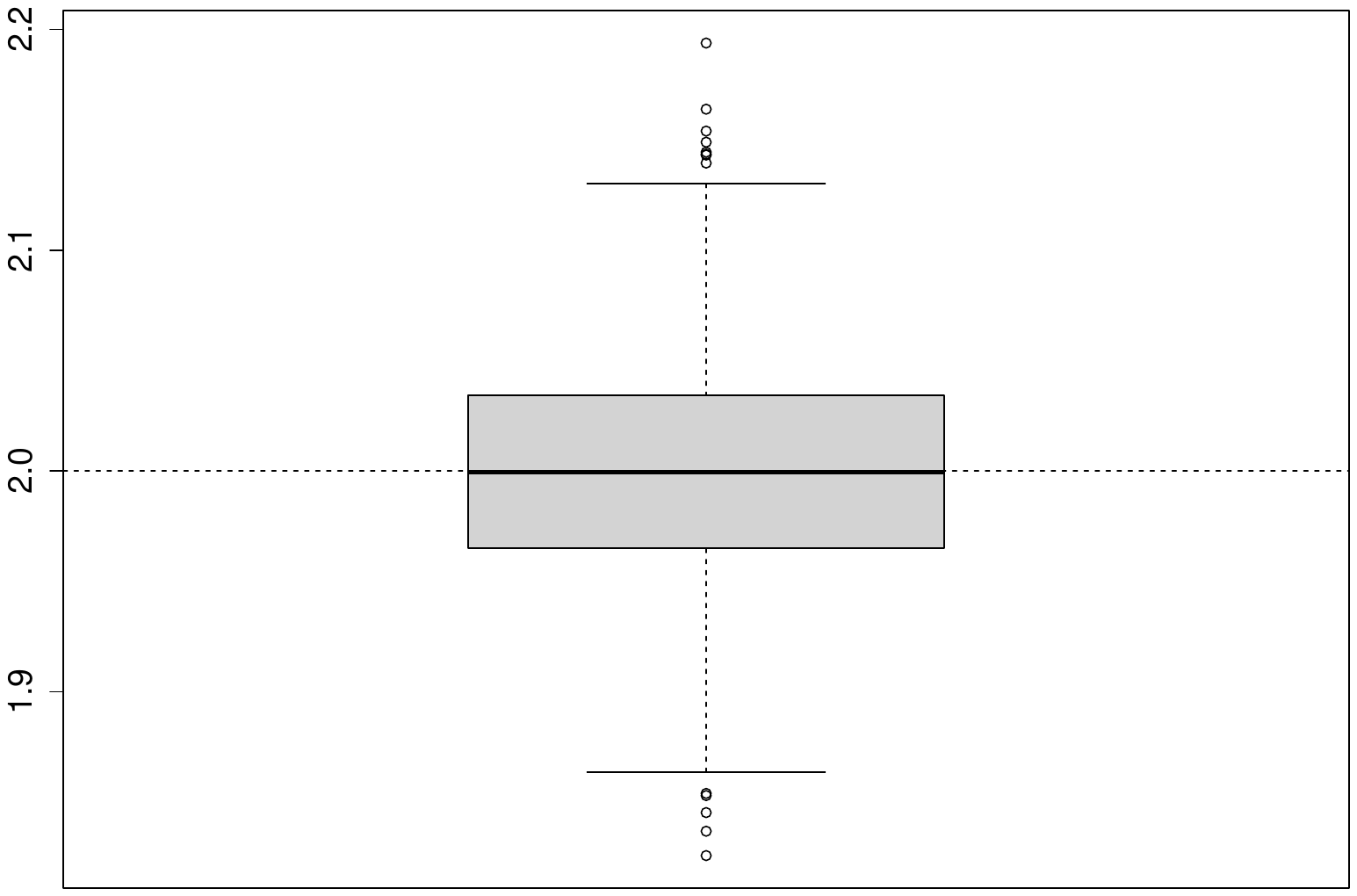}
        \caption{$\hat{\beta}_n^2$}
    \end{subfigure}
    \caption{\label{fig:boxplot_NormalExp_mu1_n4000_it100} Boxplots of $\hat{\beta}_n^1$ and $\hat{\beta}_n^2$ for 100 independent samples of size $n=4000$. Root mean squared error of $\hat{\beta}_n^1$ is $0.047$, and root mean squared error of $\hat{\beta}_n^2$ is $0.052$.}
\end{figure}

Theorem~\ref{limitDist} shows that $\sqrt{n}(\hat{\beta}_n - \beta_0)$ has asymptotically a normal distribution. This type of result can be validated by plotting the quantiles $\sqrt{n}(\hat{\beta}_n - \beta_0)$ against those of a standard normal variable. Such plot is commonly known under the name of \textit{qqplot}. Figure~\ref{fig:asympNormality} shows such qqplots for the components of $\sqrt{n}(\hat{\beta}_n - \beta_0)$. Their linear shape is much aligned with the asymptotic normality of our estimator as stated in Theorem \ref{limitDist}. 
\begin{figure}
    \centering
    \begin{subfigure}[b]{0.4\textwidth}
        \includegraphics[width=\textwidth]{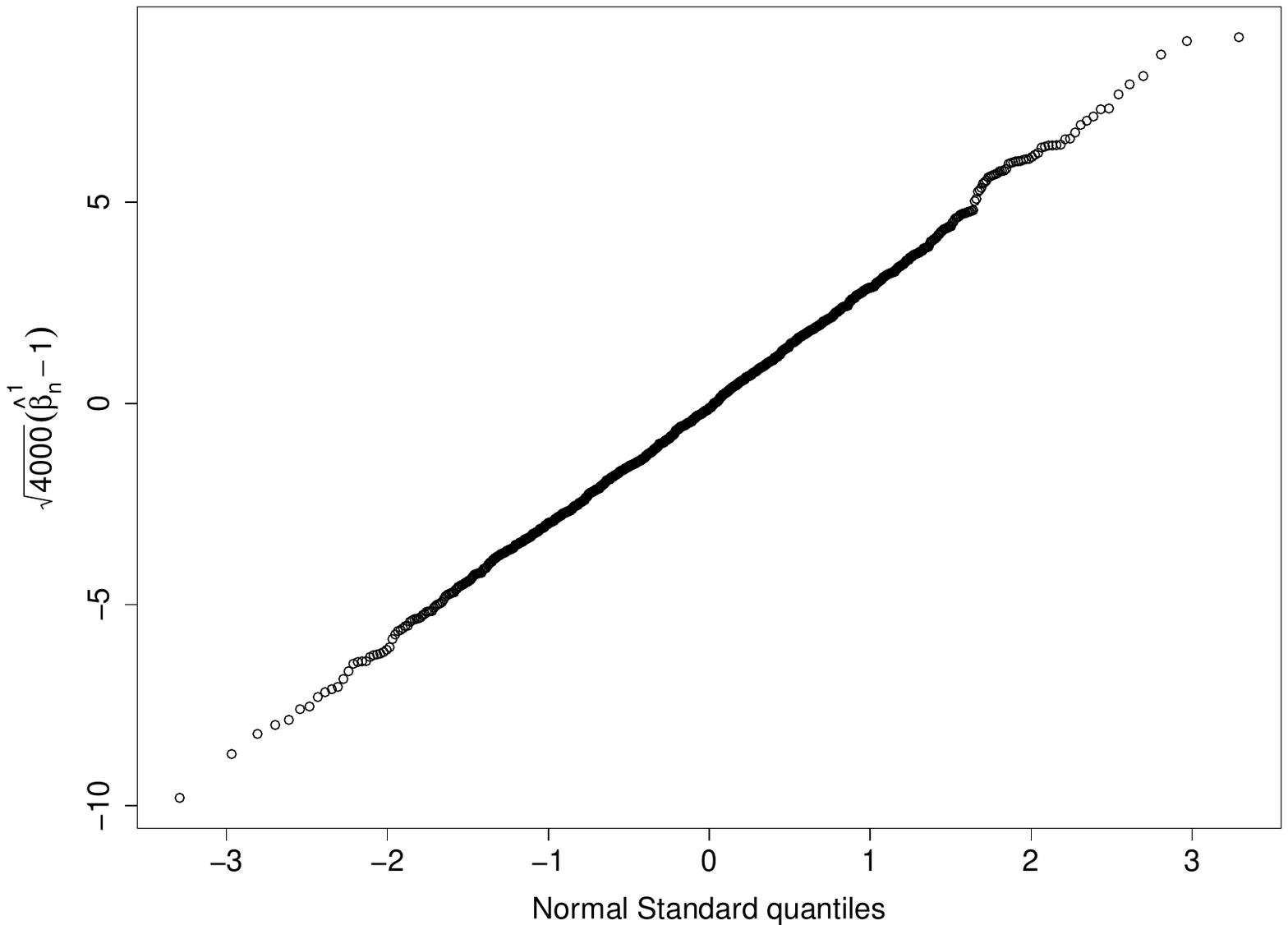}
    \end{subfigure}
    ~ 
    \begin{subfigure}[b]{0.4\textwidth}
        \includegraphics[width=\textwidth]{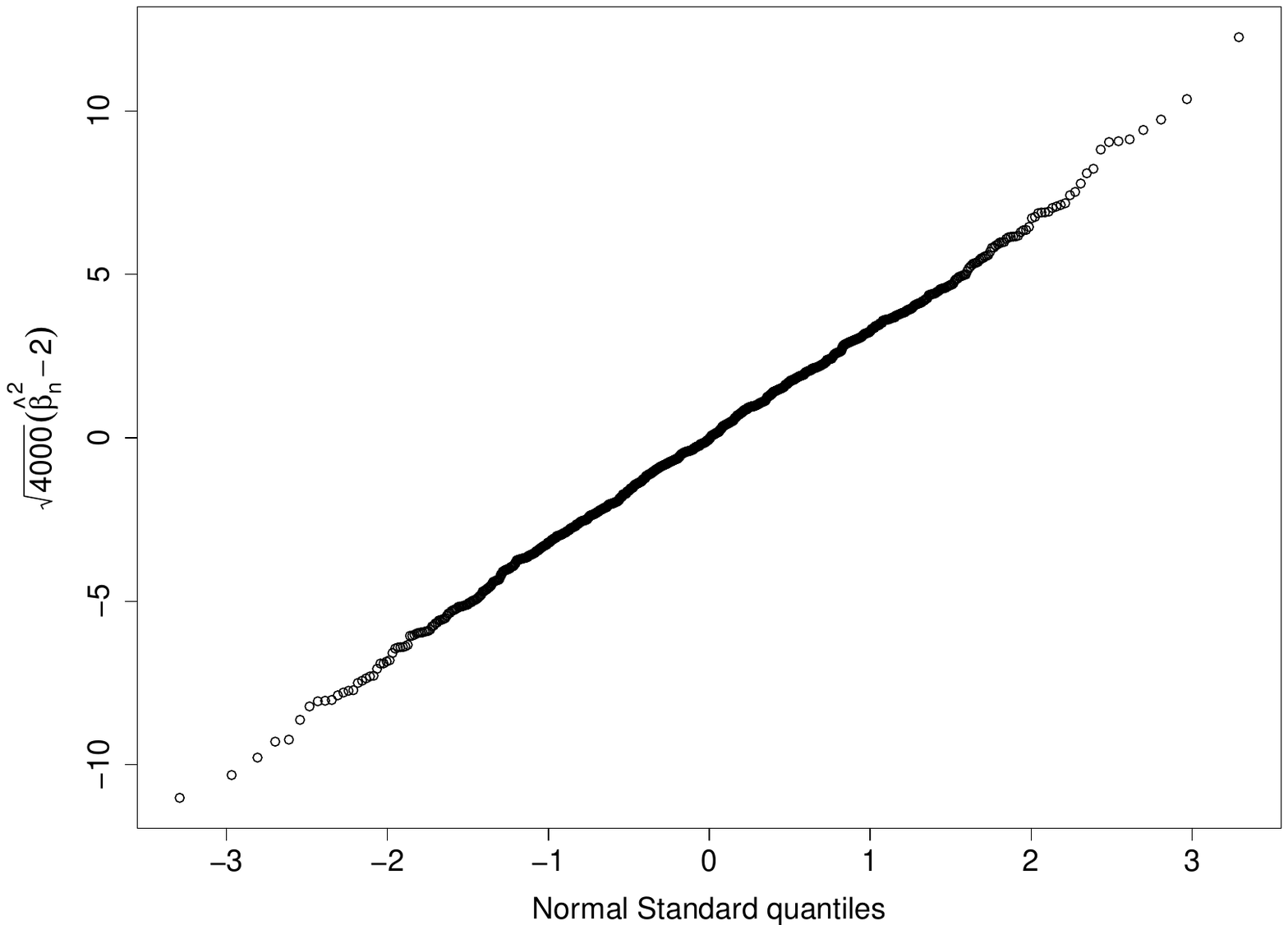}
    \end{subfigure}
    \caption{\label{fig:asympNormality} qqplot of $\sqrt{n}(\hat{\beta}_n^1 - 1)$ and $\sqrt{n}(\hat{\beta}_n^2 - 2)$ against $\mathcal{N}(0, 1)$, with $n=4000$. These plots are in line with the asymptotic result of Theorem~\ref{limitDist}.}
\end{figure}

Using unmatched data, at best, means that we only have access to the generating mechanism of $X$ and $Y$ without knowing the link between them. Lack of information about the link between the $X$'s and $Y$'s should come with a price paid on the performance of the estimate of $\beta_0$. To see how much our DLSE $\hat{\beta}_n$ suffers from this lack of knowledge, we compare the performance of $\hat{\beta}_n$ with the ordinary least square estimator $\tilde{\beta}_n$ using a matched sample of size $n=4000$ over 100 independent replications. The comparison is done in terms of absolute prediction error. For this purpose we vary the noise strength level $\epsilon\sim\mathcal{N}(0, \sigma_\epsilon^2)$ for $\sigma_\epsilon\in\{0.6, 0.8, 1, 1.2\}$ to see how much it effects the performance. In Table~\ref{tab:errorRatio}, the result is summarized in terms of the ratio of the absolute error of the OLS estimator to the absolute error of DLSE for each value of $\sigma_\epsilon$ and for each component of $\beta_0$. Clearly, $\tilde{\beta}_n$ beats $\hat{\beta}_n$ as it uses more information, but given the fact that $\hat{\beta}_n$ is oblivious to the link between the response and covariate, its performance is rather quite satisfactory.
\begin{table}[H]
\centering
\small{
\begin{tabular}{rrrrrr}
  \hline
  \multicolumn{1}{c}{} & $\sigma_\epsilon = 0.6$ & $\sigma_\epsilon = 0.8$ & $\sigma_\epsilon = 1$ & $\sigma_\epsilon = 1.2$ & $\sigma_\epsilon = 2$ \\
  \hline
 $\text{mean}(|\tilde{\beta}_n^1 - \beta_0^1|) / \text{mean}(|\hat{\beta}_n^1 - \beta_0^1|)$ & 0.14 & 0.26 & 0.32 & 0.28 & 0.61\\
  \hline
$\text{mean}(|\tilde{\beta}_n^2 - \beta_0^2|) / \text{mean}(|\hat{\beta}_n^2 - \beta_0^2|)$ & 0.18 & 0.16 & 0.18 & 0.25 & 0.55\\
\end{tabular}
}
\caption{\label{tab:errorRatio} Ratio of mean absolute error of $\tilde{\beta}_n$ to $\hat{\beta}_n$ as a function of the standard deviation of the noise.}
\end{table}
\end{example}

\begin{example}\label{normal}
We consider the following setting
\begin{eqnarray*}
    Y\stackrel{d}{=} \beta_0^T X + \epsilon,\qquad
    X\sim\mathcal{N}(0, I_2), \qquad \epsilon\sim\mathcal{N}(0, 1), \qquad\beta_0^T = [1, 2],
\end{eqnarray*}
where $I_2$ is the $2\times 2$ identity matrix. Note that in this case, $\mathcal{B}_0$ contains more than one element. In fact by Lemma~\ref{Charbeta} we know that $\mathcal{B}_0 = \{\beta\in\RR^2: \Vert\beta\Vert_{2, \Sigma} = \sqrt{5}\}$. Proposition~\ref{Gaussian} guarantees that $\Vert\hat{\beta}_n\Vert_{2, \Sigma}$ converges to $\|\beta_0\| = \sqrt 5$ in probability. The scatterplot of $\hat \beta_n$ shown in Figure~\ref{fig:scatterEx2} gives a clear illustration of this fact. The mean and standard deviation of $\Vert\hat{\beta}_n\Vert_2^2$ over the 1000 replications with sample size $n = 4000$ were found to be $5.01$ and $0.15$ respectively.
\begin{figure}[H]
  \centering
    \includegraphics[width=0.7\textwidth]{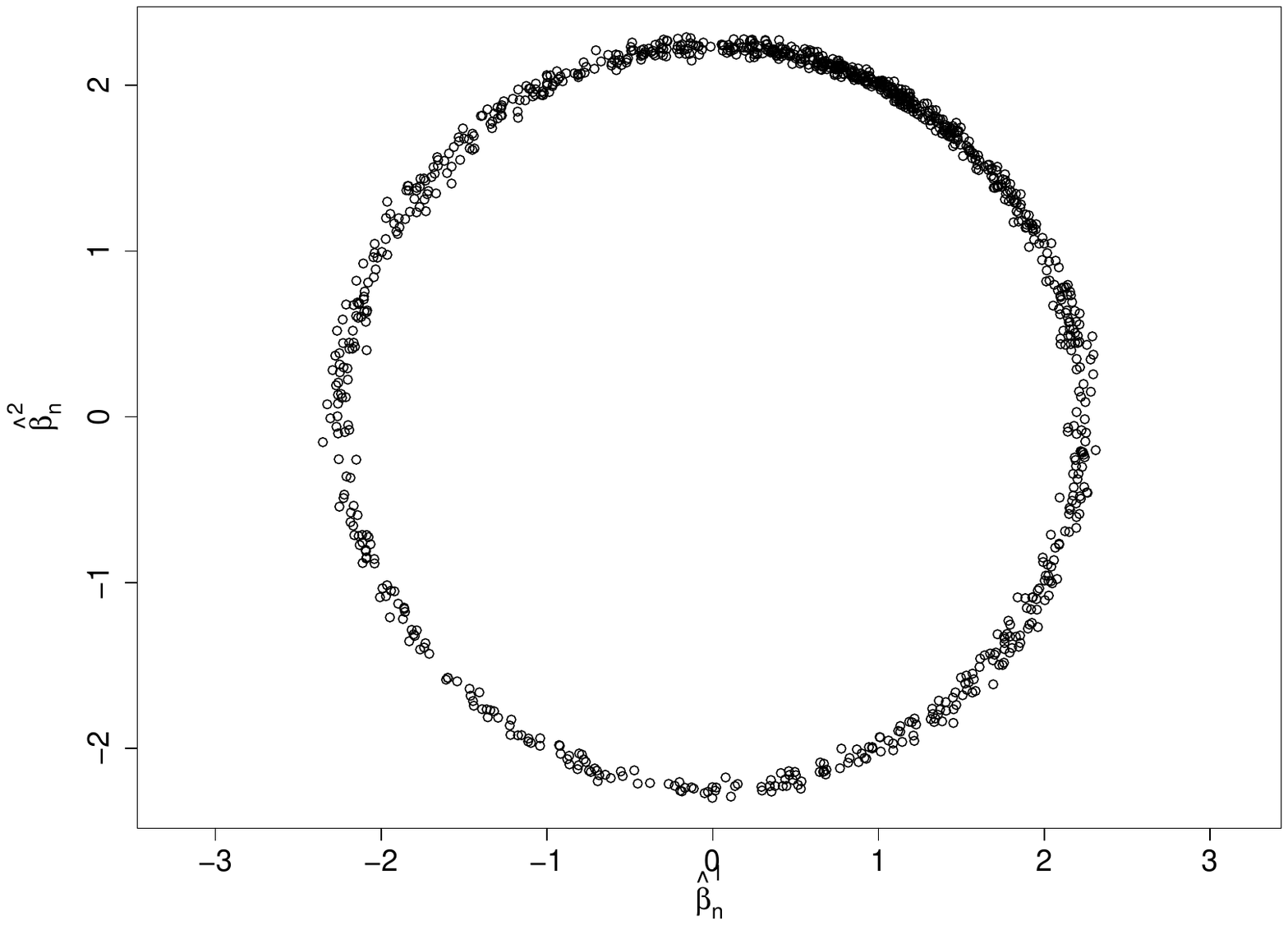}
    \caption{\label{fig:scatterEx2}Plot of $\hat{\beta}_n$ for $X\sim\mathcal{N}(0, I_2)$ and $\beta_0^T = [1, 2]$. In this case $\mathcal{B}_0 = \{\beta\in\RR^2: \Vert\beta\Vert_2 = \sqrt{5}\}$.}
\end{figure}
\end{example}

\begin{example}\label{nonunique_exp}
In this example, we generated 1000 times independent samples of $Y$ and $X$ where $Y\stackrel{d}{=} \beta_0^T X + \epsilon$ with size $n = 8000$ such that
\begin{eqnarray*}
& X = [X^1, X^2, X^3]^T\in\RR^3, \qquad X^i\sim\text{Exp}(1) \text{ for $i = 1, 2, 3$ }, \\
&\epsilon\sim\mathcal{N}(0, 1), \qquad\beta_0^T = [1, 2, -1]
\end{eqnarray*}
where $X^1,  X^2, X^3$ and $\epsilon$ are independent.  Theorem~\ref{pole} implies that $|\mathcal{B}_0| = 3! = 6$, since all permutations of $\beta_0$ lead to the same distribution. Instead of $\hat{\beta}_n$ we take a look at $\hat{\beta}_{n,\text{ordered}} = [\hat{\beta}_n^{(1)}, \hat{\beta}_n^{(2)}, \hat{\beta}_n^{(3)}]^T$ where $\hat{\beta}_n^{(i)}$ is the $i$-th smallest element of $\hat{\beta}_n$. As mentioned earlier, when $\vert \mathcal{B}_0 \vert > 1$, we do not have any asymptotic normality result but believe that $\hat{\beta}_{n,\text{ordered}}$ is asymptotically normal under some regularity assumptions. Figure~\ref{fig:exp_dim3_8000_qqplot} shows the qqplots of $\sqrt{8000}(\hat{\beta}_n^{(i)} - \beta_0^{(i)})$ against the quantiles of $\mathcal{N}(0,1)$, and which is supported by Theorem~\ref{orderedBeta}.

\begin{figure}
    \centering
    \begin{subfigure}[b]{0.3\textwidth}
        \includegraphics[width=\textwidth]{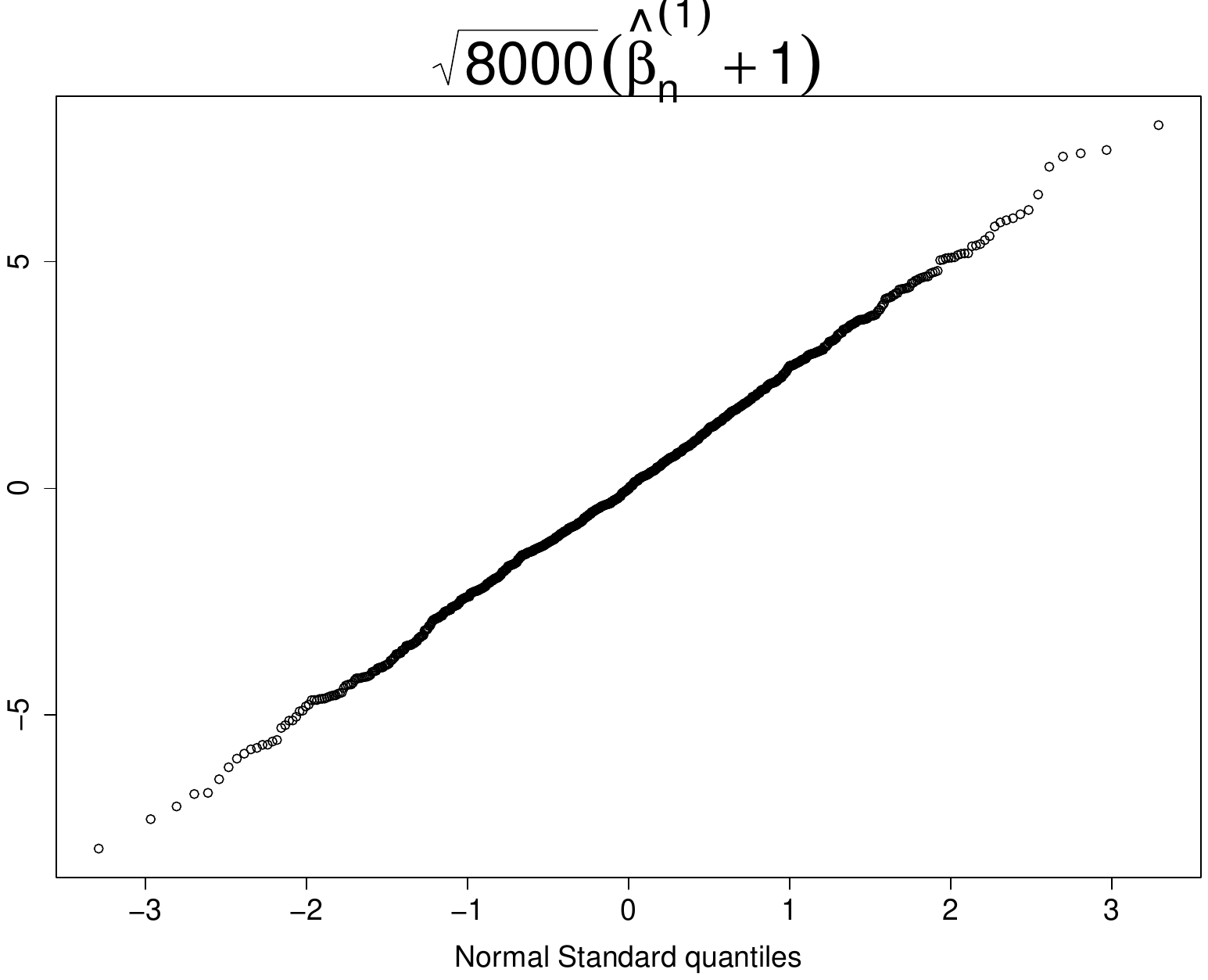}
    \end{subfigure}
    ~ 
    \begin{subfigure}[b]{0.3\textwidth}
        \includegraphics[width=\textwidth]{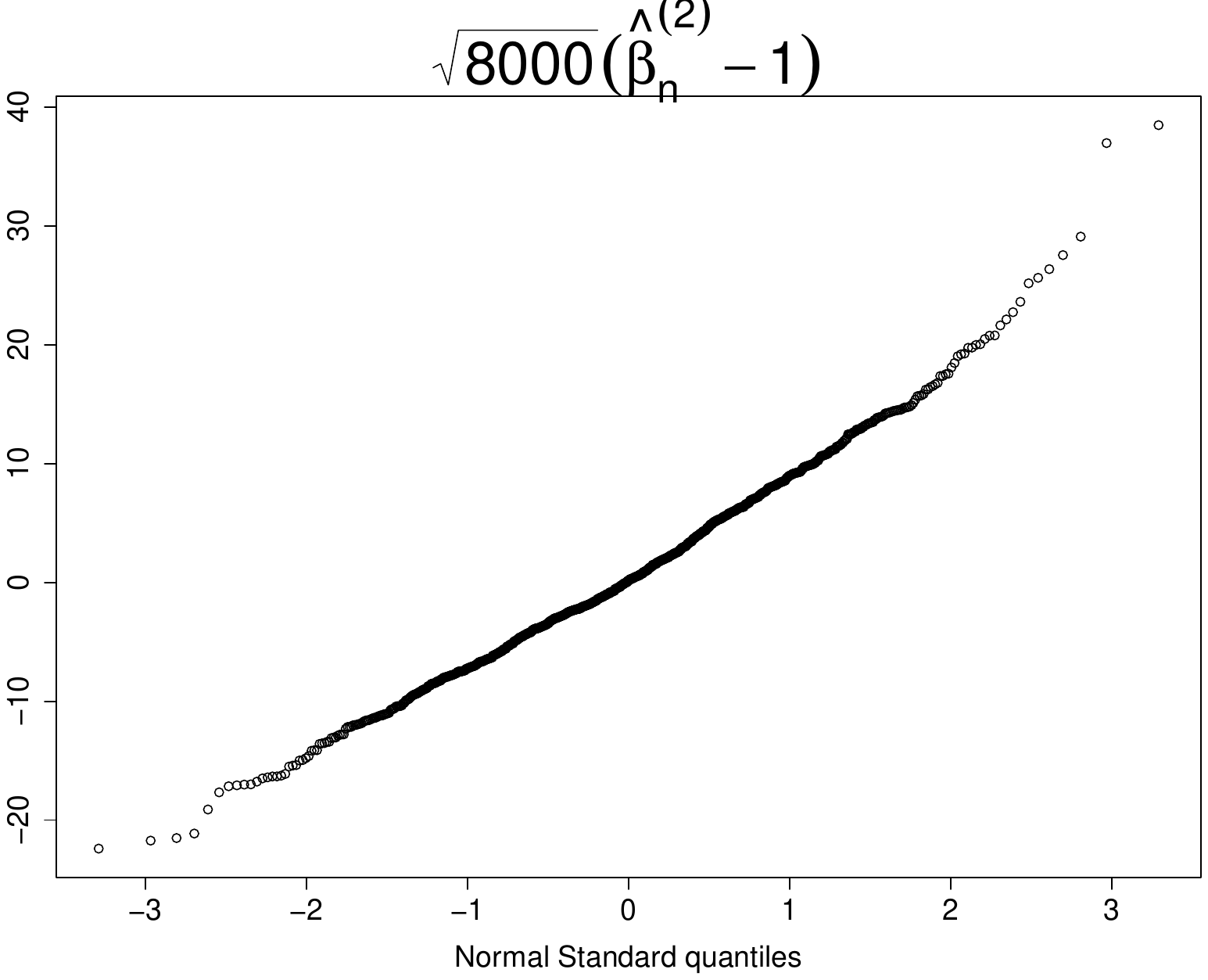}
    \end{subfigure}
    ~ 
    \begin{subfigure}[b]{0.3\textwidth}
        \includegraphics[width=\textwidth]{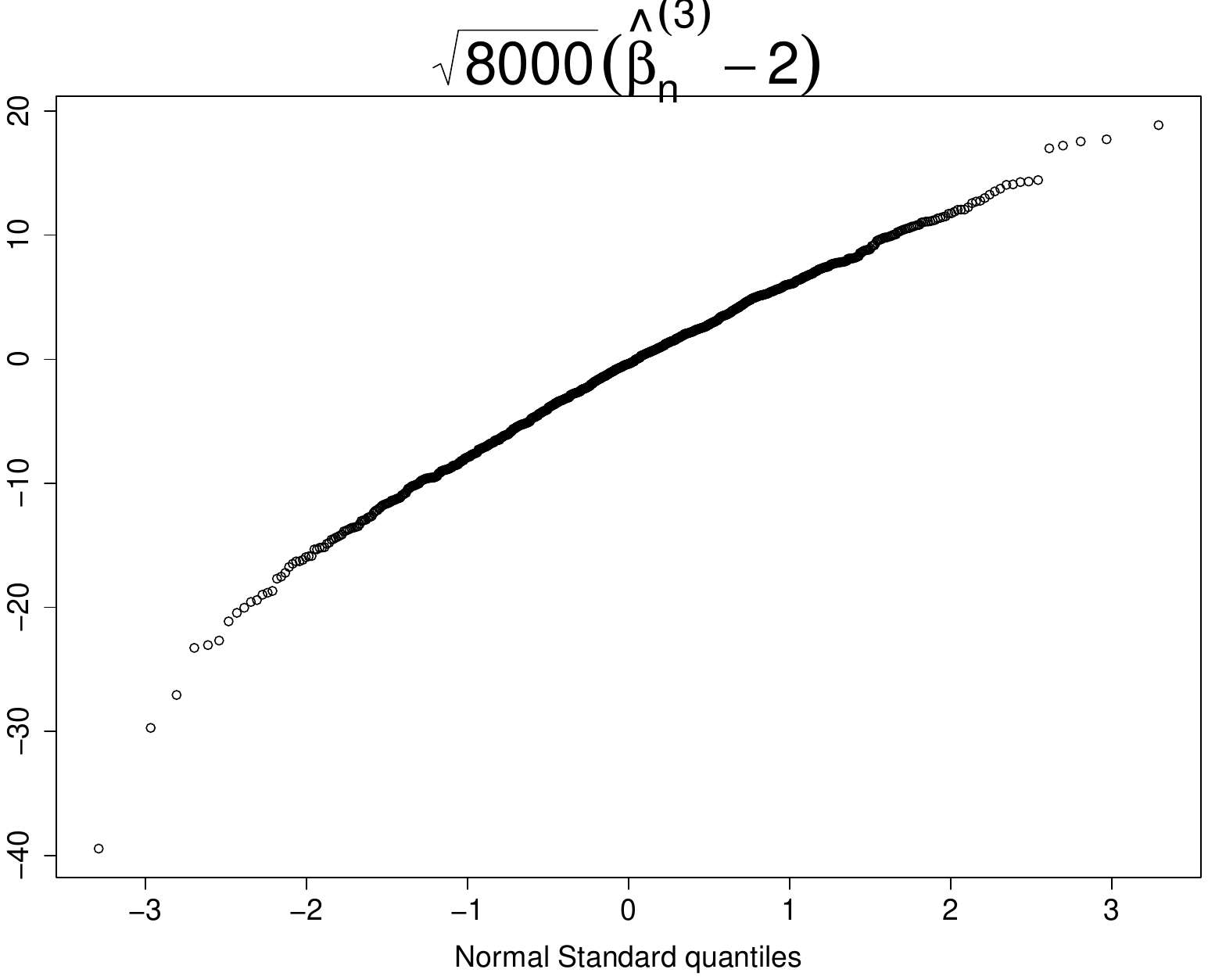}
    \end{subfigure}
    \caption{\label{fig:asympNormality} qqplot of $\sqrt{n}(\hat{\beta}_n^{(1)} + 1)$, $\sqrt{n}(\hat{\beta}_n^{(2)} - 1)$, and $\sqrt{n}(\hat{\beta}_n^{(3)} - 3)$ against $\mathcal{N}(0, 1)$, with $n = 8000$. These plots are in line with the asymptotic result of Theorem~\ref{orderedBeta}.}
\end{figure}

\begin{figure}
    \centering
    \begin{subfigure}[b]{0.3\textwidth}
        \includegraphics[width=\textwidth]{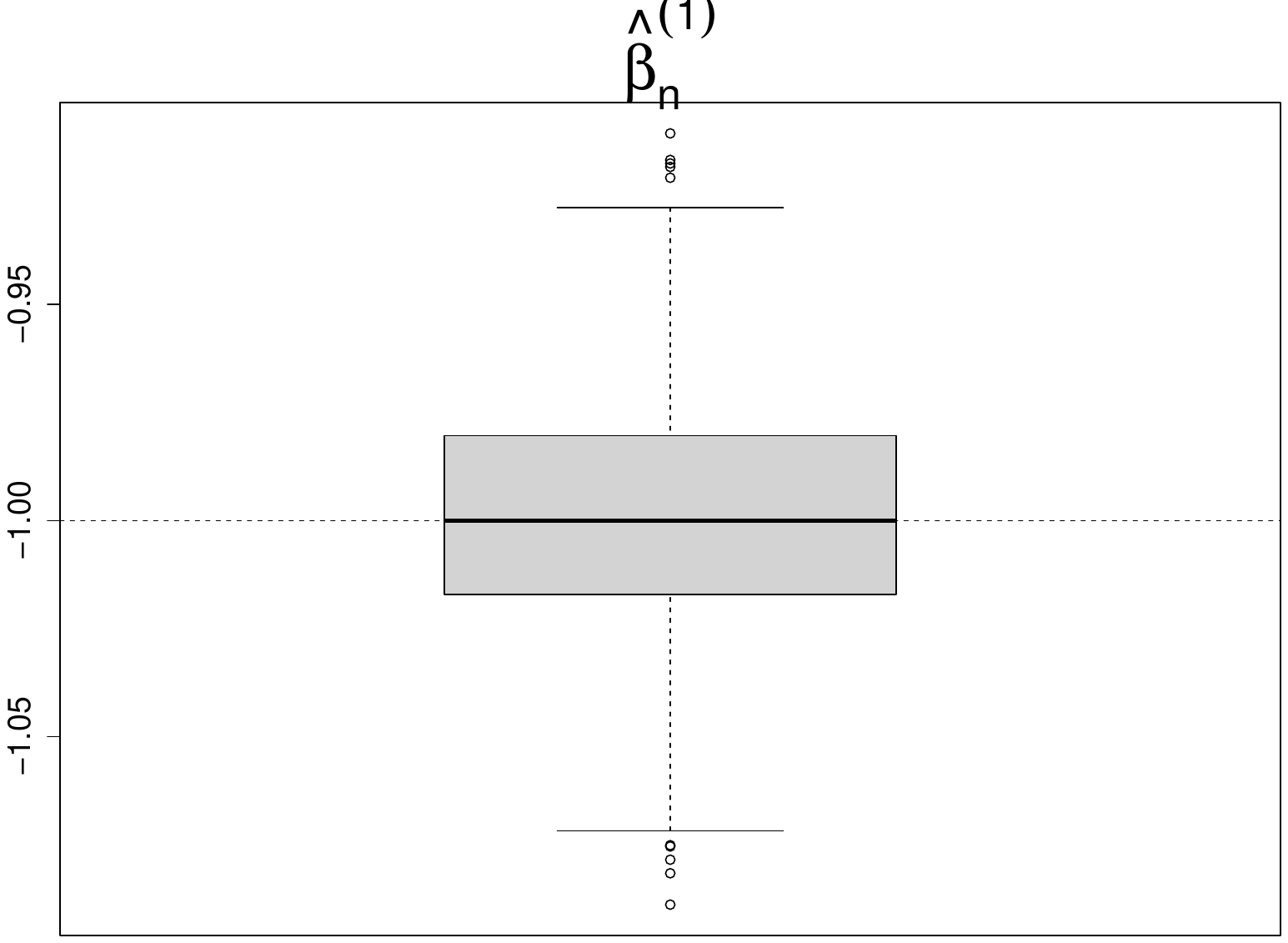}
    \end{subfigure}
    ~ 
    \begin{subfigure}[b]{0.3\textwidth}
        \includegraphics[width=\textwidth]{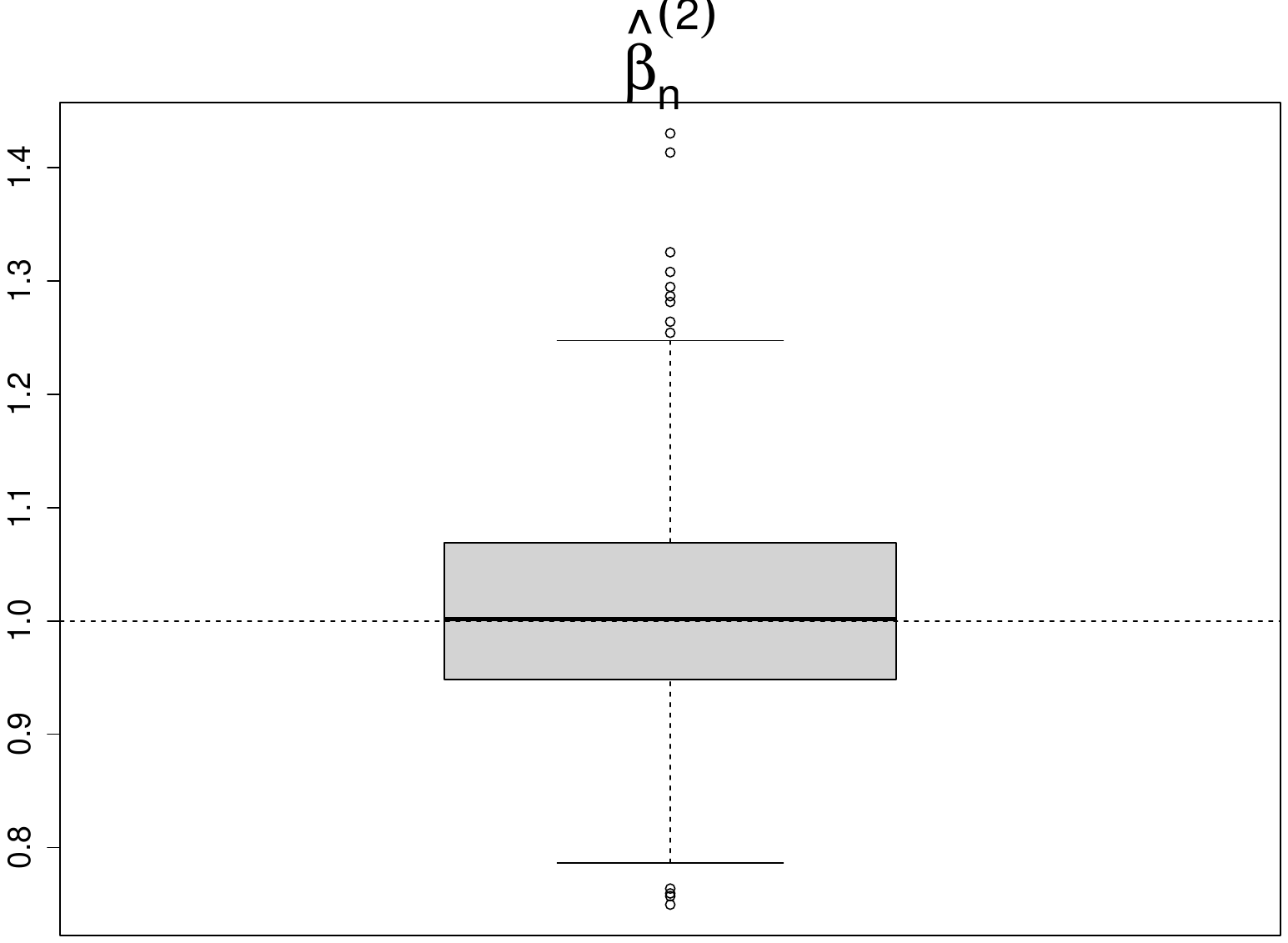}
    \end{subfigure}
    ~ 
    \begin{subfigure}[b]{0.3\textwidth}
        \includegraphics[width=\textwidth]{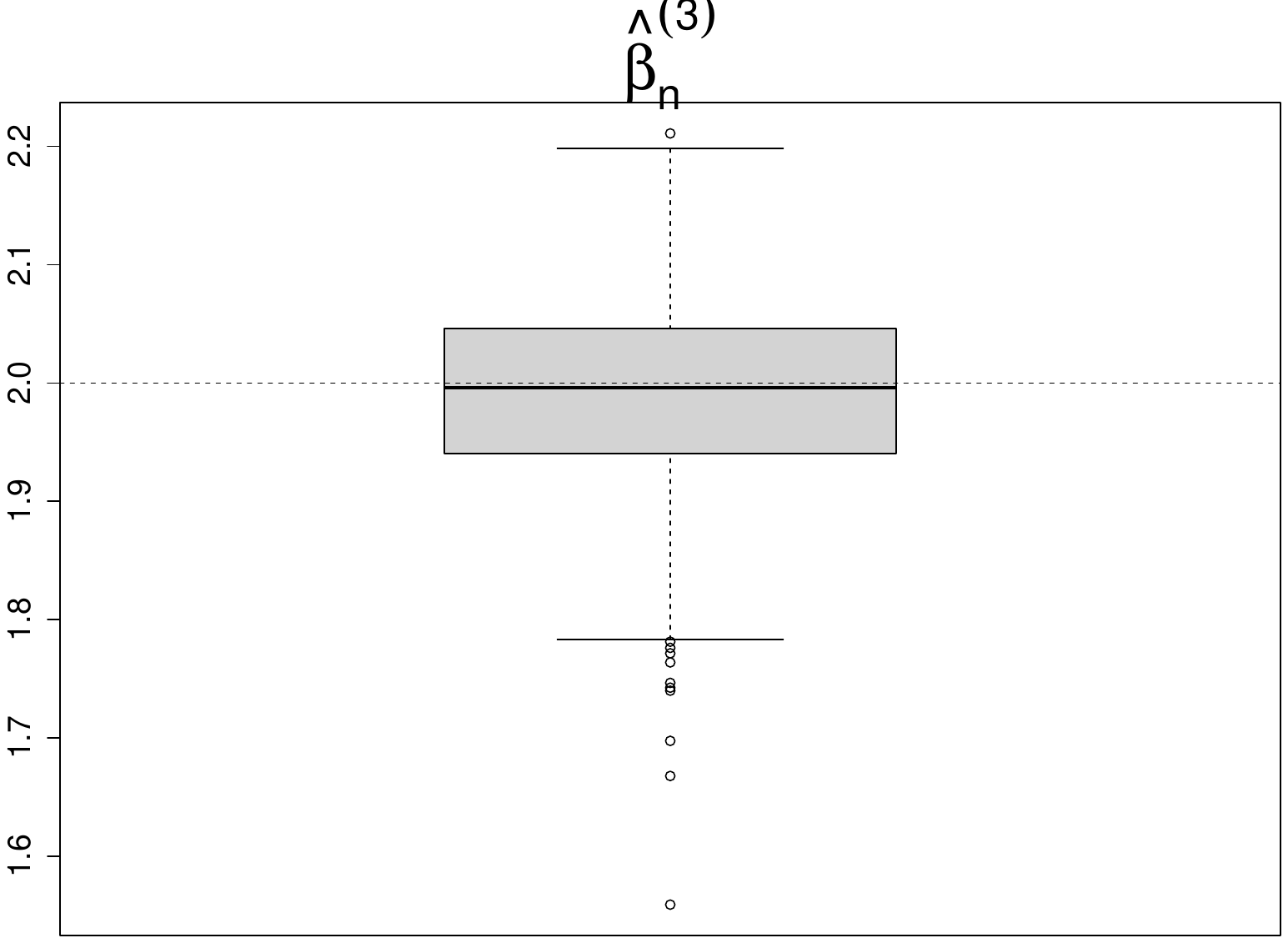}
    \end{subfigure}
    \caption{\label{fig:asympNormality} Boxplot of $\hat{\beta}_n^{(1)}$, $\hat{\beta}_n^{(2)}$, and $\hat{\beta}_n^{(3)}$ with $n=8000$ for 1000 iterations. The standard deviations are $0.02$, $0.09$, and $0.07$ respectively.}
\end{figure}

This example suggests that in cases where we do not have uniqueness, looking at the ordered version of $\hat{\beta}_n$ can provide additional information. Such a piece of information can be very valuable in the presence of a small matched data set or expert knowledge. In this case, one might even be able to recover the permutation that maps $\hat{\beta}_{n,\text{ordered}}$ to a consistent estimator of $\beta_0$.
\end{example}

\begin{example}\label{3dim}
For this example we generated 100 samples of size $n = 10^4$ from model~\eqref{unlinkedData} such that $Y\stackrel{d}{=} \beta_0^T X + \epsilon$ and 
\begin{eqnarray*}
    X\sim\mathcal{N}(0, I_3), \qquad \epsilon\sim\mathcal{N}(0, 1), \qquad\beta_0^T = [1, -1, 2].
\end{eqnarray*}
We also generate 100 samples of matched data $\{(\tilde{Y}_i, \tilde{X}_i)\}_{i=1}^m$ from model ~\eqref{matchedData} for $m\in\{10, 20, \cdots, 100\}$. Note that $\mathcal{B}_0 = \{\beta\in\RR^3: \Vert\beta\Vert_2 = \sqrt{6}\}$. We use the unlinked data for computing $\hat{\beta}_n$ and the linked data for computing the OLS estimator of $\beta_0$, $\tilde{\beta}_m$. In Table~\ref{tab:ols}, the mean and standard deviation of these estimators over the 100 replications are shown. It can be seen that $\Vert\hat{\beta}_n\Vert_2^2$ is highly concentrated around $\Vert\beta_0\Vert_2^2$ compared to the OLS estimators with smaller sample sizes.
\begin{table}[H]
\centering
\small{
\begin{tabular}{rrrrrrrrrrrr}
  \hline
  \multicolumn{1}{c}{} & \multicolumn{10}{c}{$\Vert\tilde{\beta}_m\Vert_2^2$, square norm of OLS estimate with sample size $m$} & \multicolumn{1}{c}{$\Vert\hat{\beta}_n\Vert_2^2$} \\
  \hline
m & 10 & 20 & 30 & 40 & 50 & 60 & 70 & 80 & 90 & 100 &  \\ 
mean & 6.44 & 6.23 & 6.07 & 6.02 & 6.12 & 6.04 & 6.11 & 6.25 & 6.02 & 6.11 & 6.00 \\ 
  sd & 2.03 & 1.23 & 0.87 & 0.80 & 0.66 & 0.65 & 0.58 & 0.56 & 0.52 & 0.55 & 0.16 \\ 
   \hline
\end{tabular}
}
\caption{\label{tab:ols} Mean and standard deviation of $\Vert\tilde{\beta}_m\Vert_2^2$ for values of $m\in\{10, 20, \cdots, 100\}$ for 100 independent iterations. The last column contains the mean and standard deviation of $\Vert\hat{\beta}_n\Vert_2^2$.}
\end{table}
\end{example}

\subsection{Intergenerational Mobility in the United States, 1850-1930}
The degree to which the economic status is passed along generations is an important factor in quantifying the inequality in a society. Researchers have approached this problem by studying the relationship between the father/father-in-law's income and the son/son-in-law's income. It is intuitive to assume that there is an increasing relationship between the son's income and that of the father. If we assume that this relationship is linear, the magnitude of the coefficient can quantify the existing inequality: the stronger the relationship, the more inequality. For this purpose,  we apply our method to data from 1850 to 1930 decennial censuses of the United States studied in~\cite{olivetti2015name, d2022partially} using the 1 percent IPUMS samples~\citep{ruggles2010integrated}. We follow \cite{olivetti2015name} and focus only on white father/father-in-law and son/son-in-law relationships. In this available historical data on father-son income in the United States, the link between the father/father-in-law and son/son-in-law is not available. Other studies on this data have used information on the first name to reconstruct the link between the father/father-in-law and son/son-in-law, but we only use the unmatched data and do not look at these partially reconstructed links. 

Since the exact value of the income in this data set is not available, we use the provided OCCSCORE. OCCSCORE assigns each occupation in all years a value representing the median total income (in hundreds of 1950 dollars) of all persons with that particular occupation in 1950. Therefore, using this score, we lose the within-occupation variation of the income. 

\begin{figure}[H]
  \centering
    \includegraphics[width=1\textwidth]{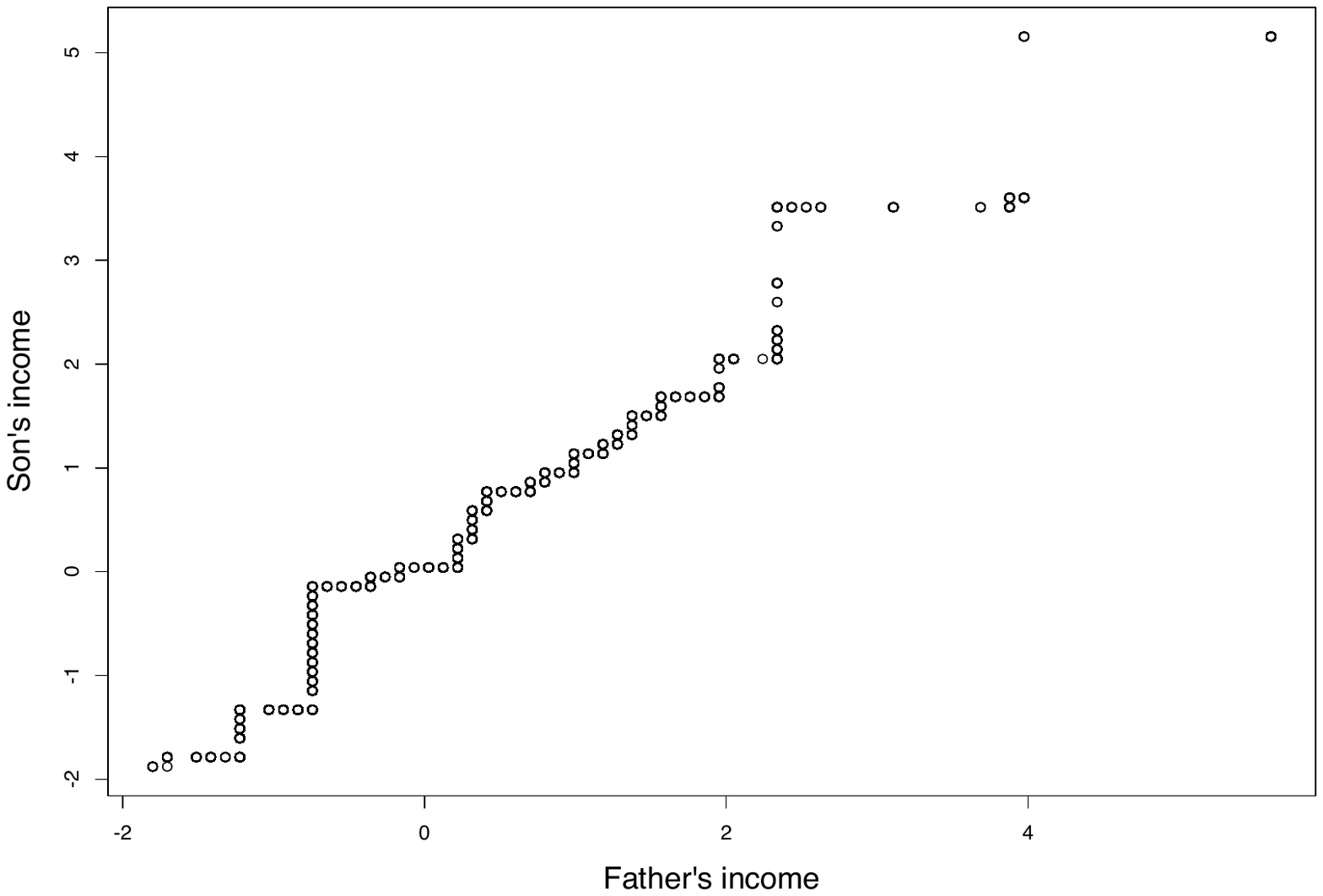}
    \caption{\label{fig:qqFS} The sorted values of $Y$ (Son's transformed income) are plotted against the sorted values of $X$ (Father's transformed income) for a sub-sample of size $10^6$.}
\end{figure}

We let $Y$ be the son/son-in-law's OCCSCORE and $X$ to be the father/father-in-law's OCCSCORE and assume that $Y \stackrel{d}{=} \beta_0 X + \epsilon$. The sample sizes in the data sets are quite large ($n_X, n_Y > 10^7$) and therefore, for computational reasons, we select a subset of size $4000$ of $X$ and $Y$ at random. We centered and normalized each sub-sample such that they have a mean of 0 and a variance of 1 but kept the same names ($Y$ and $X$) for the transformed variables. In Figure~\ref{fig:qqFS}, we show the sorted values of $Y$ (Son's income) plotted against those of $X$ (Father's income) for a sub-sample of size $10^6$.

Since the true distribution of the noise variable $\epsilon$ is unknown in this problem, we consider two different families of centred distributions,  Normal and Laplace. We consider different possible values for their scale parameters, so that the standard deviation (sd) of the noise varies in the set $\{0.1, 0.2, \cdots, 1\}$. Recall that the standard deviation of a random variable distributed as $\mathcal{N}(0, \sigma^2)$ is $\sigma$ and as $\text{Laplace}(\lambda)$ is $\sqrt{2}\lambda$. Therefore, we choose the parameters $\sigma$ and $\lambda$ accordingly. 
Note that the data is discrete as there are many repeated values for each of $X$ and $Y$.
\begin{figure}[H]
  \centering
    \includegraphics[width=1\textwidth]{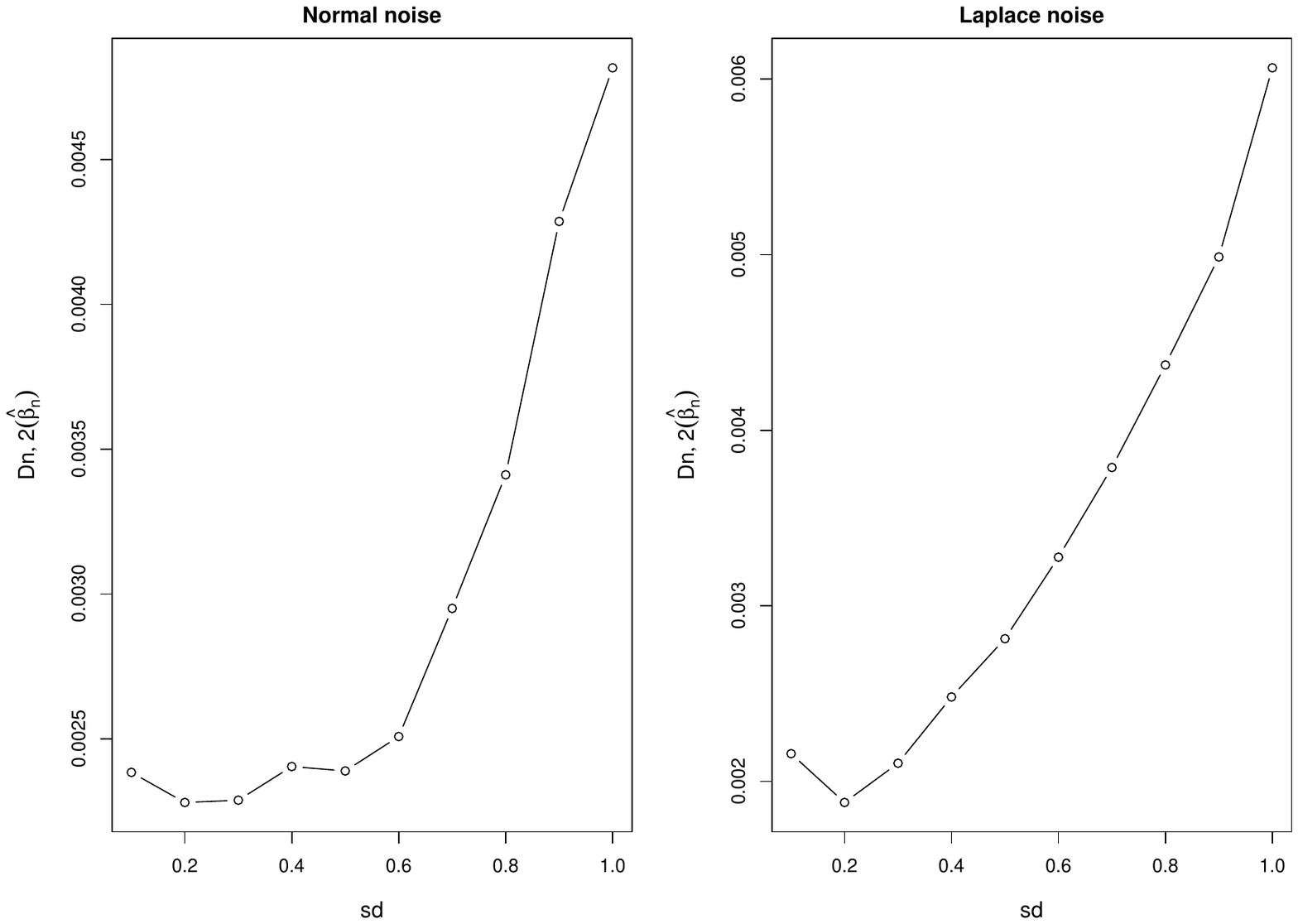}
    \caption{\label{fig:normalLaplace} Values of $\mathbb{D}_{n, 2}(\hat{\beta}_n)$ for different values of standard deviation (sd) when the noise distribution is assumed to be $\mathcal{N}(0, \sigma^2)$ (sd = $\sigma$), and $\text{Laplace}(\lambda)$ (sd = $\sqrt{2}\lambda$).}
\end{figure}
Varying $\sigma$ and $\lambda$, the respective parameters of the distributions $\mathcal{N}(0, \sigma^2)$ and $\text{Laplace}(\lambda)$, such that the standard deviation takes its value in the set $\{0.2, 0.3, \cdots, 1\}$ allows us to plot the values of $\mathbb{D}_{n, 2}(\hat{\beta}_n)$ in Figure~\ref{fig:normalLaplace}. Assuming that the noise is Gaussian, we observe that $\mathbb{D}_{n, 2}(\hat{\beta}_n)$ takes its minimum over the grid of values of $\sigma\in\{0.1, 0.2, \cdots, 1\}$ at $\sigma = 0.2$ and the optimization procedure yields $\hat{\beta}_n = 1.27$. Similarly assuming that the noise follows $\text{Laplace}(\lambda)$ with parameter $\lambda$ such that $\sqrt{2}\lambda\in\{0.1, 0.2, \cdots, 1\}$ we find the minimum to be attained for $\lambda = 0.2/\sqrt{2}$ and  $\hat{\beta}_n = 1.18$. 

Relying on the observation that $\mathbb{D}_{n, 2}(\hat{\beta}_n)$ is minimized for $\text{sd} = 0.2$ for both the Normal and Laplace distributions over 100 iterations we take independent sub-samples of size $n = 4000$ and calculate $\hat{\beta}_n$ using Normal and Laplace distribution for the noise with standard deviation equal $0.2$. For both distributions, the estimate $\hat{\beta}_n$ shows a bi-modal behaviour with positive and negative modes with the same magnitude. This is to be expected in case the covariate has a zero expectation. However, as we know that the relationship between $Y$ and $X$ should be non-decreasing, we can consider the absolute value of $\hat{\beta}_n$ in all cases. Table ~\ref{tab:mean_sd_abs} shows the values of mean and standard deviation for $|\hat{\beta}_n|$ using the Normal and Laplace distributions for the noise. Note that Theorem~\ref{limitDist} guarantees asymptotic normality of $\hat{\beta}_n$ when $|\mathcal{B}_0| = 1$ under some regularity assumptions. Such assumptions do not seem to be fulfilled in this dataset because of its discreteness. However, one may use sub-sampling ideas to create confidence intervals based on $\hat{\beta}_n$, hoping that some asymptotic normality holds. 

\begin{table}[H]
\centering
\small{
\begin{tabular}{rrr}
  \hline
  $|\hat{\beta}_n|$ & Normal (sd = 0.2) & Laplace (sd = 0.2)\\
  \hline
mean & 1.17 & 1.14 \\ 
  sd & 0.06 & 0.09 \\ 
  $95\%$ Bootstrap confidence interval & $(1.07, 1.28)$ & $(1.00, 1.26)$ \\
   \hline
\end{tabular}
}
\caption{\label{tab:mean_sd_abs} Mean and standard deviation for $|\hat{\beta}_n|$ using the Normal and Laplace as the noise distribution over 100 sub-samples of size $n = 4000$. The last row of the table shows the  $95\%$ Bootstrap confidence interval using the 100 sub-samples.}
\end{table}
In this example, the covariate $X$ is 1-dimensional, and the distribution of data is far from continuous. Still, we found this example interesting as it is the only real data problem with unmatched data that we could access to.
\subsection{Power Plant Data Set}
We consider the Power Plant data set from UCI Machine Learning Repository\footnote{\url{https://archive.ics.uci.edu/}}. The data set contains 9568 matched data points collected from a combined cycle power plant. Features consist of hourly average ambient variables Temperature (T), Ambient Pressure (AP), Relative Humidity (RH) and Exhaust Vacuum (V) to predict the net hourly electrical energy output (EP) of the plant. Assuming that EP is a linear function of T, AP, RH, and V, we run ordinary least squares using all the data points. We get $R^2 = 0.93$, which supports the linearity of the relationship. We perform 100 independent simulations in the following manner: For each simulation, we select a sub-sample of matched data of size $m = 30$. We use this subsample of matched data for both $\tilde{\beta}_m$ and estimate the distribution of $\epsilon$. Then, from the remainder of the data, we select a sub-sample of unmatched data of size $n = 4000$. This is done by selecting a sample of size $n=4000$ from the remaining $Y_i$'s and independently selecting a sample of size $n=4000$ from the remaining $X_i$'s. Since we do not have access to the population density in this case, we take the OLS estimate $\beta_{0, \text{OLS}}$ using all the data points as the ground truth. 

For the DLSE estimator $\hat{\beta}_n$, we need an estimate of the noise distribution, and for this, we use a Kernel density estimator based on the residuals of OLS $\tilde{\beta}_m$ obtained using the matched data. We use a Gaussian kernel and select the bandwidth according to~\cite{SheatherJones}.\footnote{We use function \textit{density} from R package ``stats" with hyper-parameter ``SJ" for the bandwidth and Gaussian kernel.} 

\begin{figure}[H]
  \centering
    \includegraphics[width=1\textwidth]{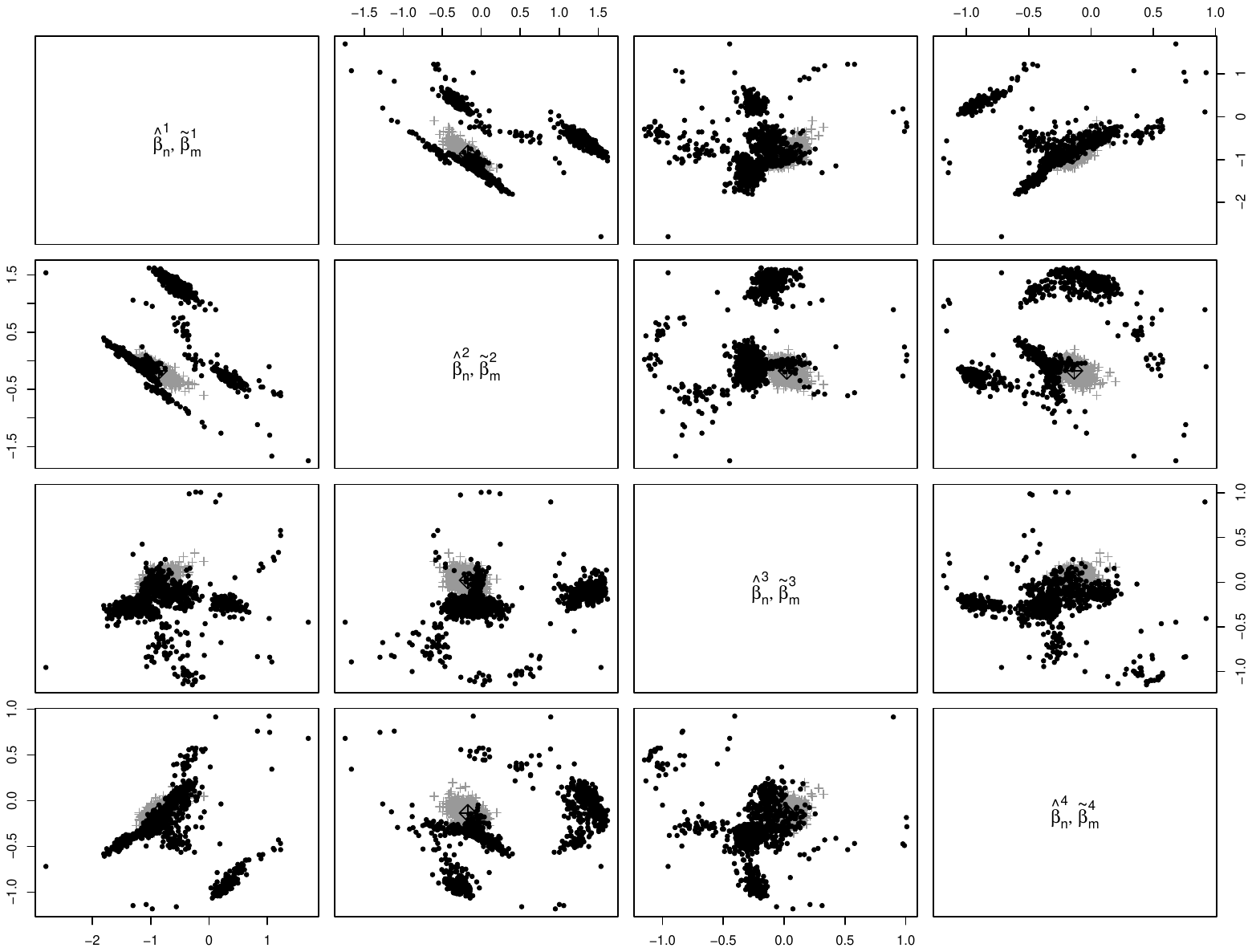}
    \caption{\label{fig:uniquenessBeta} Scatterplots of pairs of components of $\hat{\beta}_n$ (filled black point) and $\tilde{\beta}_m$ (gray plus) together with projections of $\beta_{0, \text{OLS}}$ (black diamond). Note that in each plot, one of the clusters of $\hat{\beta}_n$ is concentrated near the projection of $\beta_{0, \text{OLS}}$. In a real-world scenario, we do not have access to $\beta_{0, \text{OLS}}$ to choose the corresponding sub-clusters of $\hat{\beta}_n$. For this purpose one can use $\tilde{\beta}_m$ instead.} 
\end{figure}

Figure~\ref{fig:uniquenessBeta} depicts the scatterplots of pairs of components of the obtained DLSE $\hat{\beta}_n$ over 1000 sub-samples of size $n = 4000$ together with the projection of $\beta_{0, \text{OLS}}$ onto the corresponding sub-spaces. The grey points in the scatterplots correspond to components of $\tilde{\beta}_m$. For example, the plot in the first row and the second column is the scatter plot of $\hat{\beta}_n^1$ vs $\hat{\beta}_n^2$ (and $\tilde{\beta}_m^1$ vs $\tilde{\beta}_m^2$) which shows that the first two components of $\hat{\beta}_n$ form 3 main clusters. Therefore, the scatterplots suggest that $\hat{\beta}_n$ is not converging to a unique value. On the other hand, they show that the vectors $\hat{\beta}_n$'s are concentrated around multiple modes similar to the phenomena described in Theorem~\ref{pole}. Note that in each scatterplot, one of the clusters can be represented by the projection of $\beta_{0, \text{OLS}}$ and therefore can be recognized as the \lq\lq true cluster\rq\rq. As in reality we do not have access to the whole matched data set, one can use $\tilde{\beta}_m$ as a \lq\lq guide\rq\rq \ to pin down the right cluster.  We intend to formalize this idea more concretely in the scope of future work.

\section{Discussion and Future Research Directions}
In this paper, we have proposed an approach for making inference in a regression model in the case where the link between responses and covariates is not known. Here, our main goal is to estimate the regression vector. Note that the model we consider includes the permuted regression model studied by in~\cite{Pananjady:2018hd}, where the principal focus is to recover the permutation under which the covariates are shuffled and which is responsible for the loss of the link. The main idea that we pursued in this work is to view the responses as random variables generated from the convolution of $\beta^T_0 X$ and $\epsilon$ and find the vector which minimizes over all possible $\beta\in\RR^d$ an $\ell_2$-distance between the sample distribution of the responses and a natural estimator of the distribution of the convolution of $\beta^T X$ and $\epsilon$. As it becomes clear at the end of Section~\ref{MainResult}, the problem boils down to a deconvolution. It may come as very surprising that one is able to recover, under some identifiability conditions, the true regression vector or an ordered version thereof at the $n^{-1/2}$-rate. In fact, the convergence rates in deconvolution problems are known to be slow to very slow for smooth noise distributions. For example, if the noise is Gaussian and under minimal smoothness assumptions, then it follows from the seminal paper of~\cite{Fan91} that $(\log n)^{-1/2}$ is the optimal rate of convergence for deconvoluting the distribution of $X$ from the noise $\epsilon$ based on observations from the model $Y = X + \epsilon$. How do we reconcile this slow rate with the obtained $n^{-1/2}$ in our Theorem~\ref{limitDist}? The answer is that our problem is completely parametric, and the distribution to be deconvoluted is in a much smaller class than the one considered in~\cite{Fan91}.   

The estimation method that is chosen and implemented in this work is partially motivated by convenience to some degree. In fact, the arguments from the theory of empirical process, although technical, seem to be less cumbersome as they are, for example, for the Maximum Likelihood estimator (MLE). Finding the MLE in this problem would have been a very appealing approach, and the estimator might be even more efficient. In this case, one would find a regression vector which maximizes the log-likelihood
\begin{eqnarray*}
\beta \mapsto   \sum_{j=1}^n \log\left( n^{-1} \sum_{i=1}^n f^\epsilon(Y_j -  \beta^T X_i) \right)
\end{eqnarray*}
over $\RR^d$.  Optimization, in this case, presents harder numerical challenges because of the logarithm.  Nevertheless,  we believe that this estimator should be implemented. We leave this task to future research work, where we also plan to study the asymptotic properties of this estimator and compare its performance to the DLSE studied here. If the noise distribution is assumed to be known, we would like to stress the fact that it is possible to relax such an assumption. This can be done either by (1) estimating this distribution from a matched sample as done above for the Power Plant data set or (2) assuming that it belongs to some scaled parametric family and estimating the scale parameter, $\sigma_0$ say, together with the unknown regression vector.  This approach was followed for the inter-generational mobility data set. We conjecture that under identifiability, the resulting estimators of $\beta_0$ and $\sigma_0$ are asymptotically Gaussian.  A referee noted that it is possible to consider the case where the response $Y$ is multivariate, of dimension $p$, in which case one needs to estimate a whole regression matrix.  If the components of $Y$ are independent, then the asymptotic theory developed here can be easily extended. However, the optimization problem, already difficult for $p=1$, will impose additional numerical challenges.  In case the components are not independent, we do not think that extending our estimator to this case would be straightforward. Computation of the empirical distribution of the observed responses becomes very tedious as well as the criterion, which now involves multivariate integration.

The great potential of the method described in this work was clearly seen in the semi-supervised learning situation, where some matched data is available. We view this scenario as the best and most realistic application for two reasons:  (1)  One is able to estimate the distribution of the noise from the matched proportion and not simply impose it. (2)  The OLS found with the matched data can be used in combination with re-sampling from the unmatched data to \lq\lq guide\rq\rq \ the deconvolution least square method in finding the closest DLSE to the OLS.  The point in (2) is particularly useful in case the re-sampled data do not seem to agree on the DLSE.

Finally, we would like to point out that the linear regression model can be, of course, extended to other models, such as logistic regression for example. An extension is the high-dimensional setting, although we believe that many theoretical challenges will have to be tackled, both from theoretical and numerical perspectives. 
\acks{}
The authors thank Professor Charles H. Doss for very interesting discussions on the methodology used in the paper. The authors thank Professor Wendelin Werner for some very useful hints, which allowed us to finish the proof of Theorem 7, Professor Ashwin Pananjady and Professor Martin Wainwright for helpful comments.  

\newpage
\appendix
\section{Supplementary Material}\label{proofs}
We present proofs of all the theorems and lemmas in the following. 
\subsection{Proof of Proposition \ref{Existence}}
\begin{proof}
Our goal is to show that with probability one, there exists a vector $\beta\in\RR^d$ such that it minimizes $\DD_{n, p}$. 

For $K > 0$ consider the compact set $\{ \beta \in \RR^d: \Vert \beta \Vert \le K \}$. We show that our choice of $K$ depends only on $p$. Let $u\in S^{d-1}$ be a fixed vector. Then we have
\begin{eqnarray*}\label{eq:3}
\lim_{\lambda \to \infty} \DD_{n,p}(\lambda u) &=& \frac{1}{n^{p+1}}\sum_{j=1}^n \big \vert j - n_0 F^\epsilon(Y_{(j)}) - n_{-}\big\vert^p
\end{eqnarray*}
with $n_{-} = |\{i: u^T X_i < 0 \}|$ and $n_0 = |\{i: u^T X_i = 0\}|$. Note that $n_{-}$ and $n_0$ depend on the vector $u$. Since $X$ has a continuous distribution, $n_0\leq d$ with probability one. If $n_0\geq d+1$, then for some distinct $i_1, \cdots, i_{d+1}$, $X_{i_1}, \ldots, X_{i_{d+1}}$ are not linearly independent, which is of probability zero.

Suppose that $n_{-} \ge n/2$. Using $\vert \vert a \vert - \vert b \vert \vert \le \vert a - b \vert$, it follows that  
\begin{eqnarray*}
\sum_{j=1}^n  \big \vert j - n_0 F^\epsilon(Y_{(j)})  - n_{-}  \big \vert^p   \ge \sum_{j \le n/4}  \Big \vert   \vert j -  n_{-} \vert   - n_0 F^\epsilon(Y_{(j)})  \Big \vert^p.
\end{eqnarray*}
Since $j \le n/4$ we have $\vert j -  n_{-} \vert  \ge n/4$. By taking $n > 8d$, and following from $n_0 F^\epsilon(Y_{(j)})  \le n_0 \le d$, with probability one we have
\begin{eqnarray*}
\sum_{j=1}^n \big\vert j - n_0 F^\epsilon(Y_{(j)}) - n_{-} \big\vert^p &\ge& \sum_{j \le n/4}\Big(\vert j - n_{-}\vert - d\Big)^p\\
&\ge&\frac{n}{4}\left(\frac{n}{4} - d\right)^p\\
&\ge&\frac{n}{4}\left(\frac{n}{8}\right)^p = \frac{n^{p+1}}{2^{3p+2}} > \frac{n^{p+1}}{8^{p+1}}.
\end{eqnarray*}
This implies that for $n \geq 8d$
\begin{eqnarray*}
\frac{1}{n^{p+1}}\sum_{j=1}^n\big\vert j - n_0 F^\epsilon(Y_{(j)}) - n_{-}\big\vert^p\ge \frac{1}{8^{p+1}}.
\end{eqnarray*}
Now, for any small $\eta > 0$ and large enough $\lambda^\prime$
\begin{eqnarray*}
\DD_{n,p}(\lambda^\prime u)\ge \lim_{\lambda \to \infty} \DD_{n,p}(\lambda u)  - \eta.
\end{eqnarray*}
Let $\eta =  2^{-(3p+4)}$. Then, there exists $K_p > 0$ such that for $\lambda > K_p$
\begin{eqnarray}\label{MainIneq0}
\DD_{n,p}(\lambda u) \ge \frac{1}{16 \cdot 8^{p}}
\end{eqnarray}
with probability $1$ provided that $n \ge 8d$. Note that for $n_{-} < n/2$, similar reasoning can be applied by summing over $j \ge n/4$. Now, consider the set
\begin{eqnarray}\label{Cp}
\mathcal{C}_{p} = \Big\{\beta\in\RR^d: \Vert \beta \Vert \le  K_p \Big \}. 
\end{eqnarray}
where $K_p$ is the same constant defined above. Since the lower bound shown in (\ref{MainIneq0}) does not depend on the vector $u$,  we  conclude that  for all $n \ge 8d$
\begin{eqnarray}\label{MainIneq}
\inf_{\beta\in\mathcal{C}^c_p}\mathbb D_{n, p}(\beta)\ge\frac{1}{16 \cdot 8^{p}}.
\end{eqnarray}
On the other hand, and by using the equality $F^Y =  F^{\beta^T_0 X} \star F^\epsilon$,
\begin{align*}\label{Dnp}
    \MoveEqLeft \DD_{n,p}(\beta_0) = \int \bigg \vert \Big(F^Y_n(y) - F^Y(y)\Big) - \Big( n^{-1}\sum_{i=1}^n F^\epsilon(y - \beta^T_0 X_i) - [F^{\beta^T_0 X} \star F^\epsilon](y)\Big)\bigg\vert^p d F_n(y) \\
    & \le 2^{p-1}\ \int   \Big \vert F^Y_n(y)  -  F^Y(y)\Big \vert^p d F_n(y)  \\
    & \quad + 2^{p-1} \  \int  \bigg \vert n^{-1}  \sum_{i=1}^n F^\epsilon(y - \beta^T_0 X_i)  - [F^{\beta^T_0 X} \star F^\epsilon](y)  \bigg \vert^p  d F_n(y).
\end{align*}
Convexity of $x \mapsto \vert x \vert ^{p}$ gives us
\begin{align*}
    \DD_{n,p}(\beta_0) &\le 2^{p-1}\Vert F^Y_n - F^Y \Vert_\infty + 2^{p-1} \Big \Vert F^{\beta^T_0 X}\star F^\epsilon - F^{\beta^T_0 X}_n \star F^\epsilon \Big \Vert_\infty \\
    & \le 2^{p-1} \Vert F^Y_n - F^Y\Vert_\infty + 2^{p-1}\Big\Vert F^{\beta^T_0 X} -  F^{\beta^T_0 X}_n \Big \Vert_\infty.
\end{align*}
where $\Vert \cdot \Vert_\infty$ denotes the supremum norm, and $F^{\beta^T_0X}_n$ is the empirical distribution based on $\beta^T_0 X_i, i=1, \cdots, n$. By Glivenko-Cantelli we have
\begin{eqnarray}\label{MainConv}
 \PP\left(\omega: \exists \ n_p (\omega) \ \forall  \ n \ge n_p(\omega)  \  \DD_{n,p}(\beta_0) \le  \frac{1}{32 \cdot 8^p} \right) = 1.
\end{eqnarray}
Therefore by combining \eqref{MainIneq} and \eqref{MainConv} we have
\begin{eqnarray*}
&&\PP \Big(\omega: \exists \ n_p (\omega) \ \forall  \ n \ge \max(n_p(\omega), 8d)  \  \  \DD_{n, p} \ \textrm{admits a minimizer in}  \ \mathcal{C}^c_{p} \Big)\\
&& =   \PP \Big(\omega: \exists \ n_p (\omega) \ \forall  \ n \ge \max(n_p(\omega), 8d)  \  \  \exists \beta \in \mathcal{C}^c_p : \DD_{n, p}(\beta_0)  \ge \DD_{n, p}(\beta) \ge \frac{1}{16 \cdot 8^p}\Big) \\
&& = 0.
\end{eqnarray*}
Since $\mathcal{C}_p$ defined in \eqref{Cp} is compact and $\DD_{n, p}$ is Lebesgue-a.e continuous we have 
\begin{eqnarray*}
&& \PP \Big(\omega: \exists \ n_p (\omega) \ \forall  \ n \ge \max(n_p(\omega), 8d)  \  \  \mathbb D_{n, p} \ \textrm{admits a minimizer in}  \ \mathbb{R}^p \Big)\\
&&= \PP \Big(\omega: \exists \ n_p (\omega) \ \forall  \ n \ge \max(n_p(\omega), 8d)  \  \  \mathbb D_{n, p} \ \textrm{admits a minimizer in}  \ \mathcal{C}_p \Big) \\
&& = 1.
\end{eqnarray*}
\hfill 
\end{proof}

\subsection{A Useful Proposition and Its Proof}
We will use the following proposition in the rest of the proofs. 
\begin{proposition}\label{conv1}
For any $\hat{\beta}_n\in\mathcal{B}_{n, 2}$ we have
\begin{eqnarray*}
 \int_{\RR} \left( F_n^Y(y) - n^{-1} \sum_{i=1}^n F^\epsilon(y - \hat{\beta}^T_n X_i)   \right)^2 d F^Y_n(y) = O_{\PP}(n^{-1}).
\end{eqnarray*}
\end{proposition}

\begin{proof}
By definition of $\hat{\beta}_n$ we have 
\begin{eqnarray*}
\MoveEqLeft \int_{\RR} \left(F^Y_n(y) - n^{-1} \sum_{i=1}^n F^\epsilon(y - \hat{\beta}^T_n X_i)\right)^2 d F^Y_n(y)\le \\
&\int_{\RR} \left(F^Y_n(y) - n^{-1} \sum_{i=1}^n F^\epsilon(y - \beta^T_0 X_i)\right)^2  d F^Y_n(y)
\end{eqnarray*}
By adding and subtracting $F^Y(y)$ we have
\begin{eqnarray*}
\DD_{n, 2}(\beta_0)& = &\int_{\RR} \left(F^Y_n(y) - n^{-1} \sum_{i=1}^n F^\epsilon(y - \beta^T_0 X_i)   \right)^2  d F^Y_n(y)  \\
&\le & 2 \int_{\RR} \left(F^Y_n(y) - F^Y(y)\right)^2 d F^Y_n(y) + \\
&& 2 \int \left(F^Y(y) - n^{-1} \sum_{i=1}^n F^\epsilon(y - \beta^T_0 X_i)\right)^2d F^Y_n(y).
\end{eqnarray*}
By Dvoretzky–Kiefer–Wolfowitz inequality;~\cite{DKW} we have  
\begin{eqnarray*}
\int_{\RR} \left(F^Y_n(y) - F^Y(y)\right)^2 d F^Y_n(y) \le \Vert F^Y_n - F^Y \Vert^2_\infty  = O_{\PP}(n^{-1}).
\end{eqnarray*}
Using $F^Y(y) = \int F^\epsilon(y - \beta^T_0 x) d F^X(x)$ it follows that
\begin{eqnarray*}
&&\int \left(F^Y(y) - n^{-1} \sum_{i=1}^n F^\epsilon(y - \beta^T_0 X_i)\right)^2d F^Y_n(y)\\
&& = \int \left( \int F^\epsilon(y - \beta^T_0 x) d F^X(x) - \int F^\epsilon(y - \beta^T_0 x) d F^X_n(x) \right)^2 d F^Y_n(y)\\
&& = \int \left( \int F^\epsilon(y - z) d F^Z(z) - \int F^\epsilon(y - z) d F^Z_n(z) \right)^2 d F^Y_n(y), \  \ \text{with $Z = \beta^T_0 X$}  \\
&& =  \int\left(\int\left(F^Z(z) - F^Z_n(z)\right) f^\epsilon(y - z) dz\right)^2 d F^Y_n(y) \\
&& \le \Vert F^Z_n - F^Z \Vert^2_\infty = O_{\PP}(n^{-1}).
\end{eqnarray*}
Therefore $\DD_{n, 2}(\beta_0) = O_{\PP}(n^{-1})$. Since by definition $\DD_{n, 2}(\hat{\beta}_n)\leq\DD_{n, 2}(\beta_0)$ we conclude that
\begin{eqnarray}\label{RateDnbetahat}
\mathbb{D}_{n,2}(\hat{\beta}_n) = \int_{\RR} \left(F^Y_n(y) - n^{-1} \sum_{i=1}^n F^\epsilon(y - \hat{\beta}^T_n X_i)  \right)^2  d F^Y_n(y)  = O_{\PP}(n^{-1})
\end{eqnarray}
which implies that 
\begin{eqnarray*}
&& \int_{\RR} \left(F^Y(y) - n^{-1} \sum_{i=1}^n F^\epsilon(y - \hat{\beta}^T_n X_i)  \right)^2  d F^Y_n(y)  \\
&& \le 2 \int_{\RR} \left(F^Y(y)  - F^Y_n(y) \right)^2 d F^Y_n(y) + \\
&& 2 \int_{\RR}  \left(F^Y_n(y)  -  n^{-1} \sum_{i=1}^n F^\epsilon(y - \hat{\beta}^T_n X_i) \right)^2 d F^Y_n(y) \\
&& = O_{\PP}(n^{-1})
\end{eqnarray*}
which completes the proof. \hfill 
\end{proof}

\subsection{Proof of Theorem \ref{convNorm}}
\begin{proof}
For positive semi-definite matrix $\Gamma$ and $z \in \RR^2$, define $\Vert z \Vert_{2, \Gamma} := \sqrt{z^T \Gamma z}$. Note that $\Vert \hat{\beta} \Vert_2 = O_{\PP}(1)$ is equivalent to $ \Vert \hat{\beta}_n \Vert_{2, \Gamma} = O_{\PP}(1)$. It follows from~\cite[Theorem 5.14]{van1998asymptotic} that for any compact set $\mathcal{K}\subset\RR^d$ and $\epsilon > 0$
\begin{eqnarray*}
\lim_{n\to\infty}\PP\left(\inf_{\beta\in\mathcal{B}_0}\Vert\hat{\beta}_n - \beta \Vert_{2, \Gamma}\ge\epsilon,\hat{\beta}_n \in\mathcal{K}\right)=0.
\end{eqnarray*}
Assuming that $\Vert \hat{\beta}_n \Vert_{2, \Gamma} = O_{\PP}(1)$, for any $\eta\in (0,1)$ we can find $K_\eta > 0$ such that 
\begin{eqnarray*}
\PP(\Vert \hat{\beta}_n \Vert_{2, \Gamma} \le K_\eta) \ge 1 - \eta
\end{eqnarray*}
for all $n$. Take $\mathcal{K}$ to be the closed Euclidean ball in $\RR^d$ with centre $0$ and radius $K_\eta$. Then, 
\begin{eqnarray*}
\PP\left(\inf_{\beta\in\mathcal{B}_0}\Vert\hat{\beta}_n - \beta \Vert_{2,\Gamma}  \ge \epsilon \right) 
& = & \PP\left(\inf_{\beta\in\mathcal{B}_0}\Vert\hat{\beta}_n - \beta \Vert_{2, \Gamma}  \ge \epsilon, \hat{\beta}_n \in \mathcal{K} \right) + \\
&& \PP\left(\inf_{\beta \in \mathcal{B}_0}\Vert\hat{\beta}_n - \beta \Vert_{2, \Gamma}  \ge \epsilon, \hat{\beta}_n \notin \mathcal{K} \right) \\
& \le & \PP\left( \inf_{\beta \in  \mathcal{B}_0}  \Vert \hat{\beta}_n - \beta \Vert_{2, \Gamma}  \ge \epsilon, \hat{\beta}_n \in \mathcal{K} \right)  + \eta.
\end{eqnarray*}
Hence
\begin{eqnarray*}
\lim_{n \to \infty}\PP\left(\inf_{\beta\in\mathcal{B}_0}\Vert\hat{\beta}_n - \beta \Vert_{2, \Gamma} \ge \epsilon \right) \le \eta.
\end{eqnarray*}
Since $\eta > 0$ can be arbitrarily small, this implies
\begin{eqnarray*}
\lim_{n \to \infty} \PP\left(\inf_{\beta\in\mathcal{B}_0}\Vert\hat{\beta}_n - \beta \Vert_{2, \Gamma} \ge \epsilon \right) = 0.
\end{eqnarray*}
By the assumption of Theorem for any $\beta \in \mathcal{B}_0$ we have $\Vert \beta \Vert_{2, \Gamma} = c$. This, together with triangle inequality, implies
\begin{eqnarray*}
\vert \Vert \hat{\beta}_n \Vert_{2, \Gamma}  -  \Vert \beta \Vert_{2, \Gamma} \vert = \vert \Vert \hat{\beta}_n \Vert_{2, \Gamma}  -  c  \vert \le \Vert \hat{\beta}_n -   \beta  \Vert_{2, \Gamma},
\end{eqnarray*}
hence  
\begin{eqnarray*}
\Vert \hat{\beta}_n \Vert_{2, \Gamma} = \sqrt{\hat{\beta}^T_n \Gamma \hat{\beta}_n}  \stackrel{\PP}{\rightarrow} c
\end{eqnarray*}
as $n \to \infty$. This completes the proof. \hfill  
\end{proof}
\subsection{Proof of Lemma \ref{Charbeta}}
\begin{proof}
For a random variable $Z$ denote the characteristic function of $Z$ by $\Phi_Z$; i.e.,
\begin{eqnarray*}
\Phi_Z(t) = \EE[\exp(it Z)], t  \in \RR.  
\end{eqnarray*}
Then, $\Phi_Y(t) = \Phi_{\beta^T X}(t) \Phi_{\epsilon}(t)$ for all $t \in \RR$. Let $\mathcal{V}\subset\RR$ be a small interval containing 0 such that for all $t\in\mathcal{V}$ we have $\Phi_{\epsilon}(t) \ne 0$. Then $\beta_1, \beta_2\in\mathcal{B}_0$ and $\beta_2$ if and only if $\Phi_{\beta^T_1 X}(t) = \Phi_{\beta^T_2 X}(t)$ for all $t \in \mathcal{V}$. Since $\beta^T X \sim \mathcal{N}(0, \beta^T \Sigma \beta)$ we have
\begin{eqnarray*}
\exp\left(-\frac{t^2}{2}\beta^T_1\Sigma\beta_2 \right) = \exp\left(-\frac{t^2}{2}  \beta^T_2 \Sigma \beta_2 \right).
\end{eqnarray*}
Therefore $\beta^T_1 \Sigma \beta_1  = \beta^T_2 \Sigma \beta_2$. Note that $\Sigma$ is positive definite which implies $c = \beta^T \Sigma \beta > 0$ for $\beta\in\mathcal{B}_0$. \hfill 
\end{proof}
\subsection{Proof of Proposition \ref{Gaussian}}
\begin{proof}
Let us start with the case where $\Sigma = I_d$. Using (\ref{boundcngamma}) we have 
\begin{eqnarray*}
\small{
\lim_{n \to \infty} \PP\Big( \forall \ \beta \in \RR^d, \  \textrm{and}  \ j=1,.., n, \frac{1}{n} \sum_{i=1}^n F^\epsilon(Y_j - \beta^T X_i) \ge \int F^\epsilon(Y_j- \beta^T x) d F^X(x) - \frac{1}{16} \Big) = 1
}
\end{eqnarray*}
where
\begin{eqnarray*}
\int F^\epsilon(Y_j - \beta^T x) d F^X(x) = \int F^\epsilon(Y_j - \Vert \beta \Vert_2 z) \varphi(z) dz, \  \ \text{where $\varphi$ is the density of $\mathcal{N}(0,1)$}.
\end{eqnarray*}
This implies that
{\footnotesize{
\begin{eqnarray*}
\lim_{n \to \infty}\PP\Big(\forall \ \beta \in \RR^d, \  \text{and}  \ j=\lfloor n/5 \rfloor, \cdots, \lfloor n/4 \rfloor, \frac{1}{n} \sum_{i=1}^n F^\epsilon(Y_{(j)} - \beta^T X_i)  \ge \int F^\epsilon(Y_{(j)}- \beta^T x) d F^X(x)  - \frac{1}{16}\Big) =1.
\end{eqnarray*}
}}
For any fixed $j$ and $\delta > 0$
\begin{eqnarray*}
\int F^\epsilon(Y_{(j)}- \beta^T x) d F^X(x)  & = &  \int_{z: \ z < 0} + F^\epsilon(Y_{(j)} - \Vert \beta \Vert_2  z) \varphi(z) dz   \\
&& \int_{z: \ z > 0}  F^\epsilon(Y_{(j)} - \Vert \beta \Vert_2  z) \varphi(z) dz  \\
& \ge  & \int_{z: \ z < 0}  F^\epsilon(Y_{(j)} - \Vert \beta \Vert_2  z) \varphi(z) dz     \\
& \ge  &  \int_{\delta}^\infty  F^\epsilon(Y_{(j)}   +  \Vert \beta \Vert_2  z) \varphi(z) dz.
\end{eqnarray*}
Fix a small $\eta \in (0,1)$. Since $F^\epsilon$ is a cumulative distribution function, there exists $M_\eta > 0$ such that $F^\epsilon(t) \ge 1-\eta$ if $t > M_\eta$. This implies that for $\eta = 1/12$ and $\delta$ such $\int_\delta^\infty \varphi(z) dz \ge 3/8$ 
\begin{eqnarray*}\label{bignorm}
\Vert \beta \Vert_2 > \max((M_\eta - Y_{(j)})/\delta, 0) \Longrightarrow \int_{\delta}^\infty  F^\epsilon(Y_{(j)} + \Vert \beta \Vert_2 z) \varphi(z) dz  \ge 3(1-\eta)/8 = \frac{11}{32}.
\end{eqnarray*}
Let $M \ge \max\{M_\eta, [F^Y]^{(-1)}(1/4) + 2\}$. Since $Y_{(n/4)}$ converges almost surely to $ [F^Y]^{(-1)}(1/4)$, we have that
{\footnotesize{
\begin{eqnarray*}
\lim_{n \to \infty} \PP\Big(\frac{1}{n} \sum_{i=1}^n F^\epsilon(Y_{(j)} - \beta^T X_i)  \ge \frac{11}{32} \ \text{and} \ j=\lfloor n/5 \rfloor,\cdots, \lfloor n/4 \rfloor \ \text{and} \ \ \Vert \beta \Vert_2  \ge \frac{M - [F^Y]^{(-1)}(1/4)+1}{\delta} \Big)=1.
\end{eqnarray*}
}}
Now, note that for $\beta \in \RR^d$ such that $\int F^\epsilon(Y_{(j)}- \beta^T x) d F^X(x) \ge 11/32$ we have 
\begin{eqnarray*}
\DD_{n,2}(\beta) &\ge & \frac{1}{n} \sum_{j=\lfloor n/5 \rfloor}^{\lfloor n/4 \rfloor}   \left(F_n(Y_{(j)}) - n^{-1} \sum_{i=1}^n F^\epsilon(Y_j - \beta^T X_i)  \right)^2  \\
& \ge &  \frac{1}{n}\left(\lfloor n/4 \rfloor - \lfloor n/5 \rfloor + 1 \right)  (11/32-1/4)^2\\
&\ge & \frac{1}{20}\frac{9}{32^2}:= c.
\end{eqnarray*}  
Therefore
\begin{eqnarray*}
\lim_{n \to \infty} \PP\left( \inf_{\Vert \beta \Vert_2  \ge \frac{M - [F^Y]^{(-1)}(1/4)+1}{\delta} } \DD_{n,2}(\beta)  \ge c  \right) = 1.
\end{eqnarray*}
Now,
\begin{eqnarray*}
&& \PP\left(\Vert\hat{\beta}_n \Vert_2 \ge \frac{M - [F^Y]^{(-1)}(1/4)+1}{\delta}\right) = \\ && \PP\left(\Vert\hat{\beta}_n \Vert_2 \ge \frac{M - [F^Y]^{(-1)}(1/4)+1}{\delta}, \DD_{n,2}(\hat{\beta}_n) \ge c \right) + \\
&& \PP\left(\Vert\hat{\beta}_n \Vert_2 \ge \frac{M - [F^Y]^{(-1)}(1/4)+1}{\delta}, \DD_{n,2}(\hat{\beta}_n) < c \right)\\
&& = \PP_{n,1} + \PP_{n,2}.  
\end{eqnarray*}
Note that $\lim_{n \to \infty}\PP_{n,2} = 0$ since the event $\{\Vert \hat{\beta}_n \Vert_2  \ge  \frac{M - [F^Y]^{(-1)}(1/4)+1}{\delta}, \DD_{n,2}(\hat{\beta}_n) < c \}$ is included in the complement of the event 
\[
\left \{\inf_{\Vert \beta \Vert_2  \ge \frac{M - [F^Y]^{(-1)}(1/4)+1}{\delta} } \DD_{n,2}(\beta) \ge c \right \}.
\]
On the other hand, \eqref{RateDnbetahat} implies that $\DD_{n,2}(\hat{\beta}_n) =  O_{\PP}(n^{-1})$. Since $\PP_{n,1} \le \PP(\DD_{n,2}(\hat{\beta}_n) \ge c)$
\begin{eqnarray*}
\lim_{n\to\infty}\PP_{n,1} = 0.
\end{eqnarray*}
Thus
\begin{eqnarray*}
\lim_{n \to \infty} \PP\left( \Vert \hat{\beta}_n \Vert_2\ge\frac{M - [F^Y]^{(-1)}(1/4)+1}{\delta} \right) = 0.
\end{eqnarray*}
this implies that $\hat{\beta}_n = O_{\PP}(1)$. Now, consider the general case where $X \sim \mathcal{N}(\mu, \Sigma)$ with $\Sigma$ positive definite. Then, the model in \eqref{TheModel} can be written as 
\begin{eqnarray*}
Y \stackrel{d}{=} (\Sigma^{1/2}\beta_0)^T \Sigma^{-1/2}  X  + \epsilon  =  \gamma^T_0  \tilde{X} + \epsilon,
\end{eqnarray*}
where $\tilde{X} \sim \mathcal{N}(0, I_d)$, independent of $\epsilon$, and $\gamma_0 = \Sigma^{1/2}\beta_0$. Thus, as $\hat{\beta}_n =  \Sigma^{1/2} \hat{\gamma}_n$ and $\Vert \hat \gamma_n \Vert_2  = O_{\PP}(1)$, with $\hat{\gamma}_n$ the least squares estimator based on $Y_1, \ldots, Y_n$ and  $\tilde{X}_1,..., \tilde{X}_n$, it follows that $\Vert \hat \beta_n \Vert_2  = O_{\PP}(1)$. By Lemma~\ref{Charbeta} and Proposition~\ref{convNorm}, we have
\begin{eqnarray*}
\sqrt{\hat{\beta}^T_n \Sigma \hat{\beta}_n} \stackrel{\PP}{\rightarrow} c
\end{eqnarray*}
where $ c = \sqrt{\beta^T \Sigma \beta} > 0$ for all $\beta \in \mathcal{B}_0$ for this model. For $\widehat{\Sigma}_n$ a consistent estimator of $\Sigma$, we have
\begin{eqnarray*}
\hat{\beta}^T_n \widehat{\Sigma}_n \hat{\beta}^T_n  =   \hat{\beta}^T_n \Sigma \hat{\beta}^T_n + \hat{\beta}^T_n \left(\widehat{\Sigma}_n  - \Sigma \right) \hat{\beta}^T_n 
\end{eqnarray*}
where the second term on the right is bounded by $\max_{1 \le i, j \le d} \vert  \widehat{\Sigma}_{n, i,j} - \Sigma_{i, j} \vert \times \Vert \hat{\beta}_n \Vert^2  \stackrel{\PP}{\rightarrow} 0$ since $\Vert \hat{\beta}_n \Vert  = O_{\PP}(1)$. This concludes the proof. \hfill 
\end{proof}

\subsection{Proof of Theorem~\ref{conv2}}  
\begin{proof}
Recall that
\begin{eqnarray*}
\mathcal{D}_{2, F^Y}(\beta) =  \int\left(F^Y(y) d y - \int F^\epsilon(y-\beta^T x)d F^X(x) \right)^2 d F^Y(y).
\end{eqnarray*}
For $\beta \in \RR^d$ define
\begin{eqnarray*}
C_{\beta}(y) :=  \int F^\epsilon(y - \beta^T x)d F^X(x), \ \  \text{and} \ \ C_{n, \beta}(y) = \int  F^\epsilon(y - \beta^T x)  d F_n^X(x)
\end{eqnarray*}
the convolution density and its empirical estimator. Then, we have that
\begin{eqnarray*}\label{An&Bn}
\vert \DD_{n,2}(\beta) - \mathcal{D}_{2, F^Y}(\beta) \vert &  =  & \left \vert \int \left(F^Y_n(y) - C_{n,\beta}(y)\right)^2 d F^Y_n(y) - \int\left(F^Y(y) - C_{\beta}(y) \right)^2 dF^Y(y)   \right \vert  \notag \\
& \le & \left \vert \int \left(F^Y_n(y) - C_{n,\beta}(y) \right)^2 d(F^Y_n(y) - F^Y(y))  \right \vert \notag \\
&&  \  + \  \left \vert\int\left\{ \left(F^Y_n (y) -  C_{n, \beta}(y) \right)^2  -   \left(F^Y(y) -  C_{\beta}(y) \right)^2  \right \} d F^Y(y)   \right \vert \notag \\
& = & A_n +  B_n.
\end{eqnarray*}
Let us consider $A_n$.
\begin{eqnarray*}
0 \le A_n & \le & \left\vert \int F^Y_n(y)^2 d(F^Y_n(y) - F^Y(y))\right \vert + \left \vert \int C^2_{n,\beta}(y)  d(F^Y_n(y) - F^Y(y))\right \vert \\
&&  + \  2 \left\vert \int F^Y_n(y) C_{n,\beta}(y) d(F^Y_n(y) - F^Y(y))\right \vert.
\end{eqnarray*}
Functions $y \mapsto F^Y_n(y)^2$, $y \mapsto C^2_{n,\beta}(y)$, and $y \mapsto F^Y_n(y) C_{n,\beta}(y)$ are all non-negative, monotone and bounded above by $1$. Denote the class of real monotone functions $f$ such that $\text{Im}(f) \subseteq [0,1]$ by $\mathcal M$. We have 
\begin{eqnarray*}
0 \le A_n & \le &  4\sup_{f \in \mathcal{M}}\left\vert\int f d(F^Y_n(y) - F^Y(y)) \right \vert \\
& = &  4 n^{-1/2} \Vert \mathbb{G}_n \Vert_{\mathcal{M}}
\end{eqnarray*}

where $\mathbb G^Y_n = \sqrt{n}(\PP_n^Y - \PP^Y)$. By Theorem 2.7.5 of~\cite{aadbook} there exists a universal constant $K > 0$ such that for all $\eta > 0$ 
\begin{eqnarray*}
\log N_B(\eta,  \mathcal{M}, L_2(\PP^Y))  \le  \frac{K}{\eta}   
\end{eqnarray*}
where $N_B$ denotes the bracketing covering number. Now, by Lemma 3.4.2 of~\cite{aadbook} and using the fact that all functions in $\mathcal{M}$ are bounded above by $1$, we have that 
\begin{eqnarray*}
E[\Vert \mathbb{G}^Y_n \Vert_{\mathcal{M}}]  \lesssim J(1) \left(1 + \frac{J(1)}{\sqrt n}  \right),
\end{eqnarray*}
where for small $\delta > 0$
\begin{eqnarray*}
J(\delta) & = & \int_0^\delta\sqrt{1 + \log N_B(\eta, \mathcal{M}, L_2(\PP^Y))}d t    \\
& \le & \delta + \sqrt{K}\int_0^\delta \frac{1}{\sqrt{\eta}} d t \\
& = & \delta + 2 \sqrt{K}\sqrt{\delta}. 
\end{eqnarray*}

It follows that $\EE[\Vert \mathbb{G}^Y_n \Vert_{\mathcal{M}}] \lesssim 1$ for all $n \ge 1$. Using Markov's inequality, we conclude that 
\begin{eqnarray}\label{An}
A_n = O_{\PP}(n^{-1/2}).
\end{eqnarray}
Now, we focus on the term $B_n$. We have
\begin{eqnarray*}
\left(F^Y_n (y) -  C_{n, \beta}(y) \right)^2  & = & \left(F^Y_n (y) - F^Y(y) + F^Y(y)  -  C_\beta(y) + C_\beta(y) - C_{n, \beta}(y) \right)^2\\
& = &  \left(F^Y_n (y) - F^Y(y) \right)^2 + \left(F^Y(y) - C_\beta(y)\right)^2 +  \left(C_\beta(y) - C_{n, \beta}(y)\right)^2 \\
&& \ \ + 2 \left(F^Y_n (y) - F^Y(y)\right) \left(F^Y(y) - C_\beta(y)  \right) \\
&& \ \ + 2 \left(F^Y_n (y) - F^Y(y)\right) \left(C_\beta(y) - C_{n, \beta}(y)\right)   \\
&& \ \ + 2 \left(F^Y(y) - C_\beta(y)\right)\left(C_\beta(y) - C_{n, \beta}(y)\right).
\end{eqnarray*}
Hence
\begin{eqnarray*}
0 \le B_n  & \le & \int\left(F^Y_n (y) - F^Y(y) \right)^2 d F^Y(y) + \int\left(C_\beta(y) - C_{n, \beta}(y)\right)^2 d F^Y(y)\\
&&  +  \ \ 2 \int \left \vert F^Y_n (y) - F^Y(y) \right\vert\times\left\vert F^Y(y) -  C_\beta(y)\right\vert d F^Y(y)\\
&&  +  \ \ 2 \int\left\vert F^Y_n (y) - F^Y(y) \right\vert\times\left\vert C_\beta(y)  -  C_{n, \beta}(y)\right \vert d F^Y(y) \\
&&  +  \  \ 2 \int \left\vert F^Y(y) - C_\beta(y) \right\vert\times\left\vert C_\beta(y)  -  C_{n, \beta}(y)  \right \vert dF^Y(y)  \\
&& =  I_n + II_n(\beta) + III_n(\beta) +  IV_n(\beta) +  V_n(\beta),
\end{eqnarray*}
Note that
\begin{eqnarray*}
\EE(I_n) = n^{-1}\int F^Y(y) (1- F^Y(y)) d F^Y(y) \le n^{-1}
\end{eqnarray*}
implying that  $I_n = O_{\PP}(n^{-1})$.  Also, for all $\beta \in \mathbb  R^d$
\begin{eqnarray*}
II_n(\beta) & \le & \int \sup_{\gamma \in \RR^d} \left \vert C_{n, \gamma}(y) -  C_\gamma(y))\right\vert d F^Y(y),
\end{eqnarray*}
where
\begin{eqnarray*}
\left\vert C_{n, \gamma}(y) - C_\gamma(y)\right\vert  & = &  \vert\int F^\epsilon(y - \gamma^T x)d(F^X_n(x)  -  F^X(x))\vert\\
& = & \vert\int F^\epsilon(y -  z) d(F_n^{Z}(z) - F^Z(z))\vert  \\
& = & \vert\int (F_n^Z(z) - F^Z(z)) f^\epsilon (y -  z) dz\vert.
\end{eqnarray*}
using the change of variable $Z = \gamma^T X$, and integration by parts. Thus, it follows that 
\begin{eqnarray}\label{emp2}
\sup_{y \in \RR}\vert C_{n, \gamma}(y) - C_\gamma(y)\vert & \le &  \sup_{z \in \RR}\vert F^Z_n(z) - F^Z(z) \vert \notag  \\
& = & \sup_{z} \left\vert\int\one\{\gamma^T x \le z\} d(F^X_n(x) - F^X(x)) \right \vert.
\end{eqnarray}
Next, we show that the supremum in \eqref{emp2} is $O_{\PP}(n^{-1/2})$ independently of $\gamma$. Consider the class
\begin{eqnarray*}
\mathcal{F} = \{x\mapsto\gamma^T x - z, \gamma\in\RR^d, z \in\RR\}.  
\end{eqnarray*} 
Note that $\mathcal{F}$ is a $d+1$ dimensional vector space. By Lemma 2.6.15 of~\cite{aadbook}, $\mathcal{F}$ is a VC-subgraph of index smaller than $d+1+2 = d+3$.  Now, note that 
\begin{eqnarray*}
\int\one\{\gamma^T x \le z\} d(F^X_n(x) - F^X(x)) & = & \int\one\{\gamma^T x - z \le 0\}  d(F^X_n(x) - F^X(x)) \\
& = & \int(\phi\circ f)(x)d(F^X_n(x) - F^X(x))
\end{eqnarray*}
with $\phi(t) = \one\{t \le 0\}$ and $f(x) = \gamma^T x - z\in\mathcal{F}$. Since $\phi$ is monotone, by Lemma 2.6.18(viii) of~\cite{aadbook}, we conclude that the class $\phi\circ \mathcal{F}$ is a VC-subgraph.  Since all elements in $\phi \circ \mathcal F$ are bounded by $1$, the latter can be taken as its envelope. Thus, it follows from Theorem 2.6.7 in ~\cite{aadbook} that there exists $V \ge 2$ such that (taking $r=2$) for all $\eta \in (0,1]$ and all probability measures $\mathbb{Q}$
\begin{eqnarray*}
N(\eta, \phi\circ\mathcal{F}, L_2(\mathbb{Q}))\le K V (16 2)^V \frac{1}{\eta^{2(V-1)}}
\end{eqnarray*}
where $K > 0$ is some universal constant. Hence, there exists a constant $c \in \RR $ such that 
\begin{eqnarray*}
\sup_{\mathbb{Q}} \int_0^1 \sqrt{1 + \log N(\eta, \phi\circ\mathcal{F}, L_2(\mathbb{Q}))} d\eta & \le & c + \sqrt{2(V-1)} \int_0^1  \sqrt{\log\left(\frac{1}{\eta}\right)} d\eta \\
& = &  c + \sqrt{2(V-1)}\int_{1}^\infty\frac{\sqrt{\log t}}{t^2} d t < \infty,
\end{eqnarray*}
using the fact that
\begin{eqnarray}\label{Int}
\int_1^\infty \frac{\sqrt{\log(t)}}{t^2} d t \le \int_1^\infty \frac{1}{t^{3/2}}d x  < \infty
\end{eqnarray}
using $\log(t)\le t, t \ge 1$.    It follows from Theorem 2.14.1 in~\cite{aadbook} that 
\begin{eqnarray*}
\EE[\Vert \mathbb{G}^X_n \Vert_{\phi \circ \mathcal{F}}] \lesssim 1,
\end{eqnarray*}
where 
\[
\Vert \mathbb G^X_n \Vert_{\phi \circ \mathcal{F}} = \sup_{g\in\phi\circ\mathcal{F}} n^{1/2} \left\vert\int g(x) d(F^X_n(x) - F^X(x))\right\vert. 
\]
Using Markov's inequality, it follows that 
\begin{eqnarray*}
\sup_{g \in\phi\circ\mathcal{F}}\left\vert\int g d(F^X_n - F^X) \right \vert =  O_{\PP}(n^{-1/2}),
\end{eqnarray*}
from which we conclude that 
\begin{eqnarray}\label{boundcngamma}
\sup_{y \in\RR, \gamma\in\RR^d} \vert C_{n, \gamma}(y) - C_{\gamma}(y) \vert = O_{\PP}(n^{-1/2}).
\end{eqnarray}
Since $ II_n(\beta)\le\sup_{y \in \RR, \gamma\in\RR^d} \vert C_{n, \gamma}(y) - C_{\gamma}(y) \vert$ for all $\beta\in\RR^d$ it follows that
\begin{eqnarray*}
\sup_{\beta \in \RR^d}  II_n(\beta)  =  O_{\PP}(n^{-1/2}).
\end{eqnarray*}
By the Cauchy-Schwarz inequality, we have that 
\begin{eqnarray*}
III_n(\beta) & \le   &  2 \left( \int \left \vert \mathbb F^Y_n (y) - F^Y(y) \right \vert^2  dF^Y(y)  \right)^{1/2}  \left( \int \left \vert  F^Y(y)  -  C_\beta(y) \right \vert^2  dF^Y(y) \right)^{1/2}   \\
& \le &  2 \left( \int \left \vert \mathbb F^Y_n (y) - F^Y(y) \right \vert^2  dF^Y(y)  \right)^{1/2},
\end{eqnarray*}
since $\vert F^Y(y) - C_\beta(y)\vert \le 1$. Therefore
\begin{eqnarray*}
\sup_{\beta \in \RR^d} III_n(\beta) =  O_{\PP}(n^{-1/2}).
\end{eqnarray*}
Similarly, we have  
\begin{eqnarray*}
\sup_{\beta \in \RR^d}  IV_n(\beta)  =  O_{\PP}(n^{-1/2}).
\end{eqnarray*}
Also,
\begin{eqnarray*}
\sup_{\beta \in \RR^d}  V_n(\beta) & \le   &  2   \sup_{y \in \RR, \beta \in \RR^d}  \vert C_{n, \beta}(y)  - C_{\beta}(y)  \vert\times\int\left\vert F^Y(y) - C_\beta(y) \right\vert dF^Y(y) \\
& \le &  2   \sup_{y \in \RR, \beta \in \RR^d} \vert C_{n, \beta}(y) - C_{\beta}(y) \vert   =  O_{\PP}(n^{-1/2}).
\end{eqnarray*}
It follows that 
\begin{eqnarray*}
\sup_{\beta\in\RR^d}\vert \DD_{n,2}(\beta) - \mathcal{D}_{2, F^Y}(\beta)\vert \stackrel{\PP}{\rightarrow} 0.
\end{eqnarray*}
In the following, we show that 
\begin{eqnarray}\label{sep}
0 = \mathcal{D}_{2, F^Y}(\beta_0)  <  \inf_{\beta\not\in\mathcal{O}} \mathcal{D}_{2, F^Y}(\beta)
\end{eqnarray}
for open set $\mathcal{O}\subset\RR^d $ such that $\beta_0 \in \mathcal{O}$. Note that $F^Y  = C_{\beta_0}$, therefore
\begin{eqnarray*}
\mathcal{D}_{2, F^Y}(\beta) = \int\left(C_\beta(y) - C_{\beta_0}(y)\right)^2 d F^Y(y).
\end{eqnarray*}
Let $\mathcal{O} = B(\beta_0, r)\subset\RR^d$ be the open Euclidean Ball of center $\beta_0$ and radius $r > 0$. Suppose that $\inf_{\beta\not\in\mathcal{O}}\mathcal{D}_{2, F^Y} = 0$. This means that we can find a sequence $\{\beta_m\}_{m \ge 1}\notin B(\beta_0, r)$ such that 
\begin{eqnarray}\label{Limit}
\lim_{m \to \infty}\int\left(C_{\beta_m}(y) - C_{\beta_0}(y)\right)^2 d F^Y(y) = 0.
\end{eqnarray}
Assume there exists a subsequence $\{\beta_{m^\prime}\}$ such that $\lim_{m^\prime\to \infty}  \Vert \beta_{m^\prime} \Vert_2  = \infty$. Let $u_{m^\prime} = \beta_{m^\prime}/\Vert \beta_{m^\prime} \Vert_2$. Since the unit sphere $S^{d-1}$ is compact, there exists a subsequence of $\{u_{m^\prime}\}$, which without loss of generality we denote by $(u_{m^\prime})_{m^\prime}$, such that $\lim_{m^\prime\to\infty} u_{m^\prime} = u$ where $u\in S^{d-1}$. Note that for a fixed $y$
\begin{eqnarray*}
C_{\beta_{m^\prime}}(y) & =  & \int F^\epsilon(y - \Vert \beta_{m^\prime} \Vert_2 u^T_{m^\prime}x) d F^X(x) \\
& = & \int F^\epsilon(y - \Vert \beta_{m^\prime} \Vert_2 u^T_{m^\prime}x) \one\{u^T x = 0\} d F^X(x) \\
&& +  \int F^\epsilon(y - \Vert \beta_{m^\prime} \Vert_2 u^T_{m^\prime} x) \one\{u^T x < 0\}  d F^X(x)  \\
&& + \  \int F^\epsilon(y - \Vert \beta_{m^\prime} \Vert_2 u^T_{m^\prime}x) \one\{u^T x > 0\} d F^X(x) \\
& = &  I_{m^\prime}(y) + II_{m^\prime}(y) + III_{m^\prime}(y).
\end{eqnarray*}
Since $X$ does not belong to any affine space with probability one, we have 
\begin{eqnarray*}
\int \one\{u^T x = 0\} d F^X(x) = 0,
\end{eqnarray*}
and therefore $I_{m^\prime}=0$.  Also, for all $x$ such that $u^Tx < 0$, we have  $\Vert \beta_{m^\prime} \Vert_2 u^T_{m^\prime}x   \to -\infty$ and hence $F^\epsilon(y - \Vert\beta_{m^\prime} \Vert_2 u^T_{m^\prime}x) \to 1$ as $m^\prime\to \infty$. Similarly, for all $x$ such that $u^Tx > 0$, $F^\epsilon(y - \Vert \beta_{m^\prime} \Vert_2 u^T_{m^\prime}x) \to 0$.  By the Dominated Convergence Theorem, it follows that 
\begin{eqnarray*}
\lim_{m^\prime \to \infty}  II_{m^\prime}(y)  = 1, \ \ \text{and} \  \ \lim_{m^\prime \to \infty}  III_{m^\prime}(y)  = 0.
\end{eqnarray*}
Using the Dominated Convergence Theorem again, we conclude that 
\begin{eqnarray*}
\lim_{m^\prime\to\infty} \int \left ( C_{\beta_{m^\prime}}(y) - C_{\beta_0}(y)  \right)^2 d F^Y(y) =  \int\left(\PP(u^TX < 0) - C_{\beta_0}(y)\right)^2 d F^Y(y),
\end{eqnarray*}
which implies with \eqref{Limit} that $C_{\beta_0}(y) = F^Y(y) = \PP(u^TX < 0)$ for all $y \in \RR$, which is impossible.  We conclude that it is impossible that the sequence $(\beta_{m^\prime})_{m^\prime}$ has an Euclidean norm that tends to $\infty$. This means that it is bounded by some constant $K > 0$. Then, there exists a subsequence, w.o.l.g. we denote it by $\{\beta_{m^\prime}\}$, converging to some vector $\tilde{\beta}$. Using similar arguments as above, we have
\begin{eqnarray*}
\int\left(C_{\tilde{\beta}}(y) - C_{\beta_0}(y)\right)^2 d F^Y(y) = 0,
\end{eqnarray*}
and hence $C_{\tilde{\beta}} = C_{\beta_0}$. Since $\beta_0$ is the unique vector such that $Y \stackrel{d}{=} \beta^T_0 X + \epsilon$, we must have $\tilde{\beta} = \beta_0$. However, this is impossible because this would mean that $r = 0$. Since any open set $\mathcal{O}$ containing $\beta_0$ contains $\mathcal{B}(\beta_0, r)$ for some $r > 0$, we conclude that the separability condition in  \eqref{sep} must hold. As $\hat{\beta}_n$ minimizes $\mathbb D_{n, 2}$, it follows that all the conditions of (i) in~\cite[Corollary 3.2.3]{aadbook} are fulfilled. Therefore
\begin{eqnarray*}
\hat{\beta}_n\stackrel{\PP}{\rightarrow}\beta_0.
\end{eqnarray*} 
This concludes the proof. \hfill
\end{proof}
\subsection{Proof of Theorem~\ref{pole}}
\begin{proof}
Note that $\beta\in\mathcal{B}_0$ if and only if 
\begin{eqnarray*}
\EE\left(e^{\beta_0^T X t}\right) = \EE\left(e^{\beta^T X t}\right)
\end{eqnarray*}
for all $t\in\RR$ such that the expectations on the left and right are defined. To use a simpler notation, we write now $G$ for $G_{X^1}$. By the i.i.d. assumption, we can write that 
\begin{eqnarray*}\label{mg1sided}
G(\beta_{0, 1} t)\times\cdots\times G(\beta_{0, d} t) = G(\beta_1 t)\times\cdots\times G(\beta_d t)
\end{eqnarray*}
for all $t$ such that 
\begin{eqnarray}\label{Region}
\max_{1\le i \le d} (\beta_{0, i} t) < \lambda, \  \  \text{and}  \ \ \max_{1\le i \le d} (\beta_i t) < \lambda.
\end{eqnarray}
By ordering the components of $\beta_0$ and $\beta$ we will assume without generality that 
$\beta_{0,1}\le\ldots\le\beta_{0,d}$ and $\beta_1\le \ldots\le\beta_d$. Suppose that $\beta_{0, d}\ne\beta_d$. By symmetry, we can assume without loss of generality that $\beta_d > \beta_{0, d}$. There are 3 cases to consider. 
\begin{itemize}
\item $\beta_d > 0$. Then, let $t$ such that $t \nearrow \lambda/\beta_d$. Then, $t$ belongs to the permissible region in \eqref{Region} and 
\begin{eqnarray*}
\frac{1}{(\lambda - \beta_{0,1} t)^\alpha}\times\ldots   \times \frac{1}{(\lambda  - \beta_{0, d} t)^\alpha} =  \frac{1}{(\lambda - \beta_1 t)^\alpha}\times\ldots\times  \frac{1}{(\lambda - \beta_d t)^\alpha}K(t)
\end{eqnarray*}
for some function $K$ such that $K(t) \in (0, \infty)$ for such $t$. This implies that 
\begin{eqnarray*}
\lim_{t \nearrow\lambda/\beta_d}\frac{1}{(\lambda - \beta_{0,1} t)^\alpha}\times\ldots\times\frac{1}{(\lambda - \beta_{0,d} t)^\alpha} = \infty
\end{eqnarray*}
which is impossible. Hence, $\beta_{0,d} = \beta_d$.
  
\item $\beta_{0,d} < 0$.  Then, this implies that $\beta_{0,1}\le\ldots\le\beta_{0,d} < 0$ and $\beta_1\le\ldots\le \beta_d < 0$. We show that we have $\beta_{0,1} = \beta_1$.  Indeed, suppose that $\beta_1 > \beta_{0,1}$. Take now $t$ such that $t \beta_1  \nearrow \lambda$.  This means that $ t \searrow \lambda /\beta_{0,1} < 0$. Then, $t \beta_{0,i} \le  t \beta_{0,1}$ for $i \ge 2$ implying that $\max_{1 \le i \le d} (\beta_{0,i} t)  = \beta_{0,1} t$. Also, $\max_{1 \le i \le d} (\beta_{0,i} t)  = \beta_1 t  < \beta_{0,1} t$. Thus,  $t$ belongs to the permissible region in (\ref{Region}).  Then, we have that
\begin{eqnarray*}
\lim_{ t \searrow \lambda/\beta_{0,1}}\frac{1}{(\lambda - \beta_{0,1} t)^\alpha}\times\ldots\times\frac{1}{(\lambda - \beta_{0,d} t)^\alpha} = \infty  
\end{eqnarray*}
which implies that
\begin{eqnarray*}
\lim_{ t \searrow \lambda/\beta_{0,1}}\frac{1}{(\lambda - \beta_1 t)^\alpha}\times\ldots\times\frac{1}{(\lambda - \beta_d t)^\alpha} = \infty  
\end{eqnarray*}
which is impossible. Using a similar argument, we can show that we cannot have $\beta_{0,1} < \beta_1$. Hence, $\beta_{0,1} = \beta_1$. Also, one can show successively that $\beta_{0,i} = \beta_i $ for $i = 2, \cdots, d-1$ and finally that $\beta_{0,d} = \beta_d$.

\item  $\beta_d = 0$. Suppose that for all $i = 1, \cdots, d$ $\beta_i =0$.   Then,
\begin{eqnarray*}
\frac{1}{(\lambda - \beta_{0,1} t)^\alpha}\times\ldots\times \frac{1}{(\lambda - \beta_{0,d} t)^\alpha} H(\beta_{0,1} t) \times\ldots\times H(\beta_d t ) = \frac{1}{\lambda^{d\alpha}} H(0)^d
\end{eqnarray*}
for all $t$ such that $\max_{1\le i \le d}  (\beta_{0,i} t) < \lambda$.  It is easy to see that this implies that $\beta_{0,1} = \cdots = \beta_{0,d}$, which can be shown by assuming that this does not hold.   Now, suppose one of the coefficients $\beta_{0,i}, i < d \ne 0$. Consider then the integer $r$ such that $\beta_{r+1} = \cdots = \beta_d = 0$ and $\beta_r \ne 0$. Note that we necessarily have $\beta_1 \le\cdots\le\beta_r < 0$. Recall also that we are in the case where $\beta_{0,1}\le\cdots\le\beta_{0,d} < 0$. Then, using the same argument as in the second case, we can successively show that 
\begin{eqnarray*}
\beta_{0,1} = \beta_1, \cdots, \beta_{0,r} = \beta_r.
\end{eqnarray*}
These equalities will enable us to show that $\beta_{0, r+1}= \cdots = \beta_{0,d} = 0$. 
\end{itemize} 
Thus, in all the cases considered, the assumption $\beta_{0,d} \ne\beta_d$ leads to a contradiction. Using a recursive reasoning, we conclude that we should have $\beta_{0,i} = \beta_i$ for all $i = 1,\cdots, d$. Therefore, in general, and without assuming that the elements of $\beta_0$ and $\beta$ are ordered, $\mathcal{B}_0$ consists of all vectors $\beta$ such that for some permutation $\pi$ we have $\beta_i = \beta_{0, \pi(i)}$. Consequently $\Vert\beta\Vert_2 = \Vert\beta_0\Vert_2$ and hence $\mathcal{B}_0\subset\{\beta:\Vert\beta\Vert_2 = \Vert\beta_0\Vert_2\}$.
\end{proof}
\subsection{Auxiliary results for the proof of Theorem \ref{limitDist}}
In the sequel, we would need the following definition. For a given class of functions $\mathcal{F}$ with envelope $F$ we define
\begin{eqnarray}\label{J}
J(\delta, \mathcal{F}) = \sup_{\mathbb{Q}}\int_0^\delta\sqrt{1 + \log N(\eta \Vert F \Vert_{\mathbb{Q}, 2}, \mathcal F, L_2(\mathbb{Q})}  d\eta
\end{eqnarray}
where   $N(\nu, \mathcal F, \Vert \cdot \Vert)$ is the $\nu$-covering number of $\mathcal F$ with respect to $\Vert \cdot \Vert$. The supremum in \eqref{J} is taken over probability measure $\mathbb{Q}$ such that $\Vert F\Vert^2_{\mathbb{Q}, 2}=  \int F^2 d\mathbb{Q}  > 0$. 

To make the notation used below more compact, we shall use the classical empirical process notation 
\begin{eqnarray*}
\mathbb{G}^X_n := \sqrt{n}(\PP^X_n - \PP^X), \  \ \text{and} \  \ \mathbb{G}^Y_n := \sqrt{n} (\PP^Y_n - \PP^Y).
\end{eqnarray*}
For random vector $Z$ we define its norm in $\mathcal{F}$ as
\begin{eqnarray*}
\Vert\mathbb{G}^Z _n \Vert_{\mathcal{F}} := \sup_{f \in \mathcal{F}} \vert \mathbb{G}^Z f \vert.
\end{eqnarray*}
Additionally, we use the fact that 
\begin{eqnarray*}
\int_0^1 \sqrt{\log\left( \frac{1}{t} \right)}  dt  & =  & \int_1^\infty \frac{\sqrt{\log(x)}}{x^2} d x
\end{eqnarray*}
together with the inequality in (\ref{Int}) shown above

\begin{proposition}\label{ClassG}
Consider the classes of functions $\mathcal{I}$, $\mathcal{M}_{\beta_0}$, and $\mathcal{G}$ defined as
\begin{eqnarray*}
\mathcal{I} = \Big\{I_y(t) = \one\{t \le y\}, y \in\RR\Big \},
\end{eqnarray*}
\begin{eqnarray*}
\mathcal{M}_{\beta_0} := \Big\{f_y(x) = m_{\beta_0, y}, y \in \RR \Big \},
\end{eqnarray*}
with $m_{\gamma, y}(x) = F^\epsilon(y - \gamma^Tx), \ \gamma \in \RR^d$, and 
\begin{eqnarray*}
\mathcal{G} := \Big\{g_{\beta, y}(x) =  m_{\beta, y}(x) - m_{\beta_0, y}(x),\ \Vert\beta - \beta_0\Vert_2\le\delta, y\in\RR\Big \}.
\end{eqnarray*}
Then, we have
\begin{eqnarray*}
\EE\left[\Vert \mathbb{G}^Y_n\Vert^2_{\mathcal{I}}\right] \lesssim 1, \  \ \text{and} \ \  \EE\left[\Vert \mathbb{G}^Y_n\Vert^2_{\mathcal{M}_{\beta_0}} \right] \lesssim 1, \  \ \text{and} \ \ \EE\left [\Vert \mathbb{G}^X_n \Vert^2_{\mathcal{G}}\right] \lesssim \delta^2.
\end{eqnarray*}
\end{proposition}

\begin{proof}\textbf{(Proposition \ref{ClassG}).} We begin with the class $\mathcal{I}$. This class is a VC-class with index equal to $2$; see Example 2.6.1 in~\cite{aadbook}. Since $F \equiv 1$ is an envelope for the class, it follows from  Theorem 2.6.7 of~\cite{aadbook} that for any probability measure $\mathbb{Q}$ and $\eta \in (0,1)$
\begin{eqnarray*}
N(\eta, \mathcal{I}, L_2(\mathbb{Q})) \le 2 K (16e)^2\frac{1}{\eta^2}
\end{eqnarray*}
for some universal constant $K > 0$.  With the definition given in \eqref{J} and inequality $\sqrt{a+b} \le \sqrt a + \sqrt b $ for any $a \ge 0, b \ge 0$, it follows that
\begin{eqnarray*}
J(1, \mathcal I) & \le & 1 + \left(\sqrt{\log(2 K (16e)^2 ) \vee 0} + \sqrt{2} \int_0^1 \sqrt{\log \left(\frac{1}{\eta}\right)}\right) d\eta  < \infty
\end{eqnarray*}
by \eqref{Int}. By Theorem 2.4.1 in~\cite{aadbook}, we have that 
\begin{eqnarray*}
\EE\left[\Vert\mathbb{G}^Y_n \Vert^2_{\mathcal{I}}\right]\lesssim 1.
\end{eqnarray*}
For the class $\mathcal{M}_{\beta_0}$, we use a small adaptation of Lemma 2.6.16 in~\cite{aadbook} with $\psi = F^\epsilon$ to claim that it is a VC-class with index $2$.   The main difference between our setting and the one in that lemma is that we have $\beta^T_0 x$ instead of a univariate $x$. However, the main argument in the proof of this lemma remains the same.  Since $1$ is an envelope, we conclude using arguments similar to those above to show that 
\begin{eqnarray*}
\EE\left[\Vert\mathbb{G}^Y_n \Vert^2_{\mathcal{M}_{\beta_0}} \right]\lesssim 1.
\end{eqnarray*}
Now, we turn to the class $\mathcal{G}$, and we start by showing that it is contained in a VC-hull. Without loss of generality, suppose that $\beta^T_0 x\le\beta^T x$. Then, 
\begin{equation}\label{Mix}
\begin{split}
g(x) &= m_{\beta, y}(x) - m_{\beta_0, y}(x) \\
& =  \int\one\{t \in [y - \beta^T x, y - \beta^T_0 x]\}d F^\epsilon(t)\\
& =  \int\one\{\beta^T x + t\ge y\}\cdot\one\{\beta_0^T x + t \le y\} d F^\epsilon(t) \\ 
& =  \int\phi_{1}(\beta^T x + t - y)\cdot\phi_{2}(\beta^T_0 x + t - y) d F^\epsilon(t)
\end{split}
\end{equation}
where
\begin{eqnarray*}
\phi_{1}(u) =  \one\{u \ge 0 \},  \  \ \textrm{and}  \  \  \phi_{2}(u)  =  \one\{u \le 0\}, \  \  \ u \in \RR.
\end{eqnarray*}
Now, consider the class of functions 
\begin{eqnarray*}
\mathcal{F} =  \{x\mapsto\gamma^T x + t - y : \gamma\in\RR^d, t, y\in\RR\}.
\end{eqnarray*}
Then, $\mathcal{F}$ is a finite-dimensional space with dimension $d+1$. Lemma 2.6.15 of \cite{aadbook} implies that this class is a VC-subgraph of index smaller than $d+3$.    Let $\mathcal{H}$ be the class defined as 
\begin{eqnarray*}
\mathcal{H}  =  \{x \mapsto\phi_{1}(\beta^T x + t - y)\cdot\phi_{2}(\beta^T_0 x + t - y), \ \ \beta\in\RR^d, t, y\in\RR\}.
\end{eqnarray*}  
Let us denote $\mathcal{H}_1 = \phi_1\circ\mathcal{F}$ and $\mathcal{H}_2 = \phi_2\circ\mathcal{F}$. The key observation is that 
\begin{eqnarray*}
\mathcal{H}\subset\mathcal{H}_1 \cdot\mathcal{H}_2  = \mathcal{H}_1 \wedge\mathcal{H}_2 
\end{eqnarray*}
where $\wedge$ denotes the minimum. The equality above is true since $\phi_1$ and $\phi_2$ are indicator functions which implies that for any pair $(f_1, f_2)\in\mathcal{F}^2$ we have that  $\phi_1(f_1)\phi_2(f_2) = \phi_1(f_1)\wedge\phi_2(f_1)$. Now, It follows by (i) of Lemma 2.6.18 in~\cite{aadbook} that $\mathcal{H}_1 \wedge\mathcal{H}_2 $ is VC. This implies that $\mathcal{H}$ is VC. Given the definition of $g \in\mathcal{G}$ in (\ref{Mix}), it follows that the class $\mathcal{G}$ is a subset of the convex hull of $\mathcal H$. By the mean value Theorem, we can find  some real number $\theta^*$  (depending on $x$, $y$, $\beta_0$ and $\beta$) such that 
\begin{eqnarray*}
g(x) = (\beta-\beta_0)^T x f^\epsilon(\theta^*).
\end{eqnarray*}
Thus, using the assumption that $f^\epsilon \le M$ and the Cauchy-Schwarz inequality, it follows that 
\begin{eqnarray*}
G(x) = M \delta \Vert x \Vert_2
\end{eqnarray*}
is an envelope of $\mathcal{G}$. By Theorem 2.6.9 of \cite{aadbook}, we can find an integer $V \ge 2$ and a universal constant $K > 0$ such that we have for any probability measure $Q$ satisfying
$\Vert G\Vert_{\mathbb{Q}, 2}\in (0, \infty)$, i.e., $\int\Vert x \Vert_2 d\mathbb{Q}(x) \in (0, \infty)$
\begin{eqnarray*}
\log N\left(\eta\Vert G\Vert_{\mathbb{Q}, 2}, \mathcal{G}, L_2(\mathbb{Q})\right)\le K \left(\frac{1}{\eta}\right)^{2V/(V+2)}.
\end{eqnarray*}
Hence, with $a = V/(V+2)$
\begin{eqnarray*}
J(1, \mathcal{G})\le 1 + \sqrt K \int_0^1 \frac{1}{\eta^{a}} d\eta  < \infty
\end{eqnarray*}
since $a < 1$. Under the assumption that $\int\Vert x \Vert^2 d F^X(x) < \infty$, we are allowed to use Theorem 2.4.1 in~\cite{aadbook} to conclude that 
\begin{eqnarray*}
\EE\left[\Vert\mathbb{G}^X_n\Vert^2_{\mathcal{G}} \right]\lesssim\delta^2
\end{eqnarray*}
which completes the proof.  \hfill 
\end{proof}

\begin{proposition}\label{sumsM}
Let $M$ be the class of monotone functions $f$ such that $f\in [0,1]$. Consider the class
\begin{eqnarray*}
\mathcal{C} =  \mathcal{M} + \mathcal{M} - 2\mathcal{M} = \left\{f: f = f_1 + f_2 - 2 f_3,  (f_1,f_2,f_3) \in \mathcal{M}^3\right\}.
\end{eqnarray*} 
Then, there exists a universal constant $A > 0$ such that for all $\eta \in (0,1)$
\begin{eqnarray*}
\log N_B(\eta, \mathcal{C}, L_2(\PP))\le\frac{A}{\eta}
\end{eqnarray*}
where $N_B(\eta, \mathcal{C}, L_2(\PP))$ denotes the $\eta$-bracketing number for $\mathcal{C}$ with respect to $L_2(\PP)$.
\end{proposition}

\begin{proof}\textbf{(Proposition~\ref{sumsM}).} Fix $\eta$ and let $N = N_B(\eta, \mathcal{M}, L_2(\PP))$. For $f\in\mathcal{C}$, we can find $1 \le i, j, k \le N$ and brackets $(L_i, U_i)$, $(L_j, U_j)$ and $(L_k, U_k)$ such that 
\begin{eqnarray*}
L_i \le f_1 \le U_i, \  L_j \le  f_2  \le U_j, \  L_k \le f_3 \le U_k
\end{eqnarray*}
and $\int (U_i-L_i)^2 d\PP \le\eta^2, \int (U_j-L_j)^2 d\PP\le\eta^2$ and $\int (U_k-L_k)^2 d\PP\le\eta^2$. Then,
\begin{eqnarray*}
L_i +  L_j - 2 U_k \le f \le U_i + U_j - 2L_k 
\end{eqnarray*}
and
\begin{eqnarray*}
\int \left(U_i + U_j - 2 L_k - L_i - L_j + 2 U_k \right)^2 d\PP  & = &  \int \left( U_i  -L_j  + U_j - L_j - 2 (U_k - L_k)\right)^2 d\PP \\
& \le & 4\int (U_i-L_i)^2 d\PP + 4\int (U_j-L_j)^2 d\PP \\
&& \ + 8 \int (U_k-L_k)^2 d\PP \\
& \le & 16 \eta^2,
\end{eqnarray*}
using the fact that $(a+b+c)^2 \le 2 (a+b)^2 + 2 c^2 \le 4(a^2 + b^2) + 2 c^2$. Hence,  $N_B(4 \eta, \mathcal{C}, L_2(\PP))\le N^3$  implying that
\begin{eqnarray*}
\log N_B(4\eta, \mathcal{C}, L_2(\PP))\le 3 \log N\le\frac{3 K}{\eta}, \  \text{for some universal constant $K > 0$}
\end{eqnarray*}
where the last inequality follows from Theorem 2.7.5 of~\cite{aadbook}. It follows that
\begin{eqnarray*}
\log N_B(\eta, \mathcal{C}, L_2(\PP))\le\frac{A}{\eta}
\end{eqnarray*}
with $A = 12 K$. \hfill
\end{proof}

\begin{proposition}\label{pwmon}
Let $f$ be a real function such that $0 \le f \le M$ for some constant $M > 0$. Suppose, there exist real numbers $a < b$ such that $f$ is non-decreasing $(-\infty, a]$, non-increasing on $(a, b)$ and non-increasing on $[b, \infty)$. Let $X\in\RR^d$ be a random vector such that $\EE\Vert X\Vert_2^2 < \infty$. Let $F_n^X$ be the empirical distribution of $X_1, \cdots, X_n$ that are independent and identically distributed as $X$. Then
\begin{eqnarray*}
\sup_{(\beta, z) \in \RR^{d} \times \RR}  \left \Vert \int  x f(\beta^T x  + z) d F_n^X(x) \right \Vert_2 = O_{\PP}(1).
\end{eqnarray*}
\end{proposition}

\begin{proof}\textbf{(Proposition \ref{pwmon}).}  Consider the function class $\mathcal{K} = \{k_{\beta, z}: k_{\beta, z}(x) = \beta^T x + z, \beta\in\RR^d, z\in\RR\}$ and note that it is $d+1$ dimensional. By Lemma 2.6.15 of~\cite{aadbook}, $\mathcal{K}$ is a VC class with index $V\le d+3$. Define the functions
\begin{eqnarray*}
f_{M, a}(t) &=& f(t)\one\{t \le a\} + M\one\{t > a\}, \\
f_{M, b}(t) &=& f(t)\one\{t \ge b\} + M\one\{t < b\},\\
f_{M, a, b}(t) &=& f(t)\one\{a < t < b\} + M\one\{t \le a\}.
\end{eqnarray*}
Fix $h$ such that $\Vert h \Vert_2 \le K$. Then, we can write
\begin{eqnarray*}
\int x f(\beta^T x + z) d F_n^X(x) &=& \int x f \circ k_{\beta, z}(x)   d F_n^X(x) \\                                                        & = &  \int x f \one\{k_{\beta, z}(x) \le a\} d F_n^X(x) + \int x f \one\{a < k_{\beta, z}(x) < b\}  d F_n^X(x)  \\
&& \  +  \int x f\one\{k_{\beta, z}(x) \ge b\} d F_n^X(x)\\
& = & \int x f_{M, a} \circ  k_{\beta, z}(x)  d F_n^X(x) - M \int x \one\{k_{\beta, z}(x) > a\}  d F_n^X(x) \\
&& + \int x  f_{M, a, b} \circ  k_{\beta, z}(x) d F_n^X(x) - M \int x \one\{k_{\beta, z}(x) \le a\} d F_n^X(x)  \\
&& + \int x  f_{M, b} \circ  k_{\beta, z}(x) d F_n^X(x) - M \int x \one\{k_{\beta, z}(x) < b\} d F_n^X(x) \\
& = & \int  x f_{M, a} \circ  k_{\beta, z}(x) d F_n^X(x) + \int x f_{M, a, b} \circ  k_{\beta, z}(x) d F_n^X(x)  \\
&&  + \int  x f_{M, b} \circ  k_{\beta, z}(x) d F_n^X(x) - M \int x \one\{k_{\beta, z}(x) < b\} d F_n^X(x)  \\
& = &  T_{n, 1}(\beta, z, h) + T_{n, 2}(\beta, z, h) + T_{n, 3}(\beta, z, h) + T_{n, 4}(\beta, z, h).  
\end{eqnarray*}
Functions  $f_{M, a}$, $f_{M, b}$, $f_{M, a, b}$, and  $\one\{\cdot < b\}$ are all bounded monotone functions. Following from Lemma 2.6.18 (viii) and (vi),  Theorem 2.6.7, Theorem 2.14.1  of~\cite{aadbook} for $j \in \{1, 2, 3, 4 \}$ we have
\begin{eqnarray*}
\EE\left[\left(\sup_{\beta, z}\vert T_{n, j}(\beta ,z) \vert\right)^2\right] &\lesssim  &  \EE[\Vert X\Vert^2_2]  < \infty.
\end{eqnarray*}
Application of Chebychev's inequality yields the claim.    \hfill 
\end{proof}

\begin{proposition}\label{empcontrol}
There exists a constant $C > 0$ depending only on $\beta^0$ such that for $\delta > 0$ small enough we have that 
\begin{eqnarray*}
\EE \sup_{\Vert \beta - \beta_0 \Vert_2 < \delta}\sqrt{n}\vert(\DD_{n, 2}(\beta)  - \mathcal D_{2, F^Y}(\beta)) - (\DD_{n, 2}(\beta_0)  -  \mathcal D_{2, F^Y}(\beta_0) )\vert\le C \left(\delta + \frac{\delta}{\sqrt{n}} + \frac{1}{\sqrt n}\right)\equiv\phi_n(\delta).
\end{eqnarray*}
\end{proposition}
\begin{proof}\textbf{(Proposition \ref{empcontrol}).} We begin with rewriting $\DD_{n, 2}(\beta)$ and $\mathcal{D}_{2, F^Y}(\beta)$ for a given $\beta \in \RR^d$. For $y \in \RR$ and $\beta \in \RR^d$ define the function
\begin{eqnarray*}
m_{\beta, y}(x) = F^\epsilon(y - \beta^T x),  \   x  \in \RR^d.
\end{eqnarray*}
Then,
\begin{eqnarray*}
\DD_{n, 2}(\beta) = \int\left(F^Y_n(y) - \EE^X_n m_{\beta, y} \right)^2 d F^Y_n(y).
\end{eqnarray*}
Since $F^Y(y) = \int F^\epsilon(y-\beta^T_0 x) dF^X(x) = \EE^X m_{\beta_0, y}$ for all $y \in \RR$, we have
\begin{eqnarray*}
\mathcal{D}_{2, F^Y}(\beta) = \int\left(F^Y(y) - \EE^X m_{\beta, y} \right)^2 d F^Y(y) =  \int\left(\EE^X (m_{\beta, y} - m_{\beta_0, y})\right)^2 d F^Y(y).
\end{eqnarray*}
Therefore,
\begin{eqnarray*}
\DD_{n, 2}(\beta) & = & \int\left(F^Y_n(y) - \EE^X_n m_{\beta_0, y} - \EE^X_n (m_{\beta, y} - m_{\beta_0, y})\right)^2 d F^Y_n(y)\\
& = & \int\left(F^Y_n(y) - \EE^X_n  m_{\beta_0, y}\right)^2 d F^Y_n(y)  \\
&& - \  2  \int\left(F^Y_n(y) - \EE^X_n m_{\beta_0, y}\right)\EE^X_n (m_{\beta, y} - m_{\beta_0, y}) d F^Y_n(y)  \\
&&  +  \  \int \left(\EE^X_n (m_{\beta, y} - m_{\beta_0, y})\right)^2 d F^Y_n(y)  \\
& = & \DD_{n, 2}(\beta_0) - 2 \int\left(F^Y_n(y) - \EE^X_n m_{\beta_0, y}  \right)\EE^X_n (m_{\beta, y} - m_{\beta_0, y}) d F^Y_n(y) \\
&& +  \  \int \left(\EE^X_n (m_{\beta, y} - m_{\beta_0, y})\right)^2 d F^Y_n(y)\\
& = & \DD_{n, 2}(\beta_0)  \\
&& - 2\int\left(F^Y_n(y) - \EE^X_n m_{\beta_0, y}\right) \left(\EE^X_n(m_{\beta, y} - m_{\beta_0, y}) - \EE^X(m_{\beta, y} - m_{\beta_0, y})\right)d F^Y_n(y)  \\
&& -\ 2 \int\left(F^Y_n(y) - \EE^X_n m_{\beta_0, y}\right)\EE^X (m_{\beta, y} - m_{\beta_0, y})d F^Y_n(y) \\
&& + \ \int \left(\EE^X_n (m_{\beta, y} - m_{\beta_0, y})\right)^2 d F^Y_n(y)
\end{eqnarray*}
with 
\begin{eqnarray*}
&& \int \left(\EE^X_n (m_{\beta, y} - m_{\beta_0, y})\right)^2 d F^Y_n(y) \\
&=&  \int \left[\left(\EE^X_n(m_{\beta, y} - m_{\beta_0, y}) - \EE^X(m_{\beta, y} - m_{\beta_0, y})\right) \right]^2 d\PP^Y_n(y)  \\
&& +  \  2 \int\left(\EE^X_n(m_{\beta, y} - m_{\beta_0, y})  -  \EE^X(m_{\beta, y} - m_{\beta_0, y}) \right) \EE^X\vert m_{\beta, y} - m_{\beta_0, y} \vert d F^Y_n(y)  \\
&& +  \  \int \left(\EE^X(m_{\beta, y} - m_{\beta_0, y})\right)^2 d F^Y_n(y)  \\
& = &  \int\left[\left(\EE^X_n(m_{\beta, y} - m_{\beta_0, y}) - \EE^X(m_{\beta, y} - m_{\beta_0, y})\right)\right]^2 d F^Y_n(y)  \\
&& +  \  2\int\left(\EE^X_n(m_{\beta, y} - m_{\beta_0, y}) - \EE^X(m_{\beta, y} - m_{\beta_0, y})\right)\EE^X (m_{\beta, y} - m_{\beta_0, y})d F^Y_n(y)  \\
&& + \  \mathcal{D}_{2, F^Y}(\beta) + \int\left(\EE^X (m_{\beta, y} - m_{\beta_0, y})\right)^2 d(F^Y_n(y) - F^Y(y))
\end{eqnarray*}
These together with $\mathcal{D}_{2, F^Y}(\beta_0) = 0$ gives
\begin{eqnarray*}
&& \sqrt{n}\left\vert(\DD_{n, 2}(\beta) - \mathcal{D}_{2, F^Y}(\beta)) - (\DD_{n, 2}(\beta_0) - \mathcal{D}_{2, F^Y}(\beta_0))\right\vert  \\
&& \le 2\sqrt{n}\int\left\vert F^Y_n(y) - \EE^X_n m_{\beta_0, y}\right \vert\left\vert(\EE^X_n(m_{\beta, y} - m_{\beta_0, y}) - \EE^X(m_{\beta, y} - m_{\beta_0, y}))\right\vert d F^Y_n(y)   \\
&& + \ 2\sqrt{n}\int\left\vert F^Y_n(y) - \EE^X_n m_{\beta_0, y}\right \vert\EE^X (m_{\beta, y} - m_{\beta_0, y})d F^Y_n(y)\\
&&  + \  \sqrt{n}\int\left[\left(\EE^X_n(m_{\beta, y} - m_{\beta_0, y}) - \EE^X(m_{\beta, y} - m_{\beta_0, y})\right)\right]^2 d F^Y_n(y)\\
&& + \ 2\sqrt{n}\int\left\vert(\EE^X_n(m_{\beta, y} - m_{\beta_0, y})- \EE^X(m_{\beta, y} - m_{\beta_0, y}))\right\vert\EE^X(m_{\beta, y} - m_{\beta_0, y})d F^Y_n(y)\\
&& +  \ \sqrt{n}\left\vert\int\left(\EE^X (m_{\beta, y} - m_{\beta_0, y})\right)^2 d(F^Y_n(y)- F^Y(y))\right\vert \\
&=& 2T_{n,1}(\beta) +  2T_{n,2}(\beta)  + T_{n,3}(\beta) + T_{n,4}(\beta)  + 2T_{n,5}(\beta).
\end{eqnarray*}
Since $F^Y(y) = \EE^X m_{\beta_0, y}$, the triangle inequality gives us
\begin{eqnarray*}
&& \sup_{\beta: \Vert \beta - \beta_0 \Vert_2 \le \delta} T_{n,1}(\beta) \\
&\le & \sqrt{n}\sup_{y \in \RR}\left\vert (F^Y_n - F^Y)(y) \right \vert \times  \\ 
&&\sup_{y \in \RR, \beta: \Vert \beta - \beta_0 \Vert_2 \le \delta} \left\vert(\EE^X_n(m_{\beta, y} - m_{\beta_0, y}) - \EE^X(m_{\beta, y} - m_{\beta_0, y})) \right \vert \\
&& \ +  \ \sqrt{n}\sup_{y \in \RR} \left \vert (\EE^X_n m_{\beta_0, y}  - \EE^X m_{\beta_0, y})\right\vert\times \\ &&\sup_{y \in \RR, \beta: \Vert \beta - \beta_0 \Vert_2 \le \delta} \left \vert (\EE^X_n(m_{\beta, y} - m_{\beta_0, y}) - \EE^X(m_{\beta, y} - m_{\beta_0, y})) \right \vert \\
&& = \frac{1}{\sqrt{n}} \left(\Vert \mathbb{G}^Y \Vert_{\mathcal{I}}   +   \Vert \mathbb{G}^X_n\Vert_{\mathcal{M}_{\beta_0}} \right)  \  \Vert \mathbb{G}^X_n \Vert_{\mathcal{G}}.  
\end{eqnarray*}
Then by Cauchy-Schwarz inequality
\begin{eqnarray*}
\EE\left[\sup_{\beta: \Vert\beta - \beta_0\Vert_2\le\delta} T_{n,1}(\beta)\right] &\le & \frac{1}{\sqrt{n}}\EE\left[\left(\Vert\mathbb{G}^Y \Vert_{\mathcal{I}} + \Vert\mathbb{G}^X_n\Vert_{\mathcal{M}_{\beta_0}} \right)  \  \Vert\mathbb{G}^X_n \Vert_{\mathcal{G}}\right]  \\
& \le &  \frac{1}{\sqrt{n}}\EE\left[\left(\Vert\mathbb{G}^Y \Vert_{\mathcal{I}} + \Vert\mathbb{G}^X_n\Vert_{\mathcal{M}_{\beta_0}} \right)^2 \right]^{1/2}\EE\left[\Vert\mathbb{G}^X_n \Vert^2_{\mathcal{G}}\right]^{1/2}  \\
& \lesssim &  \frac{\delta}{\sqrt n}.
\end{eqnarray*}
We also have that
\begin{eqnarray*}
\sup_{\beta: \Vert \beta - \beta_0 \Vert_2 \le \delta} T_{n,2}(\beta) &  \le & \left(\Vert\mathbb{G}^Y_n \Vert_{\mathcal{I}} + \Vert\mathbb{G}^X_n\Vert_{\mathcal{M}_{\beta_0}}\right)\sup_{y \in \RR, \beta: \Vert \beta - \beta_0 \Vert_2\le\delta}\EE^X\vert m_{\beta, y} - m_{\beta_0, y}\vert \\
& \le & M\delta\left(\Vert\mathbb{G}^Y_n\Vert_{\mathcal{I}} + \Vert\mathbb{G}^X_n\Vert_{\mathcal{M}_{\beta_0}}\right) \EE\Vert X \Vert_2 \\
& \lesssim & \delta
\end{eqnarray*}
since $\EE\Vert X\Vert^2_2 < \infty$ implies $\EE\Vert X\Vert_2 < \infty$. This implies that 
\begin{eqnarray*}
\EE\left[\sup_{\beta: \Vert\beta - \beta_0\Vert_2\le\delta} T_{n,2}(\beta)\right] & \lesssim &\left(\EE\Vert \mathbb{G}^Y_n \Vert_{\mathcal{I}} + \EE\Vert\mathbb{G}^X_n\Vert_{\mathcal{M}_{\beta_0}}\right)\delta \\
& \le & \left(\sqrt{\EE\Vert\mathbb{G}^Y_n \Vert^2_{\mathcal{I}}} + \sqrt{\EE\Vert\mathbb{G}^X_n\Vert^2_{\mathcal{M}_{\beta_0}}}\right)\delta  \\
& \lesssim & \delta.
\end{eqnarray*}
Now, 
\begin{eqnarray*}
\sup_{\beta:\Vert\beta - \beta_0\Vert_2\le\delta} T_{n,3}(\beta) &\le &  \frac{1}{\sqrt{n}}\sup_{y \in \RR}\sup_{\beta:\Vert\beta - \beta_0 \Vert_2 \le \delta} \left \vert \sqrt n (\EE^X_n(m_{\beta, y} - m_{\beta_0, y}) - \EE^X(m_{\beta, y} - m_{\beta_0, y}))\right\vert^2\\
& \le & \frac{1}{\sqrt{n}}\Vert\mathbb{G}^X_n  \Vert^2_{\mathcal{G}}
\end{eqnarray*}
and hence
\begin{eqnarray*}
\EE\left[\sup_{\beta: \Vert\beta - \beta_0\Vert_2\le\delta} T_{n,3}(\beta)\right] & \le & \frac{1}{\sqrt{n}}\EE\Vert\mathbb{G}^X_n  \Vert^2_{\mathcal{G}}\lesssim\frac{\delta^2}{\sqrt{n}}.
\end{eqnarray*}
We have that 
\begin{eqnarray*}
\sup_{\beta: \Vert\beta - \beta_0\Vert_2\le\delta} T_{n,4}(\beta) & \le  &   \Vert \mathbb{G}^X_n\Vert_{\mathcal{G}}  \  \sup_{y \in \RR, \beta: \Vert \beta - \beta_0 \Vert_2\le\delta}\EE^X\vert m_{\beta, y}  -  m_{\beta_0, y} \vert \\
& \lesssim & \delta\Vert\mathbb{G}^X_n\Vert_{\mathcal{G}} ,
\end{eqnarray*}
which implies
\begin{eqnarray*}
\mathbb E\left[\sup_{\beta: \Vert \beta - \beta_0 \Vert_2 \le \delta} T_{n,4}(\beta)\right]  \lesssim \delta\sqrt{\EE\left[\Vert \mathbb{G}^X_n\Vert^2_{\mathcal{G}}\right]} \lesssim \delta^2.
\end{eqnarray*}
Finally, to bound $T_{n, 5}$, note that 
\begin{eqnarray*}
T_{n, 5} \le\Vert\mathbb{G}^Y_n \Vert_{\mathcal{K}}
\end{eqnarray*}
where 
\begin{eqnarray*}
\mathcal{K} = \Big\{y \mapsto k_{\beta}(y) =  (\EE^X( m_{\beta, y} - m_{\beta_0, y}))^2: \Vert\beta - \beta_0\Vert_2\le\delta\Big\}
\end{eqnarray*}
More explicitly,
\begin{eqnarray*}
k_\beta(y) &= & \left(\int\left(m_{\beta, y}(x) - m_{\beta_0, y}(x)  \right)d F^X(x)\right)^2 \\
& = & \left(\int m_{\beta, y}(x)d F^X(x)\right)^2 + \left(\int m_{\beta_0, y}(x)  dF^X(x)\right)^2 \\
&& \ - 2\left(\int m_{\beta_0, y}(x)d F^X(x)\right)\left(\int m_{\beta, y}(x)d F^X(x)\right)\\
& = & k_{1, \beta}(y) + k_{2, \beta}(y) - 2 k_{3, \beta}(y).
\end{eqnarray*}
Functions $k_{i, \beta}$ are non-decreasing and $k_{i, \beta} \in [0,1]$ for $i =1, 2, 3$. Thus, $\mathcal{K}\subset\mathcal{C}$ defined in Proposition~\ref{sumsM}. For $\nu > 0$ define 
\begin{eqnarray*}
\tilde{J}(\nu) = \int_0^\nu \sqrt{1 + \log N_B(\eta, \mathcal{K}, L_2(F^Y))} d\eta. 
\end{eqnarray*}
Since $\EE k^2_\beta \lesssim \delta^4$ and $\Vert k_\beta \Vert_\infty \lesssim \delta^2$, it follows from Lemma 3.4.2 of~\cite{aadbook} that 
\begin{eqnarray*}
\EE\Vert\mathbb{G}^Y_n \Vert_{\mathcal{K}}\le\tilde{J}(\delta^2)\left(1 + \frac{\tilde{J}(\delta^2)}{\delta^4\sqrt{n}}\right).
\end{eqnarray*}
Now, by Proposition~\ref{sumsM} and for small enough $\nu > 0$ we have
\begin{eqnarray*}
\tilde{J}(\nu) & \le & \int_0^\nu\sqrt{1 + \log N_B(\eta, \mathcal{C}, L_2(F^Y)} d\eta  \\
& \le & \nu + \sqrt{A}\int_0^\nu\eta^{-1/2} d\eta = \mu + 2\sqrt{A}\sqrt{\nu} \lesssim \sqrt{\nu}.
\end{eqnarray*}
Thus,
\begin{eqnarray*}
\EE\Vert\mathbb{G}^Y_n\Vert_{\mathcal{K}} & \lesssim & \delta\left(1 + \frac{\delta}{\delta^4 \sqrt{n}}\delta^2\right) = \delta + \frac{1}{\sqrt{n}}.
\end{eqnarray*}
Taking $\delta\in(0, 1)$, there exists a constant $D > 0$ such that 
\begin{eqnarray*}
\EE\sup_{\Vert\beta - \beta_0\Vert_2 < \delta}\sqrt{n}\left\vert(\DD_{n, 2}(\beta) - \mathcal{D}_{2, F^Y}(\beta)) - (\DD_{n, 2}(\beta_0)  -  \mathcal{D}_{2, F^Y}(\beta_0))\right\vert & \le & D\left(3\delta + 2\frac{\delta}{\sqrt{n}} + \frac{1}{\sqrt{n}}\right) \\
& \lesssim & \delta + \frac{\delta}{\sqrt{n}} + \frac{1}{\sqrt{n}},
\end{eqnarray*}
which completes the proof. \hfill 
\end{proof}
\begin{proposition}\label{d2}
Assume that $F^\epsilon$ admits a density $f^\epsilon\in C^2$ such that  $f^\epsilon$, and $\vert (f^\epsilon)^{(1)} \vert$ are bounded by some constant $M > 0$, and $\EE\Vert X\Vert^2_2 < \infty$. Also, assume that matrix
\begin{eqnarray*}
U = \int\left(\int x f^\epsilon(y - \beta^T_0 x) dF^X(x)\right) \left(\int x^T f^\epsilon(y - \beta^T_0 x) dF^X(x)\right) dF^Y(y)  
\end{eqnarray*}
is positive definite. Then, there exists a small neighborhood of $\beta_0$ and a constant $C > 0$ such that for all $\beta$ in this  neighborhood 
\begin{eqnarray*}
\mathcal{D}_{2, F^Y}(\beta) \ge C \Vert\beta - \beta_0\Vert^2_2.
\end{eqnarray*} 
\end{proposition}
\begin{proof}\textbf{(Proposition~\ref{d2}).} \ We need to show that
\begin{eqnarray*}
\psi(\beta)\equiv\mathcal{D}_{2, F^Y}(\beta), \ \beta\in\RR^d 
\end{eqnarray*}
is twice continuously differentiable in a small neighbourhood of $\beta_0$ such that
\begin{eqnarray}\label{grad}
\frac{\partial\psi(\beta)}{\partial\beta}|_{\beta_0} = \mathbf{0},
\end{eqnarray}
and the Hessian matrix $\frac{\partial^2 \psi(\beta)}{\partial\beta \partial\beta^T}|_{\beta_0}$ is positive definite. We can write
\begin{eqnarray*}
\psi(\beta) = \int(\Delta(\beta, y))^2 dF^Y(y),
\end{eqnarray*}
where
\begin{eqnarray*}
\Delta(\beta, y) = \int F^\epsilon(y - \beta^T x) dF^X(x) - \int F^\epsilon(y - \beta^T_0 x) dF^X(x), \ (\beta, y) \in \RR^d\times\RR.
\end{eqnarray*}
By assumption, the distribution function $F^\epsilon $ is continuously differentiable on $\RR$, and hence 
\begin{eqnarray*}
\frac{\partial}{\partial\beta}F^\epsilon(y - \beta^T x) = -x f^\epsilon(y - \beta^T x).
\end{eqnarray*}
Furthermore, and by the assumption of the proposition, we have $\Vert x f^\epsilon(y-\beta^T x)\Vert_2\le M\Vert x \Vert_2$ and $\EE\Vert X \Vert_2 < \infty$. Hence, for any fixed $y$, $\beta\mapsto\Vert\frac{\partial}{\partial\beta} F^\epsilon(y - \beta^T x)  \Vert\le M\Vert x\Vert_2$ which is integrable. Thus, by the Dominated Convergence Theorem, it follows that for any fixed $y$
\begin{eqnarray*}
\frac{\partial  \Delta(\beta, y) }{\partial\beta} = -\int x f^\epsilon(y - \beta^T x) dF^X(x),
\end{eqnarray*}
and hence 
\begin{eqnarray*}
\frac{\partial(\Delta(\beta, y))^2 }{\partial\beta} = -2\Delta(\beta, y)  \int x f^\epsilon(y - \beta^T x) dF^X(x).
\end{eqnarray*}
Since $\vert F^\epsilon(y - \beta^T x) - F^\epsilon(y - \beta^T_0 x) \vert\le 1$, it follows that $\vert\Delta(\beta, y)\vert\le 1$ for all $\beta, y$. Thus, 
\begin{eqnarray*}
\left\Vert\frac{\partial (\Delta(\beta, y))^2}{\partial \beta}\right \Vert_2\le  2\left\Vert\Delta(\beta, y)\int x f^\epsilon(y - \beta^T x) dF^X(x)\right\Vert_2\le M\EE\Vert X\Vert_2.
\end{eqnarray*}
From the Dominated Convergence Theorem, it follows that for any $\beta\in\RR^d$
\begin{eqnarray*}
\frac{\partial\psi(\beta)}{\partial\beta} = -2\int\Delta(\beta, y)  \left(\int x f^\epsilon(y - \beta^T x) dF^X(x)\right) dF^Y(y).
\end{eqnarray*}
In particular, \eqref{grad} holds. Now, we show that the multivariate
\begin{eqnarray*}
f(\beta, y) =  -2\Delta(\beta, y)\left(\int x f^\epsilon(y - \beta^T x) dF^X(x)\right)
\end{eqnarray*}
is continuously differentiable on $\RR$ with a positive definite Hessian at $\beta_0$. Since we have already shown that $\beta\mapsto \Delta(\beta, y)$ is differentiable, we need to show differentiability of 
\begin{eqnarray}\label{Fun}
\beta \mapsto \int x f^\epsilon(y - \beta^T x) dF^X(x).
\end{eqnarray}
Using the assumption that $f^\epsilon$ is continuously differentiable such that $\vert (f^\epsilon)^{(1)}\vert\le M$, we can use arguments that are very similar to those above to show that the function in \eqref{Fun} is indeed differentiable with gradient 
\begin{eqnarray}\label{Fun}
\beta\mapsto -\int x x^T (f^\epsilon)^{(1)}(y - \beta^T x) dF^X(x).
\end{eqnarray}
Thus, for a fixed $y\in\RR$
\begin{eqnarray*}
\frac{\partial f(\beta, y)}{\beta^T} & = & 2\Delta(\beta, y)\left(\int x x^T (f^\epsilon)^{(1)}(y - \beta^T x) dF^X(x)\right)\\
&& \ + 2\left(\int x f^\epsilon(y - \beta^T x) dF^X(x)\right)\left(\int x^T f^\epsilon(y - \beta^T x) dF^X(x)\right).
\end{eqnarray*}
For $1 \le i, j \le d$, we have that
\begin{eqnarray}\label{Boundf}
&& \Big\vert\Delta(\beta, y)\left(\int [x x^T]_{ij} (f^\epsilon)^{(1)}(y - \beta^T x) dF^X(x)\right) \notag \\
&& \ + 2 \left(\int x_i f^\epsilon(y - \beta^T x) dF^X(x)\right)  \left(\int x_j f^\epsilon(y - \beta^T x) dF^X(x)\right)\Big\vert\notag \\
&& \le M\int\vert [x x^T]_{ij}\vert dF^X(x) + 2M\left(\int\vert x_i \vert dF^X(x)\right)\left(\int\vert x_j \vert dF^X(x)\right)\notag   \\
&& = M\int\vert x_i x_j\vert dF^X(x) + 2M \left(\int\vert x_i \vert dF^X(x)\right)\left(\int\vert x_j\vert dF^X(x)\right)\notag  \\
&& \le M\EE\Vert X\Vert^2_2 + 2M\left(\EE\Vert X\Vert_2\right)^2
\end{eqnarray}
using the Cauchy-Schwarz inequality and the fact that $\vert x_i \vert \le \Vert x \Vert_2$. The term on the right of \eqref{Boundf} is a constant and hence integrable with respect to $F^Y$. Using the Dominated convergence Theorem again, it follows that $\psi$ is twice continuously differentiable and 
\begin{eqnarray*}
\frac{\partial^2\psi(\beta)}{\partial\beta\partial\beta^T} & = & 2 \int  \Delta(\beta, y)\left(\int x x^T (f^\epsilon)^{(1)}(y - \beta^T x) dF^X(x)\right)dF^Y(y) \\
&&  \ + 2\int\left(\int x f^\epsilon(y - \beta^T x) dF^X(x)\right)  \left(\int x^T f^\epsilon(y - \beta^T x) dF^X(x)\right)dF^Y(y).
\end{eqnarray*}
Therefore we have
\begin{eqnarray*}
\frac{\partial\psi(\beta)}{\partial\beta\partial\beta^T}|_{\beta_0} & =  & 2\int\left(\int x f^\epsilon(y - \beta^T_0 x) dF^X(x)\right)  \left(\int x^T f^\epsilon(y - \beta_0^T x) dF^X(x)\right) dF^Y(y) = 2 U,
\end{eqnarray*}
which is positive definite by assumption. Using Taylor expansion of $\psi$ up to the second order, and using $\psi(\beta_0) = 0$ and \eqref{grad} we can write
\begin{eqnarray*}
\psi(\beta) & = & \psi(\beta_0) + (\beta - \beta_0)^T\frac{\partial \psi(\beta)}{\partial\beta}|_{\beta_0} + (\beta-\beta_0)^T U(\beta-\beta_0) + o(\Vert\beta - \beta_0\Vert^2_2)  \\
& = & (\beta-\beta_0)^T U (\beta-\beta_0) + o(\Vert \beta - \beta_0 \Vert^2_2). 
\end{eqnarray*} 
If $\lambda_1$ is the smallest eigenvalue of $U$, then $\lambda_1 > 0$ and also $ (\beta-\beta_0)^T U (\beta-\beta_0)\ge\lambda_1\Vert\beta - \beta_0\Vert^2_2$. Also, there exist a small enough neighborhood of $\beta_0$ such that when $\beta$ is in this neighborhood we have $o(\Vert \beta - \beta_0 \Vert^2_2) > -\lambda_1\Vert\beta - \beta_0\Vert^2_2 /2$.  This implies that for $\beta$ in this neighbourhood  
\begin{eqnarray*}
\psi(\beta)  = \mathcal{D}_{2, F^Y}(\beta)\ge\frac{\lambda_1}{2}\Vert \beta - \beta_0 \Vert^2_2
\end{eqnarray*}
which completes the proof. \hfill 
\end{proof}

\begin{theorem}\label{MainTheo}
Let $K > 0$,  $h\in B(\mathbf{0}, K) = \{h\in\RR^d: \Vert h\Vert_2\le K\}$, and  $\mathbb{B}_1, \mathbb{B}_2$ two independent standard Brownian Bridges from $(0,0)$ to $(1, 0)$. Assume that $F^\epsilon$ admits a density $f^\epsilon$ which 
is continuously differentiable and that $f^\epsilon$, and $\vert (f^\epsilon)^{(1)}$ are bounded by $M > 0$. There exist an integer $m \ge 1$ and $a_1 <  a_2 < \ldots < a_m$ such that $f^\epsilon$ is monotone on $(-\infty, a_1]$, $(a_i, a_{i+1}), i =1, \ldots, m-1$ and non-increasing on $[a_m, \infty)$ and $\EE\Vert X\Vert^2_2 < \infty$. Let
\begin{eqnarray*}
\mathbb{Q}_n(h) =  n\DD_{n, 2}(\beta_0 +  h n^{-1/2}).
\end{eqnarray*}  
Then,
\begin{eqnarray*}
\mathbb{Q}_n  \Rightarrow  \mathbb{Q}
\end{eqnarray*}
in $\ell^\infty(B(\mathbf{0}, K))$, where
\begin{eqnarray*}
&& \mathbb{Q}(h) = \\
&&\int\left\{\mathbb{B}_1 \circ F^Y(y)  -  \int F^\epsilon(y - \beta^T_0 x) d\mathbb{B}_2 \circ F^X(x) + h^T \int x f^\epsilon(y - \beta^T_0 x) dF^X(x)\right\}^2 d F^Y(y) \\
&&= \int\left\{\mathbb{B}_1 \circ F^Y(y) - \int\mathbb{B}_2 \circ F^Z(z) f^\epsilon(y - z)  dz  + h^T \int x f^\epsilon(y - \beta^T_0 x) dF^X(x) \right\}^2 d F^Y(y) 
\end{eqnarray*}
with $Z = \beta^T_0 X$. Furthermore, if the $d \times d$ matrix 
\begin{eqnarray*}
U  =  \int \left(\int x f^\epsilon(y - \beta^T_0 x) dF^X(x)\right)\left(\int x^T f^\epsilon(y - \beta^T_0 x) dF^X(x)\right) dF^Y(y)
\end{eqnarray*}
is positive definite then $h^0  =  U^{-1} V$ is the unique minimizer of the process $\mathbb{Q}$ where
\begin{eqnarray*}
V & = & -\int\left (\mathbb{B}_1 \circ F^Y(y)  -  \int F^\epsilon(y - \beta^T_0 x ) d\mathbb{B}_2\circ F^X(x)\right)\int\left(x f^\epsilon(y - \beta^T_0 x) dF^X(x)\right) dF^Y(y) \\
& = &  -\int \left (\mathbb{B}_1 \circ F^Y(y)  -  \int\mathbb{B}_2 \circ F^Z(z)  f^\epsilon(y - z) dz \right)\int\left(x f^\epsilon(y - \beta^T_0 x) dF^X(x) \right) dF^Y(y)
\end{eqnarray*}
\end{theorem}

\begin{proof}\textbf{(Therem \ref{MainTheo}).}
Let $h\in B(\mathbf{0}, K)$ for a fixed $K > 0$. Using Taylor expansion, we have
\begin{eqnarray*}
F^\epsilon(y - (\beta_0 + h n^{-1/2})^T x) &=& F^\epsilon(y - \beta^T_0  x -  n^{-1/2}h^T x) \\   
&=&  F^\epsilon(y - \beta^T_0 x) - n^{-1/2} h^T x f^\epsilon(y - \beta^T_0x) + \frac{1}{2} n^{-1}(h^T x)^2(f^\epsilon)^{(1)}(\theta^*_n),
\end{eqnarray*}
for some $\theta^*_n = \theta^*(y, x, n)$ between $y - \beta^T_0 x$ and $y - \beta^T_0 x - n^{-1/2} h^T x$. Furthermore, we shall use the well-known fact that if $Z_1, \ldots, Z_n$ are i.i.d. $\sim F$ then there exists a Brownian Bridge $\mathbb{B}$ defined on the same probability space as $(Z_1, \ldots, Z_n)$ such that 
\begin{eqnarray*}
\Vert\mathbb{B}_n - \mathbb{B}\Vert_\infty \to 0
\end{eqnarray*}
almost surely, where $\mathbb{B}_n$ is defined such that for all for all $z \in \RR$ we have
\begin{eqnarray*}
F_n(z) = n^{-1}\sum_{i=1}^n\one\{Z_i \le z\} = F(z) + n^{-1/2}\mathbb{B}_n\circ F(z).
\end{eqnarray*}
Now, 
\begin{eqnarray*}
\mathbb{Q}_n(h) = \int\mathbb{K}^2_n(y, h)d F^Y_n(y)
\end{eqnarray*}
where 
\begin{eqnarray*}
\mathbb{K}^2_n(y, h) &=& n\Big\{F^Y_n(y) - \int F^\epsilon(y - \beta^T_0 x)d F^X_n(x)  + n^{-1/2}\int h^T x f^\epsilon(y - \beta^T_0x) d F^X_n(x) \\
&&  \  \ - \frac{1}{2} n^{-1}\int(h^T x)^2(f^\epsilon)^{(1)}(\theta^*_n)d F^X_n(x)  \Big\}^2 \\
&=& \Big \{\sqrt n (F^Y_n(y) - F^Y(y)) - \int F^\epsilon(y - \beta^T_0 x)d\sqrt n(F^X_n(x) - F^X(x))\\
&& \ + \int h^T x f^\epsilon(y - \beta^T_0x) d F^X_n(x) - \frac{1}{2} n^{-1/2} \int  (h^T x)^2 (f^\epsilon)^{(1)}(\theta^*_n)d F^X_n(x)\Big\}^2.
\end{eqnarray*}
Using again $F^Y(y)  = \int F^\epsilon(y - \beta_0^Tx) dF^X(x)$
\begin{eqnarray*}
\mathbb{K}^2_n(y, h) &=& \Big \{\mathbb{B}_{n,1}(F^Y(y)) - \int\mathbb{B}_{n, 2}(F^Z(z))  f^\epsilon(y - z)dz +\int h^T x f^\epsilon(y - \beta^T_0x) d F^X_n(x)  \\
&&  \ - \frac{1}{2} n^{-1/2}\int (h^T x)^2(f^\epsilon)^{(1)}(\theta^*_n)d F^X_n(x)  \Big\}^2, 
\end{eqnarray*}
where $Z=\beta^T_0 X$, and $\mathbb{B}_{n,1}\circ F^Y(y) = \sqrt{n}(F^Y_n(y) - F^Y(y)$, and $\mathbb{B}_{n,2}\circ F^X(y) = \sqrt{n}(F^X_n(y) - F^X(y))$.

Now let $Z = \beta^T_0 X$ and $z = \beta^T_0 x$. Then,
\begin{eqnarray*}
K_n(y) &=& n \Big\{F^Y_n(y) + \int f^\epsilon(y - z) F^Z_n(z) + n^{-1/2}\int h^T x f^\epsilon(y - \beta^T_0x) d F^X_n(x) \\
&&  \  \ - \frac{1}{2} n^{-1} \int(h^T x)^2(f^\epsilon)^{(1)}(\theta^*_n)d F^X_n(x)  \Big\}^2.
\end{eqnarray*}
Define now the function
\begin{eqnarray*}
\mathbb{K}^2(y, h)  =  \left(\mathbb{B}_1 \circ F^Y(y) - \int\mathbb{B}_2 \circ F^Z(z)f^\epsilon(y - z) dz  +  \int h^T x f^\epsilon(y - \beta^T_0x) dF^X(x)   \right)^2
\end{eqnarray*}
for $y\in\RR$, where $\mathbb{B}_1$ and $\mathbb{B}_2$ are two independent standard Brownian Bridges from $(0,0)$ to $(1,0)$ defined on the same probability space as $(X_1, \ldots, X_n)$ and $(Y_1, \ldots, Y_n)$ such that $\Vert\mathbb{B}_{n,1}  - \mathbb{B}_1\Vert_\infty \to 0$ and $\Vert\mathbb B_{n,2} - \mathbb{B}_2 \Vert_\infty  \to 0$ almost surely. We show now that 
\begin{eqnarray*}
\sup_{y\in\RR, \Vert h \Vert_2 \le K} \vert \mathbb{K}^2_n(y, h) - \mathbb{K}^2(y, h)  \vert = o_{\PP}(1).
\end{eqnarray*}
Let $E_1 = \mathbb{B}_{n,1} \circ F^Y - \mathbb{B}_1 \circ F^Y$ and $E_2 = \mathbb{B}_{n,2} \circ F^Z - \mathbb{B}_2 \circ F^Z$. Then,
\begin{eqnarray*}
\mathbb{K}^2_n(y, h) &=& \Big\{\mathbb{B}_1 \circ F^Y(y) - \int\mathbb{B}_2 \circ F^Z(z)f^\epsilon(y - z) dz + \int h^T x f^\epsilon(y - \beta^T_0x) d F^X(x) \\
&& \ + \ E_1(y) - \int E_2(z) f^\epsilon(y - z) dz + n^{-1/2} \int h^T x f^\epsilon(y - \beta^T_0x) d \mathbb{G}^X_n(x)\\
&&  \ - \frac{1}{2} n^{-1/2} \int (h^T x)^2 (f^\epsilon)^{(1)}(\theta^*_n) d F^X_n(x)  \Big\}^2\\
& = & \mathbb{K}^2(y, h) + \mathbb{K}(y, h)R_n(y, h) + R^2_n(y, h)
\end{eqnarray*}
where
\begin{eqnarray*}
R_n(y,h) &=& E_1(y) - \int E_2(z) f^\epsilon(y - z)dz + n^{-1/2} \int h^T x f^\epsilon(y - \beta^T_0x) d\mathbb{G}^X_n(x) \\
&& \  - \frac{1}{2} n^{-1/2} \int(h^T x)^2 (f^\epsilon)^{(1)}(\theta^*_n) d F^X_n(x).  
\end{eqnarray*}
Since $\vert\mathbb{K}(y, h)\vert\le 2 + K  M \int\Vert\Vert x \Vert_2 dF^X(x) \equiv C$ we have 
\begin{eqnarray*}
\vert\mathbb{K}^2_n(y, h) - \mathbb{K}^2(y, h)\vert\le C\vert R_n(y, h)\vert +  R^2_n(y, h)
\end{eqnarray*}
where
\begin{eqnarray*}
\vert R_n(y, h) \vert &\le & \Vert\mathbb{B}_{n,1} - \mathbb{B}_1 \Vert_\infty + \Vert\mathbb{B}_{n,2} - \mathbb{B}_2 \Vert_\infty \\
&& \ +  n^{-1/2}\left\vert\int h^T x f^\epsilon(y - \beta^T_0x) d\mathbb{G}^X_n(x)   \right\vert + \frac{M}{2} n^{-1/2} \int(h^T x)^2 d F^X_n(x).
\end{eqnarray*}
We show now that $\sup_{\Vert h\Vert_2 \le K, y \in \RR}\left\vert h^T \int x f^\epsilon(y - \beta^T_0x) d\mathbb{G}^X_n(x)\right\vert = O_{\PP}(1)$. As it is easy to adapt the proof of Proposition~\ref{pwmon} for any integer $m \ge 1$, we have that
\begin{eqnarray*}
\sup_{y\in\RR}\left\Vert\int x f^\epsilon(y - \beta^T_0x)d \mathbb{G}^X_n(x) \right \Vert_2 = O_{\PP}(1).
\end{eqnarray*}
Now, for any $h\in\RR^d:\Vert h\Vert_2\le K$ we have that  
\begin{eqnarray*}
\left\vert\int h^T x f^\epsilon(y - \beta^T_0x) d \mathbb{G}^X_n(x)\right\vert &=&  \left\vert h^T \int x f^\epsilon(y - \beta^T_0x) d \mathbb{G}^X_n(x)\right\vert \\
& \le & \Vert h\Vert_2\left\Vert\int x f^\epsilon(y - \beta^T_0x) d \mathbb{G}^X_n(x) \right\Vert_2, \\
&\le &  K\sup_{y}\left\Vert\int x f^\epsilon(y - \beta^T_0x)d\mathbb{G}^X_n(x)\right\Vert_2 \\
&=& O_{\PP}(1),
\end{eqnarray*}
where the last inequality follows from Cauchy-Schwarz. It follows that 
\begin{eqnarray*}
n^{-1/2}\sup_{h \in B(\mathbf{0}, K), y\in\RR}\left \vert \int h^T x f^\epsilon(y - \beta^T_0 x) d\mathbb{G}^X_n(x)\right\vert = O_{\PP}(n^{-1/2}).
\end{eqnarray*}
To handle the last term, we apply the SLLN to show that
\begin{eqnarray*}
\int (h^T x)^2 d F^X_n(x)\to\int(h^T x)^2 dF^X(x) < \infty,
\end{eqnarray*}
since $\int (h^T x)^2 dF^X(x)\le K^2 \vert\Vert x \Vert^2_2 dF^X(x) < \infty$ by assumption. It follows that
\begin{eqnarray*}
\left\vert\mathbb{Q}_n(h) - \int\mathbb{K}^2(y, h)d F^Y_n(y)\right\vert &=& \left\vert\int\mathbb K^2_n(y, h) d F^Y_n(y) - \int\mathbb{K}^2(y, h)d F^Y_n(y) \right\vert  \\
&\le &\sup_{y \in\RR, \Vert h \Vert_2 \le K}\vert\mathbb{K}^2_n(y, h) - \mathbb{K}^2(y, h)\vert\to 0
\end{eqnarray*}
in probability. Now, we show that 
\begin{eqnarray*}
\int\mathbb{K}^2(y, h)d F^Y_n(y)\Rightarrow\int\mathbb{K}^2(y, h)dF^Y(y) 
\end{eqnarray*}
in $\ell^\infty(B(\mathbf{0}, K))$. Here, we will apply Theorem 1.5.4 of~\cite{aadbook}. First, we show that the process 
\begin{eqnarray*}\label{Process}
\left(\int\mathbb{K}^2(y, h) d F^Y_n(y)\right)_{n, h \in B(\mathbf{0}, K)}\equiv   \left(Z_{n, h}\right)_{n, h \in B(\mathbf{0}, K)}
\end{eqnarray*} 
is asymptotically tight.  We have that for all $y \in\RR, h\in B(\mathbf{0}, K)$
\begin{eqnarray*}
K^2(h, y)\le\left(\Vert\mathbb{B}_1 \Vert_\infty + \Vert\mathbb{B}_2 \Vert_\infty +  K M \int\Vert x \Vert_2 dF^X(x)\right)^2.  
\end{eqnarray*}
Using a classical result about Brownian Bridge, it follows that for any $\delta > 0$ we can find $D_\delta > 0$ such that 
\begin{eqnarray*}
\PP(\Vert\mathbb{B}_1 \Vert_\infty +  \Vert\mathbb{B}_2 \Vert_\infty \le D_\delta) \ge 1-\delta.
\end{eqnarray*}
This implies that with $D'_\delta = (D_\delta  + K  M \int \Vert x \Vert_2 dF^x(x))^2$
\begin{eqnarray*}
\PP\left(\left(Z_{n, h}\right)_{n, h \in B(\mathbf{0}, K)}\in [0, D'_\delta]\right)  \ge 1-\delta
\end{eqnarray*}
Next, denote $Z_h = \int\mathbb{K}^2(y, h)dF^Y(y) $. We will show that for any finite integer $k \ge 1$, and  $h_1, \ldots, h_k \in B(\mathbf{0}, K)$, 
\begin{eqnarray*}
\left(Z_{n, h_1}, \ldots, Z_{n, h_k} \right) \stackrel{d}{\rightarrow} \left(Z_{h_1}, \ldots, Z_{h_k} \right) 
\end{eqnarray*}
as $n \to \infty$. We will show the stronger result that the convergence above occurs almost surely. Let $\alpha_1, \ldots, \alpha_k \in \RR$. The distribution function $F^Y$ is continuous on $\RR$ and hence 
\begin{align}\label{target}
& y \mapsto  \sum_{j=1}^k \alpha_j \mathbb{K}^2(y, h_j)\notag \\
& =\sum_{j=1}^k \alpha_j \left( \mathbb{B}_1 \circ F^Y(y)   +  \int \mathbb B_2 \circ F^Z(z) f^\epsilon(y - z)  dz + \int h_j^T x f^\epsilon(y - \beta^T_0x) dF^X(x)\right)^2
\end{align}
is continuous. In fact, it is known that the Brownian Bridge is continuous on $[0,1]$.  Also, the function 
\begin{eqnarray*}
y \mapsto\int \mathbb{B}_2 \circ F^Z(z) f^\epsilon(y - z)dz = \int\mathbb{B}_2 \circ F^Z(y-z) f^\epsilon(z) dz 
\end{eqnarray*}
is continuous at every point $y\in\RR$ by the dominated convergence Theorem: The function $y \mapsto\mathbb{B}_2 \circ F^Z(y-z) $ is continuous and 
\begin{eqnarray*}
\left\vert\mathbb{B}_2 \circ F^Z(y-z) f^\epsilon(z)\right\vert\le\Vert\mathbb{B}_2 \Vert_\infty f^\epsilon(z)
\end{eqnarray*}
is integrable. By the same theorem, we show that the function $y \mapsto \int h_j^T x f^\epsilon(y - \beta^T_0x) dF^X(x)$ is continuous: $y \mapsto f^\epsilon(y - \beta^T_0x)$ is continuous and $\left \vert h_j^T x f^\epsilon(y - \beta^T_0x) \right \vert \le \Vert h \Vert_2 \Vert x \Vert_2 M$ is integrable with respect to $F^X$.  Also, the function in (\ref{target}) is bounded  by
\begin{eqnarray*}
\sum_{j=1}^k \vert \alpha_j \vert \left( \Vert \mathbb{B}_1 \Vert_\infty + \Vert \mathbb{B}_2 \Vert_\infty  + \Vert h_j \Vert_2 M \int \Vert x \Vert_2 dF^X(x) \right)^2.
\end{eqnarray*}
Since $F^Y_n$ converges weakly to $F^Y$ almost surely and the continuity of the function,  it follows that
\begin{eqnarray*}
\sum_{j=1}^k \alpha_j Z_{n, h_j}\stackrel{a.s.}{\rightarrow}\sum_{j=1}^k \alpha_j Z_{h_j}.
\end{eqnarray*}
Since $\alpha_1, \ldots, \alpha_k$ are arbitrary it follows that  
\begin{eqnarray*}
\left(Z_{n, h_1}, \ldots, Z_{n, h_k} \right)\stackrel{a.s.}{\rightarrow}\left(Z_{h_1}, \ldots, Z_{h_k} \right) 
\end{eqnarray*}
as $n\to\infty$. This completes the proof that 
\begin{eqnarray*}
\mathbb{Q}_n \Rightarrow \mathbb{Q}
\end{eqnarray*}
in $\ell^\infty(B(\mathbb{0}, K))$. Lastly, let us put 
\begin{eqnarray*}
\mathbb{W}(y) = \mathbb{B}_1 \circ F^Y(y) - \int\mathbb{B}_2 \circ F^Z(z) f^\epsilon(y - z)  dz, \  \ \text{and} \ \  S(y)= \int x f^\epsilon(y - \beta^T_0 x) dF^X(x).
\end{eqnarray*}
Then, we can write
\begin{eqnarray*}
\mathbb{Q}(h) &=& \int\left(\mathbb{W}(y) +  h^T S(y)\right)^2 dF^Y(y)   \\
&=& \int\left(\mathbb{W}(y) +  h^T S(y)\right)\left(\mathbb{W}(y) + S(y)^T h\right)   dF^Y(y)  \\
&=& \int\mathbb{W}(y)^2  dF^Y(y) + 2\int\mathbb{W}(y) S(y)^T h dF^Y(y) +  h^T \int  S(y)S(y)^T dF^Y(y) h
\end{eqnarray*}
which is a quadratic form. It admits a unique minimum if the matrix
\begin{eqnarray*}
\int  S(y)S(y)^T dF^Y(y)  =  \int \left(x f^\epsilon(y - \beta^T_0 x) dF^X(x)\right) \left(x f^\epsilon(y - \beta^T_0 x) dF^X(x)\right)^T dF^Y(y)  = U
\end{eqnarray*}
is positive definite, in which  case the unique minimizer $h^0$ solves the equation
\begin{eqnarray*}
U h^0   = -\int \mathbb{W}(y) S(y)dF^Y(y) = V,
\end{eqnarray*}
or equivalently $h^0 = U^{-1} V$.  \hfill
\end{proof}

\subsection{Proof of Theorem~\ref{limitDist}} 
\begin{proof}
From Proposition~\ref{empcontrol} we know that there exists a constant $C > 0$ such that  
\begin{eqnarray*}
\phi_n(\delta)\le C\left(\delta + \frac{\delta}{\sqrt{n}} + \frac{1}{\sqrt{n}}\right) 
\end{eqnarray*}
for all $\delta \in (0,1)$. Now, $\delta\mapsto\phi_n(\delta)/\delta$ is decreasing. Also, if we take $r_n = \sqrt{n}$ we get 
\begin{eqnarray*}
r_n^2 \phi_n(1/r_n) = C \ n \left(\frac{1}{\sqrt{n}} + \frac{1}{n} + \frac{1}{\sqrt{n}}\right) =  C \left( 2 \sqrt{n} + \frac{1}{\sqrt{n}}\right)\le 4 C \sqrt{n}
\end{eqnarray*}
for $n$ large enough. Combining this with Proposition~\ref{d2} it follows from Theorem 3.2.5 of~\cite{aadbook} that 
\begin{eqnarray*}
\sqrt{n}(\hat{\beta}_n - \beta_0) = O_{\PP}(1).
\end{eqnarray*}
\hfill
\end{proof}
\subsection{Proof of Theorem~\ref{orderedBeta}}
\begin{proof}
The arguments are the same as those used in the proof of Theorem~\ref{limitDist} except that we replace $X$ by $\tilde{X}$ as defined in the statement of the theorem.
\end{proof}
\subsection{Proof of Proposition~\ref{improved}}
\begin{proof}
Let $\beta_{m}^0 = \tilde{\beta}_m \frac{\Vert\beta_0\Vert_{2, \widehat{\Sigma}}}{\Vert\tilde{\beta}_m\Vert_{2, \widehat{\Sigma}}}$. Note that with probability one, we have 
\[
\Vert\tilde{\beta}_m - \beta_0\Vert \geq \Vert\beta_{m}^0 - \beta_0\Vert.
\]
We can rewrite $\tilde{\beta}_{n, m}$ as
\[
\tilde{\beta}_{n, m} = \beta_{m}^0 \frac{r_n}{\Vert\beta_0\Vert_{2, \widehat{\Sigma}}}.
\]
We have 
\begin{eqnarray*}
\frac{\EE[\Vert\tilde{\beta}_{n, m} - \beta_0\Vert_{2, \widehat{\Sigma}}^2]}{\EE[\Vert\tilde{\beta}_m - \beta_0\Vert_{2, \widehat{\Sigma}}^2]} &=& \frac{\EE[\Vert\beta_{m}^0 \frac{r_n}{\Vert\beta_0\Vert_{2, \widehat{\Sigma}}} - \beta_0\Vert_{2, \widehat{\Sigma}}^2]}{\EE[\Vert\tilde{\beta}_m - \beta_0\Vert_{2, \widehat{\Sigma}}^2]} \\
&=& \frac{\EE[\Vert\beta_{m}^0 - \beta_0 +  \beta_{m}^0(\frac{r_n}{\Vert\beta_0\Vert_{2, \widehat{\Sigma}}} - 1)\Vert_{2, \widehat{\Sigma}}^2]}{\EE[\Vert\tilde{\beta}_m - \beta_0\Vert_{2, \widehat{\Sigma}}^2]} \\
&=& \frac{\EE[\Vert\beta_{m}^0 - \beta_0\Vert_{2, \widehat{\Sigma}}^2] + 2T_2 + T_3}{\EE[\Vert\tilde{\beta}_m - \beta_0\Vert_{2, \widehat{\Sigma}}^2]},
\end{eqnarray*}
where
\[
T_2 = \EE[(\frac{r_n - \Vert\beta_0\Vert_{2, \widehat{\Sigma}}}{\Vert\beta_0\Vert_{2, \widehat{\Sigma}}})(\beta_m^0 - \beta_0)^T\widehat{\Sigma}\beta_m^0
\]
and 
\[
T_3 = \EE[\frac{(r_n - \Vert\beta_0\Vert_{2, \widehat{\Sigma}})^2}{\Vert\beta_0\Vert_{2, \widehat{\Sigma}}^2}].
\]
To complete the proof, it is enough to show that there exist constants $C_1, C_2$, and $C_3$ such that 
\begin{eqnarray}\label{eq:olsBound}
\frac{C_1}{m} \leq \EE[\Vert\tilde{\beta}_m - \beta_0\Vert_{2, \widehat{\Sigma}}^2] \leq \frac{C_2}{m},
\end{eqnarray}
and 
\begin{eqnarray}\label{eq:T3}
T_3 \leq \frac{C_3}{n}.
\end{eqnarray}
With~\eqref{eq:olsBound}, and~\eqref{eq:T3} we have 
\[
\frac{2T_2}{\EE[\Vert\hat{\beta}_m - \beta_0\Vert_{2, \widehat{\Sigma}}^2]} \lesssim \sqrt{\frac{m}{n}},
\]
and
\[
\frac{T_3}{\EE[\Vert\hat{\beta}_m - \beta_0\Vert_{2, \widehat{\Sigma}}^2]} \lesssim \frac{m}{n}
\]
and therefore, we can conclude that there exists $D_1, D_2 > 0$ such that 
\begin{eqnarray*}
\frac{\EE[\Vert\tilde{\beta}_{n, m} - \beta_0\Vert_{2, \widehat{\Sigma}}^2]}{\EE[\Vert\tilde{\beta}_m - \beta_0\Vert_{2, \widehat{\Sigma}}^2]} &\leq & \frac{\EE[\Vert\beta_m^0 - \beta_0\Vert_{2, \widehat{\Sigma}}^2]}{\EE[\Vert\tilde{\beta}_m - \beta_0\Vert_{2, \widehat{\Sigma}}^2]} + D_1\sqrt{\frac{m}{n}} + D_2\frac{m}{n} \\
& < & 1 +  D_1\sqrt{\frac{m}{n}} + D_2\frac{m}{n}.
\end{eqnarray*}
Therefore for $m = o(n)$ we have 
\[
\frac{\EE[\Vert\tilde{\beta}_{n, m} - \beta_0\Vert_{2, \widehat{\Sigma}}^2]}{\EE[\Vert\tilde{\beta}_m - \beta_0\Vert_{2, \widehat{\Sigma}}^2]} < 1.
\]
To complete the proof, we need to prove~\eqref{eq:olsBound}, and~\eqref{eq:T3}. Let's start by proving~\eqref{eq:olsBound}. First note that 
\begin{eqnarray*}
\tilde{\beta}_m = \beta_0 + [\tilde{X}^T\tilde{X}]^{-1}\tilde{X}^T\tilde{\epsilon},
\end{eqnarray*}
where $\tilde{X}$ is a $m\times d$ matrix such that its rows are $\tilde{X}_i$ from the matched data. Also $\tilde{\epsilon} = [\tilde{\epsilon}_1, \cdots, \tilde{\epsilon}_m]$ with $\tilde{\epsilon}_i = \tilde{Y}_i - \beta_0^T\tilde{X}_i$. Let $M = [\tilde{X}^T\tilde{X}]^{-1}\tilde{X}^T$ then it follows that 
\begin{eqnarray*}
\Vert\tilde{\beta}_m - \beta_0\Vert_{2, \widehat{\Sigma}}^2 &=& \tilde{\epsilon}M^T\widehat{\Sigma}M\tilde{\epsilon} \\
&=& tr(\tilde{\epsilon}M^T\widehat{\Sigma}M\tilde{\epsilon}) \\
&=& tr(\widehat{\Sigma}M\tilde{\epsilon}\tilde{\epsilon}M^T),
\end{eqnarray*}
where in the above equality, we are using the fact about trace that $tr(AB) = tr(BA)$. Now let $A = M\tilde{\epsilon}\tilde{\epsilon}^TM^T$, then  we have
\begin{eqnarray*}
tr(\widehat{\Sigma}A) &=& \sum_{i=1}^d(\widehat{\Sigma}A)_{ii} \\
&=& \sum_{i=1}^d\sum_{j=1}^d \widehat{\Sigma}_{ij}A_{ji}.
\end{eqnarray*}
Note that $A$ is based on the matched sample $\{(\tilde{X}_i, \tilde{Y}_i)\}_{i=1}^m$ and $\widehat{\Sigma}$ is computed using the unmatched sample $\{Y_i\}_{i=1}^n$ and $\{X_i\}_{i=1}^n$. Therefore, $A$ and $\widehat{\Sigma}$ are independent of each other. This gives us
\begin{eqnarray*}
\EE[\Vert\tilde{\beta}_m - \beta_0\Vert_{2, \widehat{\Sigma}}^2] &=& \sum_{i=1}^d\sum_{j=1}^d \EE[\widehat{\Sigma}_{ij}]\EE[A_{ji}] \\
&=& \sum_{i=1}^d\sum_{j=1}^d \Sigma_{ij}\EE[A_{ji}] \\
&=& tr(\Sigma\EE[A])
\end{eqnarray*}
According to \cite{KleinmanAthans} we have 
\begin{eqnarray*}
\lambda_d(\Sigma)tr(\EE[A])\leq tr(\Sigma\EE[A]) \leq \lambda_1(\Sigma)tr(\EE[A]),
\end{eqnarray*}
where $\lambda_d(\Sigma)$ and $\lambda_1(\Sigma)$ are respectively the smallest and largest eigenvalues of $\Sigma$. To bound $tr(\EE[A]) = \EE[tr(A)]$, let $K = (\tilde{\epsilon}\tilde{\epsilon}^T)$ and $P = M^TM$ and note that 
\begin{eqnarray*}
\EE[tr(A)] &=& \EE[tr(M\tilde{\epsilon}\tilde{\epsilon}^TM^T)] \\
&=& \EE[tr(\tilde{\epsilon}\tilde{\epsilon}^TM^TM)] \\
&=& \EE[tr(KP)] \\
&=& \sum_{i=1}^m\sum_{j=1}^m\EE[K_{ij}P_{ji}] \\
&=& \var(\epsilon)\EE[tr(P)] \\
&=& \var(\epsilon)\EE[tr([\tilde{X}^T\tilde{X}]^{-1})] \\
&=& \frac{\var(\epsilon)}{m}\EE[tr([\frac{\tilde{X}^T\tilde{X}}{m}]^{-1})].
\end{eqnarray*}
Therefore 
\begin{eqnarray*}
\frac{\var(\epsilon)}{m} \EE[\frac{1}{\lambda_1(\frac{\tilde{X}^T\tilde{X}}{m})}]\leq\EE[tr(A)]\leq d\frac{\var(\epsilon)}{m}\EE[\frac{1}{\lambda_d(\frac{\tilde{X}^T\tilde{X}}{m})}].
\end{eqnarray*}
By Weyl's theorem, we have
\begin{eqnarray*}
\max_{j = 1, \cdots, d}|\lambda_j(\Sigma) - \lambda_j(\frac{1}{m}\tilde{X}^T\tilde{X})| \leq \Vert\Sigma - \frac{1}{m}\tilde{X}^T\tilde{X}\Vert_2.
\end{eqnarray*}
On the other hand, with high probability, we have 
\begin{eqnarray*}
\Vert\Sigma - \frac{1}{m}\tilde{X}^T\tilde{X}\Vert_2 \leq \lambda_1(\Sigma)(2\sqrt{\frac{d}{m}} + 2\delta + (\sqrt{\frac{d}{m}} + \delta)^2)
\end{eqnarray*}
for all $\delta > 0$. Putting these together gives us a high probability 
\begin{eqnarray*}
\frac{\var(\epsilon)}{m} (\lambda_1(\Sigma) - C_1\sqrt{\frac{d}{m}})^{-1}\leq\EE[tr(A)]\leq d\frac{\var(\epsilon)}{m}(\lambda_d(\Sigma) + C_2\sqrt{\frac{d}{m}})^{-1},
\end{eqnarray*}
where $C_1, C_2$ are constants which depends on $\lambda_1(\Sigma)$ and $\delta$. This gives us 
\begin{eqnarray*}
\EE[\Vert\tilde{\beta}_m - \beta_0\Vert_{2, \widehat{\Sigma}}^2] \in [\frac{\var(\epsilon)\lambda_d(\Sigma)}{m(\lambda_1(\Sigma) - C_1\sqrt{\frac{d}{m}})}, \frac{d\var(\epsilon)\lambda_1(\Sigma)}{m(\lambda_d(\Sigma) + C_2\sqrt{\frac{d}{m}})}],
\end{eqnarray*}
and this conclude \eqref{eq:olsBound}.

To show~\eqref{eq:T3} for simplicity let $\hat{r}_0 = \Vert\beta_0\Vert_{2, \widehat{\Sigma}}$ and note that
\begin{eqnarray*}
\EE[(r_n - \hat{r}_0)^2] &=& \EE[(r_n - r_0 + r_0 - \hat{r}_0)^2] \\
&\leq & 2\EE[(r_n - r_0)^2] + 2\EE[(r_0 - \hat{r}_0)^2], 
\end{eqnarray*}
where $r_0^2 = \Vert\beta_0\Vert_{2, \Sigma}^2 = \var(Y) - \var(\epsilon)$. Assuming that $m_2 = \EE[Y^2]$, $m_4 = \EE[Y^4]$ exist we have (Exercise 5.8 in \cite{CasellaBerger})
\begin{eqnarray*}
\EE[(r_n^2 - r_0^2)^2] \leq \frac{1}{n}(m_4 - \frac{n-3}{n-1}m_2^2).
\end{eqnarray*}
Hence 
\begin{eqnarray*}
\EE[(r_n - r_0)^2] &\leq & \frac{1}{r_0^2}\EE[(r_n^2 - r_0^2)^2] \\
& \leq & \frac{1}{nr_0^2}(m_4 - \frac{n-3}{n-1}m_2^2).
\end{eqnarray*}
Since $\Vert\cdot\Vert_{2, \Gamma}$ is a continuous function of $\Gamma$ and since $\Vert\widehat{\Sigma} - \Sigma\Vert = O(n^{-1})$ we have
\begin{eqnarray*}
\EE[(r_0 - \hat{r}_0)^2] &=& \EE[(\Vert\beta_0\Vert_{2, \Sigma} - \Vert\beta_0\Vert_{2, \widehat{\Sigma}})^2] \\
&\lesssim & \frac{1}{n}.
\end{eqnarray*}
Furthermore for large enough $n$ we have $\hat{r}_0/r_0 \geq 1/2$, therefore
\begin{eqnarray*}
T_3 = \EE[\frac{(r_n - \hat{r}_0)^2}{\hat{r}_0^2}]
= \EE[\frac{(r_n - \hat{r}_0)^2}{r_0^2}\frac{r_0^2}{\hat{r}_0^2}] 
\lesssim \frac{1}{n},
\end{eqnarray*}
and this completes the proof of \eqref{eq:T3}.
\hfill
\end{proof}
\subsection{Uniqueness in Example~\ref{normalExp}}
We show that for $Y \stackrel{d}{=} \beta_0^T X + \epsilon$ such that 
\begin{equation*}
X = [X^1, X^2]^T\in\RR^2, X^1\sim\mathcal{N}(0, 1), X^2\sim\text{Exp}(1),
\epsilon\in\mathcal{N}(0, 1),\beta_0^T = [1, 2]
\end{equation*}
we have $\mathcal{B}_0 = \{\beta_0\}$. We prove this by contradiction. 
\begin{proof}
Assume $|\mathcal{B}_0| > 1$. Take $\beta\neq\tilde{\beta}\in\mathcal{B}_0$. Let $\beta = [\beta_1, \beta_2]^T$ and $\tilde{\beta} = [\tilde{\beta}_1, \tilde{\beta}_2]^T$. Note that in this case both $\beta$ and $\tilde{\beta}$ lead to the same distribution for $Y$ and hence using the moment generating function of $Y$ we have 
\begin{eqnarray}\label{MGequality}
\exp(\mu\beta_1 t + \frac{\beta_1^2 t^2}{2})\frac{1}{1 - \beta_2t} = \exp(\mu\tilde{\beta}_1 t + \frac{\tilde{\beta}_1^2 t^2}{2})\frac{1}{1 - \tilde{\beta}_2 t}
\end{eqnarray}
for all $t$ such that $\beta_2 t < 1$ and $\tilde{\beta}_2 t < 1$. First suppose $\beta_2 = 0$. By rewriting \eqref{MGequality}, we have 
\begin{eqnarray*}
\exp(\mu(\beta_1 - \tilde{\beta}_1) t + \frac{(\beta_1^2 - \tilde{\beta}_1^2) t^2}{2}) = \frac{1}{1 - \tilde{\beta}_2 t}
\end{eqnarray*}
for all $t$ such that $\tilde{\beta}_2 t < 1$. This equality holds if and only if $\tilde{\beta}_2 = 0$ and $\beta_1 = \tilde{\beta}_1$. By symmetry, we reach the same conclusion assuming $\tilde{\beta}_2 = 0$. Now assume $\beta_2\neq 0$ and $\tilde{\beta}_2\neq 0$. Without loss of generality assume $\beta_2 < \tilde{\beta}_2$. Then we rewrite \eqref{MGequality} as
\begin{eqnarray}\label{MGequality1}
\exp(\mu(\beta_1 - \tilde{\beta}_1) t + \frac{(\beta_1^2 - \tilde{\beta}_1^2) t^2}{2})\frac{1 - \tilde{\beta}_2 t}{1 - \beta_2 t} = 1
\end{eqnarray}
for all $t$ such that $\beta_2 t < 1$ and $\tilde{\beta}_2 t < 1$. We study this within the following cases:

\textit{Case 1:} $\beta_2 < \tilde{\beta}_2 < 0$. In this case we have $|\tilde{\beta}_2| < |\beta_2|$ and hence the equality in ~\eqref{MGequality1} holds for $t > -1/|\beta_2|$. Now note that 
\begin{eqnarray*}
\lim_{t \searrow -1/|\beta_2|} \exp(\mu(\beta_1 - \tilde{\beta}_1) t + \frac{(\beta_1^2 - \tilde{\beta}_1^2) t^2}{2})\frac{1 - \tilde{\beta}_2 t}{1 - \beta_2 t} = \infty,
\end{eqnarray*}
which is contradiction with \eqref{MGequality1}.

\textit{Case 2:} $\beta_2 < 0 < \tilde{\beta}_2$. In this case \eqref{MGequality1} must hold for $-1/|\beta_2| < t < 1/\tilde{\beta}_2$. Taking the limit $t \searrow -1/|\beta_2|$ again we reach contradiction. 

\textit{Case 3:} $0 < \beta_2 < \tilde{\beta}_2$. In this case \eqref{MGequality1} must hold for $t < 1/\tilde{\beta}_2$. Taking the $t \nearrow 1/\beta_2$ implies 
\begin{eqnarray*}
\lim_{t \nearrow 1/\tilde{\beta}_2|} \exp(\mu(\beta_1 - \tilde{\beta}_1) t + \frac{(\beta_1^2 - \tilde{\beta}_1^2) t^2}{2})\frac{1 - \tilde{\beta}_2 t}{1 - \beta_2 t} = 0,
\end{eqnarray*}
which is again a contradiction. Therefore we must have $\beta_2 = \tilde{\beta}_2$ which implies
\begin{eqnarray*}
\exp(\mu(\beta_1 - \tilde{\beta}_1)t + \frac{(\beta_1^2 - \tilde{\beta}_1^2)t^2}{2}) = 1
\end{eqnarray*}
for all $\beta_2 t < 1$ which means
\begin{eqnarray*}
\mu(\beta_1 - \tilde{\beta}_1)t + \frac{(\beta_1^2 - \tilde{\beta}_1^2)t^2}{2} = 0.
\end{eqnarray*}
Since $\mu \neq 0$ this implies $\beta_1 = \tilde{\beta}_1$.
\end{proof}
\bibliography{unlinkedreg}
\end{document}